%     Ce ficher est du plain TeX.   
%     Schema stable pour Riccati
%     Fran\c{c}ois Dubois, aout 2000, 22, 26 mai, 16 sept  2002.
%     version Linux du 28 mai 2002.
%     nouvelle version fd le 08 sept, 05, 06 octobre 2005.
%     modifs minimes  fd le 4, 9, 13, 18 janvier, 5 mai  2011.
%%%%%%%%%%%%%%%%%%%%%%%%%%%%%%%%%%%%%%%%%%%%%%%%%%%%%%%%%%%%%%%%%% 
%                       thanks to Gabriel Turinici
\input epsfx.tex
%                       thanks to Gabriel Turinici
%%%%%%%%%%%%%%%%%%%%%%%%%%%%%%%%%%%%%%%%%%%%%%%%%%%%%%%%%%%%%%%%%%  
\overfullrule=0pt
%  \nopagenumbers  

%taille d'agrandissement  :
\magnification=1400
\hsize=12.2cm
\vsize=16.0cm
% marge gauche : \hoffset=-1cm
\hoffset=-0.4cm
\voffset=.8cm
\baselineskip 16 true pt

\font \smcaps=cmbx10 at 12 pt

% \font \pbf=cmb10   scaled 200

\font \gcaps=cmbx10 scaled 1333
\font \smcaps=cmcsc10 
\font \pecaps=cmcsc10 at 9 pt

%  Pagination d'apres Raymond S\'eroul pages 231 et 68 
\newtoks \hautpagegauche  \hautpagegauche={\hfil}
\newtoks \hautpagedroite  \hautpagedroite={\hfil}
\newtoks \titregauche     \titregauche={\hfil}
\newtoks \titredroite     \titredroite={\hfil}
\newif \iftoppage         \toppagefalse   
\newif \ifbotpage         \botpagefalse    
\titregauche={\pecaps    Fran\c{c}ois Dubois and Abdelkader  Sa\"{\i}di }
\titredroite={\pecaps    Homographic  scheme  for Riccati equation  }
\hautpagegauche = { \hfill \the \titregauche  \hfill  }
\hautpagedroite = { \hfill \the \titredroite  \hfill  }
\headline={ \vbox  { \line {  
\iftoppage    \ifodd  \pageno \the \hautpagedroite  \else \the
\hautpagegauche \fi \fi }     \bigskip  \bigskip  }}
\footline={ \vbox  {   \bigskip  \bigskip \line {  \ifbotpage  
\hfil {\oldstyle \folio} \hfil  \fi }}}

%  petite boullette

%  nombres mathematiques 
\def\N{{\rm I}\! {\rm N}}
\def\R{{\rm I}\! {\rm R}}

%  indice bas
\def\ib#1{_{_{\scriptstyle{#1}}}}

% les valeurs absolues, les modules et les nornes
\def\abs#1{\mid \! #1 \! \mid }

\def\mod#1{\setbox1=\hbox{\kern 3pt{#1}\kern 3pt}%
\dimen1=\ht1 \advance\dimen1 by 0.1pt \dimen2=\dp1 \advance\dimen2 by 0.1pt
\setbox1=\hbox{\vrule height\dimen1 depth\dimen2\box1\vrule}%
\advance\dimen1 by .1pt \ht1=\dimen1
\advance \dimen2 by .01pt \dp1=\dimen2 \box1 \relax}
\rm

% les valeurs absolues, les modules et les nornes
\def\abs#1{\mid \! #1 \! \mid }

\def\mod#1{\setbox1=\hbox{\kern 3pt{#1}\kern 3pt}%
\dimen1=\ht1 \advance\dimen1 by 0.1pt \dimen2=\dp1 \advance\dimen2 by 0.1pt
\setbox1=\hbox{\vrule height\dimen1 depth\dimen2\box1\vrule}%
\advance\dimen1 by .1pt \ht1=\dimen1
\advance \dimen2 by .01pt \dp1=\dimen2 \box1 \relax}

\def\nor#1{\setbox1=\hbox{\kern 3pt{#1}\kern 3pt}%
\dimen1=\ht1 \advance\dimen1 by 0.1pt \dimen2=\dp1 \advance\dimen2 by 0.1pt
\setbox1=\hbox{\kern 1pt  \vrule \kern 2pt \vrule height\dimen1 depth\dimen2\box1
\vrule
\kern 2pt \vrule \kern 1pt  }%
\advance\dimen1 by .1pt \ht1=\dimen1
\advance \dimen2 by .01pt \dp1=\dimen2 \box1 \relax}

%  carre de fin de dmonstration
\def\sqr#1#2{{\vcenter{\vbox{\hrule height.#2pt \hbox{\vrule width .#2pt height#1pt 
\kern#1pt \vrule width.#2pt} \hrule height.#2pt}}}}
\def\square{\mathchoice\sqr64\sqr64\sqr{4.2}3\sqr33} 

% indice bas
\def\ib#1{_{_{\scriptstyle{\!#1}}}}

%  fleche propre sur les vecteurs

%  Star propre 

%%%%%%%%%%%%%%%%%%%%%%%%%%%%%%%%%%%%%%%%%%%%%%%%%%%%%%%%%%%%%%%%%%%%%%%%%%%%%%%%%%%%%%
%%%%%%%%%%%%%%%%%%%%%%%%%%%%%%%%%%%%%%%%%%%%%%%%%%%%%%%%%%%%%%%%%%%%%%%%%%%%%%%%%%%%%%

$~$ 

 \bigskip 

\centerline{\gcaps  Homographic  scheme  for Riccati equation 
\footnote {$ \, ^{^{\scriptstyle \bf \square}}$}
{\rm Rapport CNAM-IAT n$^{\rm o}$  338-2K, 29 ao\^ut 2000.
\smallskip  \vskip-1pt 
Rapport de recherche    n$^{\rm o}$  2000-32 du laboratoire IRMA de
l'Universit\'e
\smallskip  \vskip-1pt 
Louis Pasteur \`a   Strasbourg.   Edition du 5 mai 2011. }}
 
\bigskip  \bigskip \bigskip

\centerline { \smcaps Fran\c{c}ois Dubois~$^{\rm ab}$    
and    Abdelkader Sa\"{\i}di~$^{\rm c}$  }

\bigskip \bigskip

\centerline { ~$^{\rm a} \,\,$  Conservatoire National des Arts et M\'etiers } 

\centerline {  15 rue Marat, F-78 210, Saint  Cyr l'Ecole, France. }  

\smallskip    

\centerline { ~$^{\rm b} \,\,$  Applications Scientifiques du Calcul Intensif  }  

\centerline {  b\^at. 506, BP 167, F-91~403  Orsay  Cedex, France.}

\smallskip    

\centerline { ~$^{\rm c} \,\,$   Universit\'e Louis Pasteur } 

\centerline { Institut de Recherche  Math\'ematique Avanc\'ee }

\centerline { 7 rue Ren\'e Descartes, F-67 084 Strasbourg, France. }

\smallskip 
\centerline {   dubois@asci.fr, saidi@math.u-strasbg.fr. }

\bigskip  \bigskip

\noindent {\bf Abstract} \smallskip 
In this paper we present a numerical scheme for the resolution of matrix
Riccati equation, usualy used in control problems. The scheme is unconditionnaly stable
and the solution is definite positive at each time step of the resolution. We prove
the convergence in the scalar case and present several numerical experiments for
classical test cases.

\smallskip \noindent    {\bf Keywords}: 
control problems, ordinary differential equations, stability.

\smallskip \noindent    {\bf AMS classification}:
 34H05,     49K15, 65L20, 93C15.  
 
\bigskip  

\noindent {\bf Contents}  

\setbox21=\hbox{$\qquad$ Introduction}
\setbox22=\hbox{$\qquad$ Scalar Riccati equation}
\setbox23=\hbox{$\qquad$ Matrix Riccati equation}
\setbox24=\hbox{$\qquad$ Numerical experiments }
\setbox25=\hbox{$\qquad$ Conclusion }
\setbox26=\hbox{$\qquad$ Acknowledgments }
\setbox27=\hbox{$\qquad$ References. }
\setbox55= \vbox {\halign{#&#\cr
$\qquad \quad$ 1) & \box21 \cr   $\qquad \quad$ 2) & \box22 \cr  
$\qquad \quad$ 3) & \box23 \cr   $\qquad \quad$ 4) & \box24 \cr
$\qquad \quad$ 5) & \box25 \cr   $\qquad \quad$ 6) & \box26 \cr 
$\qquad \quad$ 7) & \box27 \cr  }}
\smallskip \noindent $ \box55 $

~ 

%%%%%%%%%%%%%%%%%%%%%%%%%%%%%%%%%%%%%%%%%%%%%%%%%%%%%%%%%%%%%%%%%%%%%%%%%%%%%%%%%%%%%%
%%%%%%%%%%%%%%%%%%%%%%%%%%%%%%%%%%%%%%%%%%%%%%%%%%%%%%%%%%%%%%%%%%%%%%%%%%%%%%%%%%%%%%
\vfill \eject 
\noindent  {\smcaps 1) $ \quad$ Introduction.}

\bigskip \noindent $\bullet \quad$
We study the optimal control of a differential linear system

\smallskip \noindent  (1.1) $\qquad \displaystyle
{{{\rm d}y}\over{{\rm d}t}} \,\,=\,\, A \, y \,+\, B \, v\,\,, $

\smallskip \noindent 
where the state variable $\, y(\bullet)\, $ belongs to $\,\R^n\, $ and the control 
variable $\,v(\bullet)\, $ belongs to $\,\R^m\, $, with $n$ and $m$ be given integers~:

\smallskip \noindent  (1.2) $\qquad \displaystyle
y(t) \in \R^n \,, \quad v(t) \in \R^m \,.$

\smallskip \noindent 
Matrix $\,A\, $ is composed by $n$ lines and $n$ columns and matrix $\,B\, $ contains 
$n$
lines and $m$ columns. Both matrices $\,A\, $ and $\,B\, $ are fixed relatively to the evolution in
time. Ordinary differential equation (1.1) is associated with an initial condition

\smallskip \noindent  (1.3) $\qquad \displaystyle
y(0) = y_0   $

%%%%%%%%%%%%%%%%%%%%  fd le 08 septembre 2005  %%%%%%%%%%%%%%%%%%%%
\toppagetrue  
\botpagetrue    
%%%%%%%%%%%%%%%%%%%%  fd le 08 septembre 2005  %%%%%%%%%%%%%%%%%%%%

\smallskip \noindent  
with $\,y_0\, $ be given in $\,\R^n $. Morever the solution of system (1.1) (1.3)
is parame\-trized by the function  $\, v(\bullet)\, \, $ and instead of the short notation
$\,y(t) $, we can set more precisely

\smallskip \noindent  (1.4) $\qquad \displaystyle
y(t)\,\, =\,\, y(t\,;\, y_0,v(\bullet)) \,.$

\smallskip \noindent
 The control problem consists of finding the minimum $\,u(t)\, $ of some quadratic 
functional $\,J(\bullet)\, $:

\smallskip \noindent  (1.5) $\qquad \displaystyle
J(u(\bullet)) \,\, \leq \,\, J(v(\bullet)) ,\quad \forall \, v(\bullet) \,.$

\smallskip \noindent
The functional $\,J(\bullet)\, $ depends on the control variable function 
$\,v(\bullet) $, is additive relatively  to the time and represents the coast 
function. We set classically~:

\smallskip \noindent  (1.6) $\,\,\,  \displaystyle
J(v(\bullet))\,\,=\,\, {1 \over 2} \int_{0}^{T} (Q y(t),y(t)) {\rm d}t + {1 \over
2} \int_{0}^{T} (R v(t),v(t)) {\rm d}t +{1 \over 2}(D y(T),y(T)) \,.$

\smallskip \noindent
Functional $\,J(\bullet)\, $ is parametrized by the horizon $\,T>0\, $, the symmetric
semi-definite positives  $n$ by $n$ constant  matrices $\,Q\, $ and $\,D\, $~: 

\smallskip \noindent  (1.7) $\qquad \displaystyle
(Qy,y)\, \geq \, 0 ,\quad \forall \,  y \in \R^n,\quad y\,  \neq \, 0 \, ,$
\smallskip \noindent  (1.8) $\qquad \displaystyle
(Dy,y)\, \geq \, 0 ,\quad \forall \,  y \in \R^n ,\quad y\,  \neq \, 0 \,.$

\smallskip \noindent 
and the symmetric definite positive $\,m$ by $\,m\, $ constant matrix $\,R\, $~:

\smallskip \noindent  (1.9) $\qquad \displaystyle
(Ru,u) > 0 ,\quad \forall \,  u \in \R^m,\quad u\,  \neq \, 0 \,.$

\bigskip \noindent  $\bullet \quad$
Problem (1.1) (1.3) (1.5) (1.6) is a classical mathematical
modeling of linear quadratic loops for dynamical systems in automatics (see {\it e.g.}
Athans, Falb [AF66], Athans, Levine and Levis [ALL67],  Kawakernaak-Sivan [KS72],
Faur\-re Robin [FR84], Lewis [Le86]).  When control function  $\,v(\bullet)\, $  is
supposed to be squarely integrable $\,(v(\bullet)\in L^2(]0,T[,\R^m))\, $ then the
control problem (1.1) (1.3)  (1.5) (1.6) has a unique solution $\,u(t)\, $ (see for
instance Lions [Li68])~:

\smallskip \noindent  (1.10) $\qquad \displaystyle
u(t)\, \in  L^2(]0,T[,\, \R^m) \,.$

\smallskip \noindent 
When there is no constraint on the control variable the minimum of functional
$\,J(\bullet)\, $  is characterized by the condition~:

\smallskip \noindent  (1.11) $\qquad \displaystyle
{\rm d}J(u) \, {\scriptstyle \bullet} \,  w \,\, = \,\,0 \, ,\quad \forall \,  w  \in  
L^2(]0,T[,\, \R^m)\,,$

\smallskip \noindent
which is not obvious to derive because $\,y(\bullet)\, $ is a function  of 
$\,v(\bullet) $. We introduce the  differential equation (1.1) as a constraint
between $\,y(\bullet)\, $ and $\,v(\bullet)\, $ with  the associated Lagrange
multiplayer $\,p\, $. We set~:

\smallskip \noindent  (1.12) $\qquad \displaystyle
{\cal{L}}(y,v;p) \,\, =\,\, J(v)\,-\,\int_{0}^{T} \biggl( p,{{{\rm d}y}\over{{\rm
d}t}}\, \,-\,Ay\,-\,Bv \biggr)  \, {\rm d}t \,$

\smallskip \noindent
and the variation of $\,{\cal{L}}(\bullet)\, $ under an infinitesimal
variation $\,\delta y(\bullet)\, $, $\,\delta v(\bullet)\, 
$ and $\,\delta p(\bullet)\,$ of  the other variables can be conducted as follow~:

%           Calcul des variations 

\setbox20=\hbox{$  \displaystyle \delta {\cal{L}} $}
\setbox21=\hbox{$ \displaystyle   \int_{0}^{T}(Qy,\delta y)\, {\rm d}t \,+\, 
\int_{0}^{T} (Rv,\delta v)\, {\rm d}t \,\,+\,\, (Dy(T),\delta y(T))  $}
\setbox22=\hbox{$ \qquad \displaystyle  \,-\, \,\,  \int_{0}^{T}(\delta p,{{{\rm
d}y}\over{{\rm d}t}} - Ay - Bv) \, {\rm d}t \,\,-\,\, \int_{0}^{T}\bigl(\,p\,,\, 
{{{\rm d}\delta y}\over{{\rm d}t}}\,-\,A\, \delta y\, -\, B\, \delta v)\,{\rm d}t \, $}
\setbox23=\hbox{$ \qquad \displaystyle  \int_{0}^{T}(Qy,\delta y)\,{\rm d}t \,\,+\,\,
\int_{0}^{T} (Rv,\delta v)\,{\rm d}t  \,+\,  \bigl( Dy(T),\delta y(T) \bigr) $}
\setbox24=\hbox{$ \qquad \displaystyle  \,-\,\, 
\int_{0}^{T} \bigl(\delta p,{{{\rm d}y}\over{{\rm d}t}}\,-\,Ay\,-\,Bv \bigr) \,{\rm
d}t\,\,-\,\, \Bigl[ p\, \delta  y \Bigr]_{0}^{T} \,\,+\,\, \int_{0}^{T} \bigl( {{{\rm d}p}
\over {{\rm d}t}},\delta y  \bigr) \, {\rm d}t  \,$}
\setbox25=\hbox{$ \qquad \displaystyle  \,+\,  \int_{0}^{T} \bigl( A^{\rm 
\displaystyle t}
p,\delta y  \bigr) \,  {\rm d}t \,\, +\,\,  \int_{0}^{T} \bigl( B^{\rm \displaystyle t}
p,\delta v  \bigr) {\rm d}t \,, \quad $ and}
\setbox40= \vbox {\halign{#&#&# \cr $\box20$  & = &  $\box21$ \cr & & $\box22$ \cr
 & = & $\box23$ \cr  & & $\box24$ \cr  & & $\box25$ \cr}}
\smallskip \noindent $ \box40 $

\setbox20=\hbox{$\displaystyle  \delta{\cal{L}} \,\,=\,\, \int_{0}^{T} 
\bigl( {{{\rm d}p}  \over {{\rm d}t}} \,+\,Qy\,+\, A^{\rm \displaystyle t}p \,,\, 
\delta y \bigr) \, {\rm d}t \,+\, \int_{0}^{T} \bigl( Rv \,+\,B^{\rm \displaystyle t}
p,\delta v \bigr) \, {\rm d}t   \, $}
\setbox21=\hbox{$\displaystyle  \qquad \quad   -\,\,\,  \int_{0}^{T} 
\bigl( \delta p,{{{\rm d}y}\over{{\rm d}t}}\,-\,Ay\,-\,Bv \bigr) \, {\rm d}t
\,\,+\, \,\bigl(  Dy(T)\,-\,p(T),\delta y(T) \bigr)  \, $}
\setbox23= \vbox {\halign{#\cr \box20 \cr \box21 \cr}}
\setbox24= \hbox{ $\vcenter {\box23} $}
\setbox25=\hbox{\noindent (1.13) $\, \, \left\{  \box24 \right. $}
\smallskip \noindent $ \box25 $

\smallskip \noindent
because $\,\delta y(0) \,\,=\,\, 0\, $ when the initial condition (1.3) is always
satisfied by the function $\,(y \,+\, \delta y)(\bullet)\, $.

\bigskip \noindent  $\bullet \quad$
The research of a minimum for $\,J(\bullet)\, $ (condition (1.11)) can be
rewritten under the form of research of a saddle point for
Lagrangien $\,{\cal{L}}\, $ and we deduce from (1.13) the evolution equation for the 
adjoint variable~:

\smallskip \noindent  (1.14) $\qquad \displaystyle
{{{\rm d}p} \over {{\rm d}t}}\,+\,A^{\rm \displaystyle t} p\,+\,Q\, y\,\,=\,\,0 \,,$

\smallskip \noindent 
the final condition when $\,t\,\,=\,\,T\, $,

\smallskip \noindent  (1.15) $\qquad \displaystyle
p(T) \,\,=\,\, D \, y(T)											\,$

\smallskip \noindent
and the optimal command in terms of the adjoint state $\,p(\bullet)\, $~:

\smallskip \noindent  (1.16) $\qquad \displaystyle
R \, u(t) \,+\, B^{\rm \displaystyle t} \, p(t) \,\,=\,\, 0 \, . $

\smallskip \noindent 
We observe that the differential system (1.1) (1.14) joined with the initial
condition (1.3) and the final condition (1.15) is coupled through
the initial optimality condition (1.16).  In practice, we need a linear feedback
function of the state  variable $\,y(t)\, $ instead of the adjoint variable $\,p(t)\,
$.  Because adjoint state $\,p(\bullet)\, $ linearily depends on state variable
$\,y(\bullet)\, $ we set classically~:

\smallskip \noindent  (1.17) $\qquad \displaystyle
p(T) \,\,=\,\, X(T-t) \,\,   y(t),\quad 0\,  \leq \, t \,  \leq \, T\,,$

\smallskip \noindent  
with a symmetric $n$ by $n$ matrix $\,X(\bullet)\, $, positive definite for $\,t >
0\,\, $ (see {\it e.g.} [KS72] or [Le86]).

\smallskip \noindent  (1.18) $\qquad \displaystyle
X(t)\,$ is a symmetric $\,\, \displaystyle n\times n\,\,$  definite positive
matrix, $ \,\,\, t>0 \,.$

\smallskip \noindent 
The final condition (1.15) is realised for each value $\,y(T)\, $, then we have the
following condition~: 

\smallskip \noindent  (1.19) $\qquad \displaystyle
X(0) \,\,=\,\, D \,, $

\smallskip \noindent 
and introducing the representation (1.17) in the differential equation
(1.14) and (1.1) we obtain~: 

\smallskip \noindent   $ \displaystyle
-{{{\rm d}X} \over {{\rm d}t}}(T-t)\, y(t)\,+\,X(T\,-\,t) \bigl[ \,
Ay(t)\,+\,Bu(t)  \, \bigr] \,\,  +\,\, A^{\rm \displaystyle t}\,
X(T-t)\, y(t)\,+\,Qy(t)\,\,=\,\,0 \,.$

\smallskip \noindent 
We replace the control $\,u(t)\, $ by its value obtained in relation (1.16) 
and we deduce~:

\smallskip \noindent   $ \displaystyle
-{{{\rm d}X} \over {{\rm d}t}}(T-t) \,+\,X(T-t)  \bigl[ \, A\,+\,B(-R^{-1})\,
B^{\rm \displaystyle t}\, X(T-t)   \, \bigr] \,+\,A^{\rm \displaystyle t}\,
X(T-t)\,+\,Qy(t)\,\,=\,\,0 \,.$

\smallskip \noindent
This last equation is realised for each state value $\,y(t)\, $. Replacing $\,t\, $ by
$\,T-t\, $ in this equation, we get~:

\smallskip \noindent  (1.20) $\qquad \displaystyle
{{{\rm d}X} \over {{\rm d}t}}\,\, -\,\, \bigl( \, XA\,+\,A^{\rm \displaystyle t}\,
X \, \bigr) \,\,+\, \,XBR^{-1}BX\,-\,Q\,\,=\,\,0 \,,$

\smallskip \noindent
which defines the Riccati equation associated with the control problem 
(1.1) (1.3) (1.5) (1.6).

\bigskip \noindent  $\bullet \quad$
In this paper we study the numerical approximation of differential system
(1.19) (1.20).  Recall that datum matrices $\,Q\,, \,D\, $ and $\,K$, with $\,K\,$
defined according to~:

\smallskip \noindent  (1.21) $\qquad \displaystyle
K \,\,=\,\, B\, R^{-1}\, B^{\rm \displaystyle t} \,,$

\smallskip \noindent
are $\,n \times n\, $ symmetric matrices, with  $\,Q\, $ and $\,D\, $ semi-definite
positive and $\,K\, $ positive definite ;  datum matrix $\,A\, $ is an $\,n\, $ by $n$
matrix  without any other condition and the unknown matrix $\,X(t)\, $ is symmetric.
We have the following property (see {\it e.g.} [Le86]).

\bigskip 
\bigskip \noindent  {\bf Proposition 1.  Positive definitness of the solution
of Riccati equation.}

\noindent 
Let $\,K, \, Q, \, D, A\, $ be given $\,n \times n\, $ matrices with $\,K,\, Q, \,
D\, $ symmetric matrices,  $\,Q\, $ and  $\,D\, $ positive matrices and $\,K\, $ a
definite positive matrix. Let $\,X(\bullet)\, $ be the solution of the Riccati
differential equation~:
\smallskip \noindent  (1.22) $\qquad \displaystyle
{{{\rm d}X} \over {{\rm d}t}}\,-\, \bigl( \, XA\,+\,A^{\rm \displaystyle t}\, X \, \bigr)
\,\, +\, \,X\,K\, X\,-\,Q\,\,=\,\,0     $
\smallskip \noindent 
with initial condition (1.19). Then $\,X(t)\, $ is well defined for each $\,t \, \geq
\, 0\, $, is symmetric and for each  $\,t > 0\, $, $\,X(t)\, $  is definite positive
and tends to a definite positive matrix $\,X_{\infty}\, $ as t tends to infinity~:
\smallskip \noindent  (1.23) $\qquad \displaystyle
X(t)\, \longrightarrow \, X_{\infty} \quad if \quad t \, \longrightarrow \, \infty
\,.$
\smallskip \noindent 
Matrix $\,X_\infty\, $ is the unique positive symmetric matrix which is solution of the
so-called algebraic Riccati equation~:
\smallskip \noindent  (1.24) $\qquad \displaystyle
-(XA\,+\,A^{\rm \displaystyle t}X)\,\, +\, \,XKX\,-\,Q\,\,=\,\,0\,.$

\bigskip \noindent  $\bullet \quad$
As a consequence of this proposition it is usefull to simplify the feedback command
law  (1.16) by the associated limit command obtained by taking $\,t \longrightarrow
\infty\, $, that is~:
\smallskip \noindent  (1.25) $\qquad \displaystyle
v(t) \,\,=\,\, -R^{-1} \, B^{\rm \displaystyle t} \, X_{\infty} \, y(t) \, ,$

\smallskip \noindent 
and the differential system (1.1) (1.25) is asymptotically stable (see {\it e.g.} [Le86]). 
The pratical computation of matrix $\,X_{\infty}\, $ with direct methods is not
obvious and we refer e.g to [La79] for a description of the state of the art. 
If we wish to compute directly a numerical solution of instationnary Riccati equation
(1.22), classical methods for ordinary differential equations like e.g the forward
Euler method~:

\smallskip \noindent  (1.26) $\qquad \displaystyle
{{1} \over {\Delta t}}(X_{j+1}\,-\,X_j)\,+\,X_j \, K \, X_j\,\, -\,\,(A^{\displaystyle
t} \, X_j\, +\,X_j A)\,-\,Q\,\,=\,\,0 \,,$

\smallskip \noindent 
or Runge Kutta method as we will see in what follows fail to maintain positivity
 of the iterate $\,X_{j+1}\, $ at the order $(j+1)$~:
 
\smallskip \noindent  (1.27) $\qquad \displaystyle
(X_{j+1} \, x\,,\, x) \, > \, 0,\quad \forall \,  x \in \R^n,\quad x\,  \neq \,0 \, , $

\smallskip \noindent   
even if $\,X_j\,$ is positive definite 
if time step $\,\Delta t > 0\, $ is not small enough, see {\it e.g.} [Sa97].
Morever,  the stability constraint (1.27) is not classical and there is at our
knowledge no simple way to determine a priori if time  step $\,\Delta t\, $ is
compatible or not with condition (1.27).

\bigskip \noindent  $\bullet \quad$
In this paper, we propose a method for numerical integration for Riccati equation
(1.22) which maintains condition (1.27) for each time step $\,\Delta
t >0\, $. We detail in second  paragraph the simple case of scalar Riccati equation
and prove the convergence for this particular case ;  under some constraints on
parameters, the scheme is monotonous and remains at the order one of precision,  as
suggested by results of Dieci and Eirola [DE96]. We present  the homographic scheme
in the general case in section $3$ and establish its  principal property ~: for each
time step and without explicit constraint on the time step  $\,\Delta t ,\, $ the
numerical scheme  defines  a symmetric positive definite matrix. We propose and
present four numerical test cases in section~$4$.

%%%%%%%%%%%%%%%%%%%%%%%%%%%%%%%%%%%%%%%%%%%%%%%%%%%%%%%%%%%%%%%%%%%%%%%%%%%%%%%%%%%%%%
%%%%%%%%%%%%%%%%%%%%%%%%%%%%%%%%%%%%%%%%%%%%%%%%%%%%%%%%%%%%%%%%%%%%%%%%%%%%%%%%%%%%%%
\bigskip  %  \bigskip
\noindent  {\smcaps 2) $ \quad$ Scalar Riccati equation.}

\smallskip   \noindent $\bullet \quad$
When the unknown is a scalar variable, we write Riccati equation under the following
form~:

\smallskip \noindent  (2.1) $\qquad \displaystyle
{{{\rm d}x} \over {{\rm d}t}}\,+\,kx^2\,-\,2ax\,-\,q\,\,=\,\,0 \,,$

\smallskip \noindent 
with

\smallskip \noindent  (2.2) $\qquad \displaystyle
k>0,\quad q\, \geq \, 0\,,$

\smallskip \noindent 
and an initial condition~:

\smallskip \noindent  (2.3) $\qquad \displaystyle
x(0)\,\,=\,\,d,\quad d\, \geq \, 0 ,\quad a^2 \,+\, q^2 > 0 \,.$

\smallskip \noindent 
We approach the ordinary differential equation (2.1) with a finite difference
scheme of the type proposed by Baraille [Ba91] for hypersonic chemical cinetics and
independently with the ``family method''  proposed by  Cariolle [Ca79] and studied by
Miellou [Mi84]. We suppose that time step $\,\Delta t\, $ is given strictly positive.
The idea that we have proposed in [Du93], [DS95] and [DS2k] is to write the
approximation  $\,x_{j+1}\, $ at time step $\,(j+1)\Delta t\, $ as a rational fraction
of $\,x_j\, $ with positive coefficients. We decompose first the real number
$\,a\, $ into positive and negative parts~:

\smallskip \noindent  (2.4) $\qquad \displaystyle
a\,\,\, =\, \,\, a^+ \,\,-\,\,  a^- \,,  \qquad  a^+\,\, \geq \, 0 \,, \quad 
\,a^-\, \geq \, 0 \,, \quad a^+\,a^-\,\,=\,\,0\,,$

\smallskip \noindent 
and factorise the product $\,x^2\, $ into the very simple form~:

\smallskip \noindent  (2.5) $\qquad \displaystyle
(x^2)\ib{j+1/2} \,\,=\,\, x\ib{j} \, x\ib{j+1} \,.$

%%%%%%%%%%%%%%%%%%%%%%%%%%%%%%%%%%%%%%%%%%%%%%%%%%%%%%%%%%%
\bigskip \noindent  {\bf Definition 1. $\quad$  Numerical scheme in the scalar case.}

\noindent
For resolution of the scalar differential equation (2.1), we define our numerical 
 scheme by the following relation~: 
\smallskip \noindent  (2.6) $\qquad \displaystyle
{{x\ib{j+1}\,-\,x\ib{j}} \over {\Delta t}} \,+\, k \, x\ib{j} \,
x\ib{j+1}\,-\,2a^+\,x\ib{j}\,\, +\, \,2a^- \,x\ib{j+1}\,\, -\,\, q\,\,=\,\,0\,.$

\bigskip  \noindent   $\bullet \quad$
The scheme (2.6) is implicit because some equation has to be solved in order to
compute  $\,x_{j+1}\, $ whereas $\,x_j\, $ is supposed to be given. In the case of our
scheme this equation is {\bf  linear} and the solution $\,x_{j+1}\, $ is directly
obtained with scalar scheme (2.6) by the homographic relation~:

\smallskip \noindent  (2.7) $\qquad \displaystyle
x_{j+1}\,\,=\,\,{{{\bigl( \, 1\,+\,2a^+\,\Delta t \, \bigr) \, x_j \,\,+\,\, q \,
\Delta t}} \over  {k\,\Delta t\,x_j \,+\,(1\,+\,2a^-\,\Delta t)}} \,.$

\bigskip \bigskip  \noindent  {\bf Proposition 2.  Algebraic properties of the
scalar homographic scheme.}

 \noindent 
Let $\,(x_j) \ib{j \in \N}\, $ be the sequence defined by initial condition~:
\smallskip \noindent  (2.8) $\qquad \displaystyle
x_0\,\,=\,\,x(0)\,\,=\,\,d \,$

\smallskip \noindent 
and recurrence relation (2.7). Then sequence $\,(x_j) \ib{j \in \N}\, $ is
globally defined and remains positive for each time step. 
\smallskip \noindent  (2.9) $\qquad \displaystyle
x_j\, \geq \, 0,\quad \forall \,  j\in \N, \quad \forall \,  \Delta t > 0\,.$

\smallskip \noindent  $\bullet \quad$
If $\,\Delta t\, >\, 0\, $ is chosen such that~:
\smallskip \noindent  (2.10) $\qquad \displaystyle
1\,+\,2|a|\Delta t \,-\, k\, q\, \Delta t^2 \,  \neq \, 0 \,,$
\smallskip \noindent 
then $\,(x_j)\ib{j\in \N}\, $ converges towards the positive solution
$\,x^*\, $ of the  ``algebraic Riccati equation''~:
\smallskip \noindent  (2.11) $\qquad \displaystyle
k \, x^2 \,-\, 2 a \, x \,-\, q \,\,=\,\, 0 \,$
\smallskip \noindent  
and explicitly computed according to the relation 
\smallskip \noindent  (2.12) $\qquad \displaystyle
x^* \,\,=\,\, {{1} \over {k}} \bigl( \, a\,+\,\sqrt{a^2 \,+\, kq}Ê\, \bigr)\,.$
\smallskip  \noindent  $\bullet \quad$
In the exceptional case where $\,\Delta t > 0\, $ is chosen such that
(2.10) is not satisfied, then the sequence $\,(x_j)\ib{j\in \N}\, $ 
is equal to the constant $\displaystyle \,{{1 \, + \, 2 a^{+} \, \Delta t} \over 
{ k \, \Delta t}}\, $  for $\,j \geq  1\, $ and the scheme (2.7) cannot be used for the
approximation of Riccati equation (2.1).

\bigskip \noindent  {\bf Proof of  proposition 2.}
\smallskip \noindent $\bullet \quad $ 
The proof of the relation (2.9) is a consequence of the fact that the recurrence 
relation (2.7) defines an homographic function f~:

\smallskip \noindent  (2.13) $\qquad \displaystyle
x_{j+1}\,\,=\,\,f(x_j) \,$
\smallskip \noindent  (2.14) $\qquad \displaystyle
f(x)\,\,=\,\,{{\alpha x \,+\, \beta} \over {\gamma x \,+\, \delta}} $

\smallskip \noindent 
with positive coefficients $\, \alpha\,, \beta \,, \gamma \,, \delta$~: 
\smallskip \noindent  (2.15) $\qquad \displaystyle
\alpha\,\,=\,\,1\,+\,2a^+\,\Delta t\,,\quad 
\beta\,\,=\,\,q \, \Delta t \,,\quad 
\gamma\,\,=\,\,k \, \Delta t \,,\quad  \delta\,\,=\,\,1\,+\,2a^-\,\Delta t \,.$

\smallskip \noindent 
Because $\,x_0 \,\,=\,\, d \, \geq \, 0\, $, it is then clear that 
$\,x_j \, \geq \, 0\, $ for each $\,j \, \geq \, 0\, $
and relation (2.9) is established. The homographic function $\,f(\bullet)\, $ is
a constant scalar equal to $\displaystyle \,{{\alpha}\over{\gamma}} \,\,=\,\,
{{\beta}\over{\delta}}\, $  when the determinant~:

\smallskip \noindent  (2.16) $\qquad \displaystyle
\Delta \,\,=\,\, \alpha \delta \,-\, \beta \gamma $

\smallskip \noindent 
is equal to zero. It is the case when condition (2.10) is not realised. When (2.10)
holds, we have simply~:

\smallskip \noindent  (2.17) $\qquad \displaystyle
{{f'(x)} \over {f(x)}} \,\,=\,\, {{\alpha \, \delta \,\,-\,\, \beta \,\gamma} \over 
{(\gamma  \,x \,\,+\,\, \delta)\, (\alpha \, x \,+\, \beta)}} $

\smallskip \noindent 
and $\,f\, $ is a bounded monotonic function on the interval $\,]0,\infty[\, $. 
Let $\,x^*\, $ be the positive solution of equation~:

\smallskip \noindent  (2.18) $\qquad \displaystyle
f(x) \,\,=\,\, x\,.$

\smallskip \noindent 
It is elementary to observe that $\,x^*\, $ is given by the relation (2.12). Let
$\,x_-\, $ be defined by~:

\smallskip \noindent  (2.19) $\qquad \displaystyle
x_-\,\,=\,\,-{{ q} \over {k \, x^*}} \,\, = \,\, x^* \,-\, 2 \, {{ \sqrt{a^2 \,+\,
kq}}\over{k}} Ê\, \, ,   $

\smallskip \noindent 
 the other root of equation (2.18). We set~:

\smallskip \noindent  (2.20) $\qquad \displaystyle
u_j\,\,=\,\,{{x^*\,-\,x_j} \over {x_j\,-\,x_-}}\,.$

\smallskip \noindent 
Then we have~:
% 	                                           Calcul de u j+1  en fonction de u j

\smallskip \noindent  (2.21) $\qquad \displaystyle
u_{j+1} \,\,=\,\, {{\displaystyle 
{{\alpha \,x^*\,+\,\beta} \over  {\gamma \, x^*\,+\,\delta}} \,\,-\,\, 
{{\alpha \, x_j\,+\,\beta} \over {\gamma \, x_j\,+\,\delta}}} \over  
{\displaystyle {{\alpha \, x_j\,+\,\beta} \over {\gamma \,x_j\,+\,\delta}} \,\,-\,\, 
{{\alpha \, x_-\,+\, \beta} \over {\gamma \,x_-\,+\, \delta}}}}
\,\,\,=\,\,\, 
{{\gamma \, x_-\,+\, \delta} \over {\gamma \,x^*\ +\, \delta}} \,\, 
{{x_j \,-\, x^*} \over {x_- \,-\, x_j}} \,\,, $

\smallskip \noindent 
if (2.10) holds. Replacing $\,x_-$ by its value given from (2.19) and using
(2.21),  we obtain~:
 
\smallskip \noindent  (2.22) $\qquad \displaystyle
u_{j+1}\,\,=\,\,{{\delta \, x^* - \beta} \over {\alpha \, x^* \,+\, \beta}} \,\, u_j
\,.$

\smallskip \noindent 
The sequence $\,(u_j)\ib{j\in \N}\, $ is geometric and the ratio  $\,\displaystyle
{{u_{j+1}} \over {u_j}} \, $ has always an absolute value less than $1$. Effectively
we have on one hand~:

\smallskip \noindent  $ \qquad  \qquad  \displaystyle
- \bigl( \, \alpha \, x^*\,+\,\beta \, \bigr) \,\,\, 
\leq \,\,\, \delta \, x^* - \beta \, $

\smallskip \noindent 
 because $\quad \,\, \,\alpha \,+\, \delta
\,\,=\,\, 2\,(1\,+\,|a|\Delta t)\,\,\, $ is strictly positive and $\,x^*\, $ is
positive and on the other hand~: 

\smallskip \noindent  $ \qquad  \qquad  \displaystyle
\delta \, x^* - \beta \,  \,\,\leq \,\,\,
\alpha x^* \,+\, \beta  \,  $ 

\smallskip \noindent 
because $\qquad\,\, (\alpha \,-\, \delta) \, x^*
\,+\, 2 \beta \,\,=\,\, 2\Delta t \, (ax^*\,+\,q) \,\,=\,\, \Delta t \, ( k (x^*)^2
\,+\, q )\, \,\,$ which is a positive number. The absolute value of $\,\displaystyle
{{u_{j+1}} \over {u_j}}\, $ is not exactly equal to $1$ because  $\,x^* > 0\, $ and
$\,a \, x^* \,+\, q \,\,=\,\, x^* \, (k \, x^* \,-\, a) \,\,=\,\, x^*
\, \sqrt{a^2\,+\,k\,q} \,\,>\,\, 0\,\,$  according to  (2.10). Then  $\,u_j\, $ is
converging to zero and $\,x_j\, $ toward $\,x^*\, $ that completes the proof.   $
\hfill \square \kern0.1mm $

\bigskip 
\bigskip \noindent  {\bf Theorem 1. $\quad$  Convergence of the scalar scheme. }

 \noindent 
We suppose that the data $\,k,a,q\, $ of Riccati equation satisfy (2.2) and
(2.10) and  that the datum $\,d\, $ associated with the initial condition
(2.3) is relatively closed to $\,x^*\, $, ie~:
\smallskip \noindent  (2.23) $\qquad \displaystyle
-{{1} \over {k\, \tau}} \,+\, \eta \, \,\,  \leq \, \, \,  d - x^* \, \, \,  \leq \,
\,\, C \,,$
\smallskip \noindent 
where $\,C\, $ is some given strictly positive constant $\,(C > 0)\, $, $\,x^*\, $ 
calculated according to relation (2.12) is the limit in time of the Riccati equation,
$\,\tau\, $ is defined from data $\,k, \, a, \, q\, $ according to~:
\smallskip \noindent  (2.24) $\qquad \displaystyle
\tau\,\,=\,\,{{1} \over {2\,\sqrt{a^2\,+\,kq}}} \,,$       
\smallskip \noindent 
and $\,\eta\, $ is some constant chosen such that 
\smallskip \noindent  (2.25) $\qquad \displaystyle
0 \, < \, \eta \, < \, {{1} \over {k\,\tau}} \,. $
\smallskip  \noindent  $\bullet \quad$
We denote by  $\,x(t;d)\, $  the solution of differential equation (2.1)
with initial condition (2.3). Let $\,(x\ib{j}(\Delta t\,;\,d\ib{\Delta}))\ib{(j\in
\N)}\, $ be the solution of the  numerical scheme defined at the relation (2.7) and
let $\,\, d\ib{\Delta} \,\,$ be the initial condition~: \smallskip \noindent  (2.26)
$\qquad \displaystyle x_0(\Delta t\,;\,d\ib{\Delta})\,\,=\,\,d\ib{\Delta} \, .$   
\smallskip \noindent 
We suppose that the numerical initial condition $\,\,  d\ib{\Delta} > 0 \,\,  $
satisfies a  condition analogous to (2.23)~:

\smallskip \noindent  (2.27) $\qquad \displaystyle
- {{1} \over {k\,\tau}}\,+\,\eta \, \,\,  \leq \,\,\, d\ib{\Delta} \, - \, x^* \, \, 
\, \leq \,\,\, C \,,$

\smallskip \noindent 
with $\,C\, $ and $\,\eta > 0\, $ equal to the constant introduced in (2.23) 
and satisfying (2.25) . 
 
\smallskip \noindent  $\bullet \quad$
Then the  approximated value  $\,(x_j(\Delta t \,;\, d\ib{\Delta}))\ib{j \in \N}\, $ is
arbitrarily closed to the exact value  $\,x(j \Delta t \,;  \, d)\, $ for each $\,j\,
$  as  $\,\Delta t \longrightarrow 0 \, $ and $\,d\ib{\Delta} \longrightarrow d\, $. 
More precisely, if $\,a \,  \neq \, 0\, $ we have the  following estimate for the
error at time equal to $\,j \Delta t\, $~:
\smallskip \noindent  (2.28) $\, \displaystyle
\abs{ x(j \Delta t\,;\, d) \,-\, x_j(\Delta t \,;\, d\ib{\Delta})}  \leq 
A\,(\Delta t \,+\, \abs{ d \,-\, d\ib{\Delta}}),\,\,\, \forall \,  j \in \N \, , \,  0
\,<\,  \Delta t \,  \leq \,  B $

\smallskip \noindent 
with constants $\,A > 0, B > 0\, $, depending on data $\,k, a, q, \eta\, $ but
independent on time step $\,\Delta t > 0\, $ and iteration $\,j\, $. 

\smallskip \noindent  $\bullet \quad$
If $\,a \,\,=\,\, 0\, $, the scheme is second order accurate in the following sense~:
\smallskip \noindent  (2.29) $\,\,  \displaystyle
|x(j \Delta t;d) \,-\, x\ib{j}(\Delta t;d\ib{\Delta})|   \leq   A\,(\Delta t^2 \,+
\, \abs{d \,-\, d\ib{\Delta}}),\,\,\, \forall \, j \,  \in \N \, , \,  0 < \Delta t \, 
\leq \, B      \,$
\smallskip \noindent 
with constants $\,A\, $ et $\,B\, $ independent on time step $\,\Delta t\, $ and 
iteration $j$.

%  	                           Remarque

\bigskip \noindent  {\bf Remark. } 
\smallskip \noindent   $\bullet \quad $ 
A direct application of the Lax theorem (see {\it e.g.} [La74]) for numerical scheme
associated to   ordinary differential equations is not straightforward because both
Riccati equation and the numerical scheme are nonlinear. Our proof is based on the
fact that this problem is algebrically relatively simple.

%% \bigskip  
\bigskip \noindent  {\bf Lemma 1. $\quad$  }

 \noindent 
Let $\, d\ib{\Delta} \,$ be some discrete initial condition and
$\, x( t \,;\, d\ib{\Delta})\, $ be the exact solution of Riccati differential
equation (2.1) associated to initial condition (2.26). We set 

\smallskip \noindent  (2.30) $\quad \displaystyle
y(t\,;\,d\ib{\Delta})\,\,\,=\,\,\,x(t\,;\,d\ib{\Delta})\,-\,x^*  \,. $

\smallskip \noindent 
Then we have 
\smallskip \noindent  (2.31) $\quad \displaystyle
y(t\,;\,d\ib{\Delta})\,\,\, =\,\,\, {{(d\ib{\Delta}\,-\,x^*) \,\, e^{\displaystyle
-{t/\tau}}} \over {1\,+\,k\, \tau \, (d\ib{\Delta}\,-\,x^*)\, (1\,-\,e^{\displaystyle
-{t/\tau}})}}  \,.$		

\bigskip \noindent  {\bf Proof of Lemma 1.}
\smallskip \noindent $\bullet \quad $ 
The real function $\,\, \R \ni t \longmapsto y(t\,;\,d\ib{\Delta}) \in \R \, \, $
introduced in relation (2.30) satisfies the equation  $ \quad \displaystyle
{{{\rm d}y}\over{{\rm d}t}}  \,+\,ky^2\,+\,{{1} \over {\tau}} \, y\,\,=\,\,0
\,, \quad $						 that can be rewritten under the form~:

\smallskip \noindent  $ \displaystyle
{{{\rm d}y} \over {k\,y^2\,+\, {\displaystyle {{y} \over {\tau}}}}}
\,\,\equiv\,\,\tau \, \Bigl( {{{\rm d}y} \over {y}} \,\,-\, \,
{{k\, \tau} \over {1\,+\,k\, \tau \,  y}} \, {\rm d}y \, \Bigr)
\quad = \quad -{\rm d}t \,. $

\smallskip \noindent 
After integration, we have~:

\smallskip \noindent  $ \displaystyle
{{y(t\,;\,d\ib{\Delta})} \over {1\,+\,k\, \tau \, y(t,d\ib{\Delta})}} \,\,\,=\,\,\,
{{y(0\,;\,d\ib{\Delta})} \over  {1\,+\,k\, \tau \, y(0,d\ib{\Delta})}}\,
e^{\displaystyle -{t/\tau}} \,,$

\smallskip \noindent 
giving simply~:

\smallskip \noindent $ \displaystyle
y(t\,;\,d\ib{\Delta})\,\,\, =\,\,\, {{(d\ib{\Delta}\,-\,x^*) \,\, e^{\displaystyle
-{t/\tau}}} \over {1\,+\,k\, \tau \, (d\ib{\Delta}\,-\,x^*)\, (1\,-\,e^{\displaystyle
-{t/\tau}})}}  \,$				

\smallskip \noindent 
{ \it i.e.}  relation (2.31). Then Lemma 1 is established.  $ \hfill \square \kern0.1mm $

% \vfill \eject  
\bigskip  \bigskip \noindent  {\bf Lemma 2. $\quad$  }

  \noindent 
Let $\, d\ib{\Delta} \,$ be some discrete initial condition associated to (2.26) and
$\, \epsilon\ib{j} \,$  the error between the exact solution $\, x(j \Delta t
\,;\, d\ib{\Delta})\, $ of the differential equation  and the solution of the
 numerical scheme $\, x\ib{j}(\Delta t \,;\, d\ib{\Delta})\, $~: 
\smallskip \noindent  (2.32) $\quad \displaystyle
\epsilon\ib{j}  \,\,=\,\,  x(j \Delta t\,;\,d\ib{\Delta}) \,\, -\, \, x_j(\Delta
t\,;\,d\ib{\Delta})  \,. $
\smallskip \noindent
Let $\, x^* \,$ and $\, x_{-} \,$ defined in (2.12) and (2.19) be the two   fixed
points of the homographic function $\, f({\scriptstyle \bullet}) \,$ introduced in
(2.14) and (2.15). We set 
\smallskip \noindent  (2.33) $\quad \displaystyle
h(\xi) \,\,=\,\, {{x^{*} \,+\, x_{-} \, \xi} \over {1 \,+\, \xi}}  \,, \quad \xi  > 
\,-1 \,\, $
\smallskip \noindent 
and we introduce on one hand  the function $\, t \longmapsto z(t) \,$ relative to  the
exact solution~: 
\smallskip \noindent  (2.34) $\quad \displaystyle
z(t) \,\,=\,\, {{x(t\,;\,d\ib{\Delta}) \,-\, x^*} \over {x_- \,-\, x(t,d\ib{\Delta})}}
\,\,$
\smallskip \noindent 
and on the other hand a new sequence $\, u_j  \,$ by the same algebraic relation~: 
\smallskip \noindent  (2.35) $\quad \displaystyle
u_j\,\,=\,\,{{x\ib{j}(\Delta t \,;\, d\ib{\Delta})  \,-\,  x^*} \over {x_-
\,-\, x\ib{j}(\Delta t \,;\, d\ib{\Delta})  }}\,\, .$
\smallskip \noindent 
Then we have the following estimate~: 
\smallskip \noindent  (2.36) $\quad \displaystyle
\abs{ \epsilon\ib{j} }  \, \, \,  \leq \, \,  \,  \abs { h' \bigl( z(0) \bigr)  }
\,\, \,  \abs {  u_j \,-\, z(j \Delta t)}  \,. $

\bigskip \bigskip  \noindent  {\bf Proof of Lemma 2.}
\smallskip \noindent $\bullet \quad $ 
We have constructed function $\, h ({\scriptstyle \bullet})\, $ in order to have   $
\,\,  \displaystyle h \bigl( z(\theta) \bigr) \,=\, z \bigl( h(\theta) \bigr)
\,,\,\, $ $\forall \, \theta \in \R \,. \, $ Then we have 

\smallskip \noindent  $ \displaystyle
x(t\,;\,d\ib{\Delta}) \,\,\,\quad =\,\, h \bigl( z(t) \bigr) \,\,\,\,\quad 
{\rm for \, each} \,\, t > 0$
\smallskip \noindent  $ \displaystyle
 x_j(\Delta t\,,\, d\ib{\Delta}) \,\,=\,\, h(u_j) \,\, \qquad {\rm  for \, each } \,\,
j \, \geq \, 0 \,,  $

\smallskip \noindent 
with function $\, z({\scriptstyle \bullet})\, $ introduced at relation (2.34) and
sequence  $\,u_j\, $  in (2.35). Then $\,\epsilon\ib{j} \, $ can be
rewritten with the help of this function$\,h({\scriptstyle \bullet})\, $ and we have~:

\smallskip \noindent  (2.37) $\quad \displaystyle
\epsilon\ib{j}   \,\, = \,\,    h \bigl( z(j  \Delta t) \bigr) \,-\, h(u_j)  \, 
  \leq \,    \sup_{\displaystyle \xi \in [z(j \Delta t) \,,\,u_j ]}
\,    \abs { h' \bigl( z(\xi) \bigr)  }  \,\,  \abs { z(j \Delta t) \,-\,u_j } . $

\bigskip \bigskip  \noindent  $\bullet \quad$
We note that due to (2.22), the sequence  $\,u_j\, $ is a geometric converging
one, then we have~: $ \displaystyle \quad 
- \abs {u_0}  \,\, \leq \,\, u_j  \,\, \leq \,\,  \abs {  u_0 } \quad  $ for each $\,j
\, \geq \, 0\, . \, $ Moreover thanks to relation (2.31) of Lemma 1,  we have the
following calculus~: 

% \vfill \eject   %%%%%%%%%%%%%%%%%       ajout 4 janvier 2011  fd 

\smallskip \noindent  $ \displaystyle
z(t) \,\,= \,\, {{y(t\,;\,d\ib{\Delta})}\over{x_- - x^* - y(t\,;\,d\ib{\Delta})}} \,$

\smallskip \noindent  $ \displaystyle \qquad \,\,\, 
= \,\, {{  {{\displaystyle (d\ib{\Delta}\,-\,x^*) \,\, e^{\displaystyle
-{t/\tau}}} \over {\displaystyle 1\,+\,k\, \tau \, (d\ib{\Delta}\,-\,x^*)\,
(1\,-\,e^{\displaystyle -{t/\tau}})}}  } \over{  \displaystyle -\, {{1}\over{k \,\tau}}
\,-\, {{(d\ib{\Delta}\,-\,x^*) \,\, e^{\displaystyle -{t/\tau}}} \over {1\,+\,k\, \tau
\, (d\ib{\Delta}\,-\,x^*)\, (1\,-\,e^{\displaystyle -{t/\tau}})}}   }} \,$

\smallskip \noindent  $ \displaystyle \qquad \,\,\, 
= \,\,   \,-\, {{k \, \tau \, (d\ib{\Delta} \,-\, x^*) } \over {1 \,+\, k \, \tau
\, (d\ib{\Delta} \,-\, x^*)}} \, e^{\displaystyle -{t/\tau}} \, \,, $

\smallskip \noindent
thus satisfy clearly~:  $ \displaystyle \quad 
-\,\abs{z(0)} \,\,   \leq \,\,  z(j \Delta t)  \,\, \leq \,\,  \abs {z(0)} \quad \,,
\quad \forall \,  j \, \geq \, 0 \,.$ We remark also that~:

\smallskip \noindent  $ \displaystyle
-1 \,\, < \,\,   z(0) \,\,=\,\, u_0 \,\,=\,\, \,- \, {{k \, \tau \, (d\ib{\Delta} \,-\,
x^*)} \over {1 \,+\, k \, \tau \, (d\ib{\Delta} \,-\, x^*)}} \, .$

\bigskip \noindent $\bullet \quad $ 
Morever, $\, |h'(\xi)|\, $ is a decreasing function for $\,\xi > -1.$  Then by the
two points finite difference Taylor formula and the above statements, we deduce from
(2.37) the estimate (2.36) and  Lemma 2 is established.  $ \hfill \square \kern0.1mm $

\bigskip  \bigskip \noindent  {\bf Lemma 3. $\quad$  }

  \noindent 
Let $\, \tau \,$ be defined in (2.24). We introduce $\, \alpha\,$ and $\, \beta \,$
according to the following relations 
\smallskip \noindent  (2.38) $\quad \displaystyle
\alpha  \,\,=\,\,  \displaystyle  {{1} \over {2 \tau}} \,-\, \abs{a} $
\smallskip \noindent  (2.39) $\quad \displaystyle
\beta \,\,=\,\,\displaystyle {{1} \over {2\tau}} \,+\,  \abs{a} \,$
\smallskip \noindent
and when $\, \Delta t \,> 0 \,$ we define function $\, \varphi( \Delta t ) \,$
according to 
\smallskip \noindent  (2.40) $\quad \displaystyle
\varphi(\Delta t) \,\,=\,\, 1\,+\, {{\tau}  \over {\Delta t}} \, {\rm log} \, \Bigl(
\, {{1 \,-\, \alpha \, \Delta t \, } \over {1 \,+\, \beta \, \Delta t \, }} \, \Bigr)
\,.\, $
\smallskip \noindent
Then with notations introduced in Lemma 1 and the following one 
\smallskip \noindent  (2.41) $\quad \displaystyle
\theta_j \,\, = \,\, {{j \Delta t}\over{\tau}} \,,\, $
\smallskip \noindent
we have 
\smallskip \noindent  (2.42) $\quad \displaystyle
u_j \,-\, z(j \Delta t) \,\,=\,\, - z(0) \, \, e^{\displaystyle \,- \theta_j}  \,
\biggl[
\, 1 \,-\, {\rm exp}  \Bigl( \theta_j  \,  \varphi(\Delta t) \, \Bigr)  \,
\biggr] \,.$

\bigskip \noindent  {\bf Proof of Lemma 3.}
\smallskip \noindent $\bullet \quad $ 
We  study the difference that majorate the right hand side of (2.37). From
relation (2.21), we know that sequence $\, u_j\,$ is geometric and more precisely~:

\smallskip \noindent  $ \displaystyle
u_j \,\, = \,\, \Bigl(   {{\gamma \, x_- \,+\, \delta} \over {\gamma \, x^*
\,+\, \delta}}  \Bigr)^j \, z(0) \,$ 

\smallskip \noindent
with $\, \gamma \,$ and $\, \delta \,$ computed in (2.15)  $ \displaystyle \quad 
\gamma \,\,=\,\, k \, \Delta t \,,  \quad  \delta\,\,=\,\,1\,+\,2a^-\,\Delta t \,.$
We have from equation (2.11) and defintion (2.24) of variable $\, \tau \,$~: 

\setbox20=\hbox{$   \displaystyle  {{\gamma \, x_- \,+\, \delta} \over
{\gamma \, x^* \,+\, \delta}} $}
\setbox29=\hbox{$ \displaystyle   {{ k\, \Delta t \,\, {{ \displaystyle 1}\over{
\displaystyle k}} \,\, \Bigl(  a \,-\, \displaystyle  {{1} \over {2 \tau}} \, \Bigr)  
\,\, +\,\,  1 \,+\, 2 a^- \Delta t} \over {  k\, \Delta t \,\,  \,\, {{
\displaystyle 1}\over{ \displaystyle k}} \,\, \Bigl(  a \,+\, \displaystyle  {{1}
\over {2 \tau}} \, \Bigr)   \,\,  \,+\, 1 \,+\, 2 a^-\, \Delta t}} $}
\setbox21=\hbox{$ \displaystyle   {{ \Delta t \, \Bigl( a^+ \, - \, a^- \,-\,
\displaystyle  {{1} \over {2 \tau}} \Bigr) \,\, +\,\,  1 \,+\, 2 a^- \Delta t} \over {
\Delta t \,   \Bigl(  a^+ \, - \, a^-  \,+\,  \displaystyle  {{1} \over {2 \tau}}
\Bigr ) \,+\, 1 \,+\, 2 a^-\, \Delta t}} $}
\setbox22=\hbox{$   \displaystyle   {{1 \,+\, \Delta t \, \bigl(  \abs{a} \,-\,
 \displaystyle  {{1} \over {2\tau}}  \bigr)  } \over {1 \,+\, \Delta t \, \bigl( 
\abs{a}  \,+\, \displaystyle  {{1} \over {2\tau}} \bigr) }} $}
\setbox23=\hbox{$  \displaystyle  {{1 \,-\, \alpha \, \Delta t } \over {1 \,+\,
\beta \, \Delta t }} \,\, . $}
\setbox40= \vbox {\halign{#&#&# \cr $\box20$  & = &  $\box29$ \cr  & = &  $\box21$
\cr  & = &  $\box22$ \cr  & = &  $\box23$ \cr  }}
\smallskip \noindent $ \box40 $

\setbox20=\hbox{$   \displaystyle  {{\gamma \, x_- \,+\, \delta} \over
{\gamma \, x^* \,+\, \delta}} $}
\setbox23=\hbox{$  \displaystyle  {{1 \,-\, \alpha \, \Delta t } \over {1 \,+\,
\beta \, \Delta t }} \,\, . $}
\setbox40= \vbox {\halign{#&#&# \cr $\box20$  & = &  $\box23$ \cr  }}
% \smallskip \noindent $ \box40 $

\smallskip \noindent 
We can now  write~:

\smallskip \noindent  $ \displaystyle
u_j \,-\, z(j \Delta t) \,\,=\,\, \,-z(0) \, \biggl[ \,  e^{\displaystyle {-j\Delta
t/\tau}} \, \,-\,\, \Bigl(  {{1 \,-\, \alpha \, \Delta t } \over {1 \,+\,
\beta \, \Delta t }}  \Bigr)^j \, \biggr] \,$
\smallskip \noindent  $ \displaystyle \qquad \qquad \qquad \,
=\,\, \,-z(0) \,\,  e^{\displaystyle - j \Delta t / \tau } \,  \biggl[ \, 1 \,- \, 
{\rm exp} \Bigl( \theta_j \,+\, j \, {\rm log}  \Bigl( {{1 \,-\, \alpha \, \Delta t \,
} \over {1 \,+\, \beta \, \Delta t \, }} \, \Bigr)  \, \Bigr)  \, \biggr] \,$
\smallskip \noindent  $ \displaystyle \qquad \qquad \qquad \,
=\,\, \,-z(0) \,\,  e^{\displaystyle - \theta_j } \,  \biggl[ \, 1 \,- \, 
{\rm exp} \Bigl( \theta_j \, \Bigl( 1 \,+Ê{{\tau}\over{\Delta t}}  {\rm log}  \Bigl(
{{1 \,-\, \alpha \, \Delta t \, } \over {1 \,+\, \beta \, \Delta t \, }} \, \Bigr)  \,
\Bigr)  \,\Bigr)\, \biggr] \,.$

\smallskip \noindent 
and relation (2.42) is proven.   Lemma 3 is established.  $ \hfill \square \kern0.1mm $

\bigskip  \bigskip  % \vfill \eject 
\noindent  {\bf Lemma 4. $\quad$  }

 \noindent 
Let $\, \varphi({\scriptstyle \bullet})\,$ be defined by relation (2.40). We suppose
that time step $\, \Delta t \,$ satisfies the condition 

\smallskip \noindent  (2.43) $\quad \displaystyle
0 \,\, < \,\, \Delta t \,\, \leq \,\,{{\tau}\over{2}} \,. \, $

\smallskip \noindent
Then we have 

\setbox10=\hbox{$  \displaystyle  \! \varphi (\Delta t) \, -\, \tau \, \Bigl[  
{{\beta^2 \,-\, \alpha^2} \over {2}} \Delta t \,  - {{1} \over {3}}(\alpha^3
+ \beta^3) \Delta t^2 \, \Bigr] \!  $}

\smallskip \noindent  (2.44) $\quad \displaystyle
\mod  { \box10} \,\,  \leq \,\, {{1} \over {2}}(\alpha^4 + \beta^4)\,\,  \tau  \, 
\Delta t^3  \,\,.$

\bigskip \noindent  {\bf Proof of Lemma 4.}
\smallskip \noindent $\bullet \quad $ 
We have the following elementary  calculus~: 
\smallskip \noindent $ \displaystyle 1  \,-\, \xi \,+\, \xi^2  \,-\, \xi^3 \,\, \leq
\,\, {{1}\over{1+\xi}} \,\, \leq \,\,  1 \,-\,  \xi \,+\, \xi^2  \,-\, \xi^3 \,+\, \xi
^4 \qquad {\rm if} \,\, \abs{\xi} \, < \,1 \,$ 

\smallskip \noindent
and we deduce after integration~: 

\smallskip \noindent $ \displaystyle 
- {{\xi^4}\over{4}} \,\, \leq \,\, \Biggl[ \,  {\rm log}  ( 1 \,+\, \xi) \,-\,
\biggl( \xi \,-\, {{\xi^2} \over {2}} \,+\, {{\xi^3} \over {3}} \biggr) \, \Biggr]
\,\, \leq \,\, - {{\xi^4}\over{4}} \,+\, {{\xi^5}\over{5}} \,\, \leq \,\,
{{\xi^4}\over{2}} \qquad {\rm if} \,\, \abs{\xi} \, < \,{1\over2} \,. \,$ 

\smallskip \noindent
Then we have 

\setbox10=\hbox{$  \displaystyle {\rm log}  ( 1 \,+\, \xi) \,-\, \biggl( \xi \,-\,
{{\xi^2} \over {2}} \,+\, {{\xi^3} \over {3}} \biggr)  $}

\smallskip \noindent  (2.45) $\quad \displaystyle
\mod  { \box10}  \,\,  \leq \, \,\, {{1} \over {2}} \,  |\xi|^4   \qquad  {\rm if} \,
\, |\xi| \,  \leq \, {{1} \over {2}} $

\smallskip \noindent
and we use this estimation with  $\, \,\xi  \,=\, - \alpha  \Delta t \, $  
and   $\,\, \xi  \,=\, \beta  \Delta t .\,  $ We remark that $\, \displaystyle
\abs{a} \, \leq \, {{1}\over{\displaystyle 2\,  \tau}} \,$ then $\,\displaystyle 
\beta \, \leq \, {{1}\over{\displaystyle   \tau}} \,$ and we deduce that  when
condition (2.43) is satisfied, we have 

\smallskip \noindent$ \displaystyle
\Delta t \, \, \leq \, \,{{\tau}\over{2}} \, \, \leq \, \,  {{1}\over{2 \, \beta}}
\,\, \leq \,\, {{1}\over{2 \, \alpha}} \,,\,\,\, $ then   $ \displaystyle \quad 
\alpha \, \Delta t \,  \leq \, \displaystyle {{1} \over {2}}\,\,\,  $ and
 $\, \,\beta \, \Delta t \,  \leq \, \displaystyle {{1} \over {2}}\,. \,\, $

\smallskip \noindent 
We deduce from  (2.45) 

\setbox10=\hbox{$  \displaystyle  {\rm log}  \Bigl( {{1 \,-\, \alpha \, \Delta t \, }
\over {1 \,+\, \beta \, \Delta t \, }} \, \Bigr) \,-\, \Bigl[\, (- \alpha \,-\,
\beta)\Delta t \,+\, {{\beta^2 \,-\, \alpha^2} \over {2}}\Delta t^2 \,+\, 
 {{1} \over {3}}(- \alpha^3 \,-\,
\beta^3) \Delta t^3    \Bigr] \, $}
\smallskip \noindent $  \displaystyle
\mod  { \box10} \,\, \leq $ 
\smallskip \noindent $ \qquad  \qquad \qquad  \qquad \qquad   \qquad   \qquad  
\qquad  \qquad \qquad   \qquad  \displaystyle
\leq \,\, {{1} \over {2}}(\alpha^4 \,+\, \beta^4) \, \Delta t^4 \,.$

\bigskip \noindent  $\bullet \quad$
As a  consequence of the previous inequality, multiplying
the above expression by $\,\displaystyle {{\tau} \over {\Delta t}}$,  we obtain~:

% \vfill \eject 

\setbox10=\hbox{$  \displaystyle  \varphi (\Delta t) \, -\, 1
\,-\, {{\tau} \over {\Delta t}} \, \Bigl[  (- \alpha \,-\, \beta)\Delta t \,+\, 
{{\beta^2 \,-\, \alpha^2} \over {2}}\Delta t^2 \,+\,  {{1} \over {3}}(- \alpha^3 \,-\,
\beta^3) \Delta t^3 \, \Bigr] \, $}
\smallskip \noindent $   \displaystyle
\mod  { \box10} \,\, \leq $ 
\smallskip \noindent $ \qquad  \qquad \qquad  \qquad \qquad   \qquad   \qquad  
\qquad  \qquad \qquad   \qquad  \displaystyle
\leq \,\, \, {{1} \over {2}}(\alpha^4 \,+\, \beta^4) \, \tau \, \Delta t^3 \,, $

\smallskip \noindent 
and due to the fact that $\,\, \, (- \alpha \,-\, \beta) \, \tau \,=\, -1 \,$  we
have  since 
$\,\, \displaystyle \tau \,  \leq \, {{1} \over {\beta}} \,  \leq \, {{1} \over
{\alpha}}\,, $

\smallskip \noindent  $  \displaystyle \!  \! \varphi (\Delta t) \, -\, \Bigl[ \,  
{{\beta^2 \,-\, \alpha^2} \over {2}} \, \tau \, \Delta t \, - \, {{1} \over
{3}}(\alpha^3 + \beta^3) \, \tau \, \Delta t^2 \,  \Bigr]  
\,\,  \leq \,\, {{1} \over {2}}(\alpha^4 + \beta^4)\,  \tau \, 
\Delta t^3  \,.$

\smallskip \noindent 
and relation (2.44) is proven.   Lemma 4 is established.  $ \hfill \square \kern0.1mm $

\bigskip  \bigskip \noindent  {\bf Lemma 5. $\quad$  }

 \noindent 
We suppose that $\, a \ne 0 . \,$  If time step $\, \Delta t \,$ satisfy 
\smallskip \noindent  (2.46) $\quad \displaystyle
0 \,\, < \,\, \Delta t \,\, \leq \,\,{\rm inf} \, \Bigl( \, {{\tau}\over{2}} \,,\,
\abs{a}\, \tau^2 \, \Bigr) \,,\,$
 
\smallskip \noindent 
we have 
\setbox10=\hbox{$   \displaystyle   \, \varphi (\Delta t) \,-\,   \abs{a} \, 
\Delta t  \,  $}
\smallskip \noindent  (2.47) $\quad \displaystyle
 \mod{\box10} \,\, \leq \,\, {7\over12} \, \, \abs{a} \,  \Delta t \,. \, $

\bigskip \noindent  {\bf Proof of Lemma 5.}
\smallskip \noindent $\bullet \quad $ 
If $\,a \,  \neq \, 0$, we observe that~:  $\,\,  \beta^2 \,-\, \alpha^2 \,\,=\,\,
\displaystyle {{2|a|} \over {\tau}} \, > \, 0 \, ,  $ and inequality (2.46)  suppose
that $\,\Delta t$ has been chosen such that   $\quad \displaystyle
\Delta t \,  \leq \,\,  \abs{a}  \, \tau^2 \, . $ Due to previous computations, we
have the following set of estimations~: 

\setbox10=\hbox{$   \displaystyle   \Bigl[  \,  {{\beta^2 \,-\, \alpha^2} \over {2}}
\, \tau \, \Delta t \,  \,-\,  \varphi (\Delta t) \, \Bigr]  $}
\setbox11=\hbox{$  \displaystyle  \mod{\box10} \,\,   $}
\setbox12=\hbox{$ \,\,  \displaystyle \Bigl( 
{{1} \over {3}} \bigl( \alpha^3 + \beta^3 \bigr) \, \Delta t \, \,+\,\, 
{{1} \over {2}} \bigl( \alpha^4 + \beta^4  \bigr) \Delta t^2  \Bigr) \,\,
\tau \, \Delta t $}
\setbox40= \vbox {\halign{#&#&# \cr $\box11$  & $\leq$ &  $\box12$ \cr }}
\smallskip \noindent $ \box40 $
\setbox11=\hbox{$ \qquad \qquad \qquad  \qquad $}
\setbox13=\hbox{$ \,\,  \displaystyle \biggl[ {{\Delta t} \over {3}} \Bigl( 
{{1} \over {4 \tau^3}} \,+\, {{3 a^2} \over {\tau}} \Bigr) \,+\,  {{\Delta t^2}
\over {2}}  \Bigl(  {{1} \over {8 \tau^4}} \,+\, {{3 a^2} \over {\tau^2}} \,+\, 2 a^4
\Bigr) \biggr] \,\, \tau \, \Delta t $}
\setbox14=\hbox{$ \,\,  \displaystyle \Bigl( 
{{\Delta t} \over {3 \tau^3}} \,+\, {{\Delta t^2} \over {2\tau^4}} \Bigr)  \,\, \tau
\, \Delta t $}
\setbox40= \vbox {\halign{#&#&# \cr $\box11$  & = &  $\box13$ \cr  & $\leq$ &  $\box14$
\cr  }}
\smallskip \noindent $ \box40 $

\smallskip $ \qquad \qquad \qquad \qquad \qquad \qquad \qquad \qquad $
because $ \displaystyle \,\,  \abs{a} \, \,  \leq \ {{1} \over {2 \tau}} \,\,$ due to
(2.24)

\setbox11=\hbox{$ \qquad \qquad \qquad  \qquad $}
\setbox15=\hbox{$ \,\,  \displaystyle 
 {{{\Delta t}^2} \over {\tau^2}}  \,\biggl(  {{1} \over {3}} \,+\, {{1} \over {2}} 
 |a| \,  \tau \biggr)  \,   $}
\setbox16=\hbox{$ \,\,  \displaystyle 
{{{\Delta t}^2}  \over {\tau^2}}  \,\biggl(  {{1} \over {3}} \,+\, {{1} \over {4}} 
\biggr)  \, \qquad = \qquad {7 \over 12}  {{{\Delta t}^2} \over {\tau^2}} \, $}
\setbox17=\hbox{$ \,\,  \displaystyle   \quad {7 \over 12}  \,\, 
\abs{a} \, \Delta t  \, $}
\setbox40= \vbox {\halign{#&#&# \cr    $\box11$ & $\leq$ &  $\box15$   \cr   & $\leq$
&  $\box16$   \cr & $\leq$ &  $\box17$   \cr  }}
\smallskip \noindent $ \box40 $

\smallskip \noindent 
and relation (2.47) is proven.   $ \hfill \square \kern0.1mm $

%%%%%%%%%%%%%%%%%%%%%%%%%%%%%%%%%%%%%%%%%%%%%%%%%%%%%%%%%%%%%%%%%%%%%%%%%%%%%%%%%%
\bigskip \noindent  {\bf Lemma 6. $\quad$  }

 \noindent 
We suppose that $\, a \,=\,  0 . \,$  If time step $\, \Delta t \,$ satisfy condition
(2.43),  we have 
\setbox10=\hbox{$   \displaystyle    \, \varphi (\Delta t) \,+\, 
{{1}\over{12\,\tau^2}}\, \Delta t^2 \,  $}
\smallskip \noindent  (2.48) $\quad \displaystyle
 \mod{\box10} \,\, \leq \,\, {{1}\over{32\,\tau^2}} \, \Delta t^2  \,  $

% \vfill \eject 

\bigskip \noindent  {\bf Proof of Lemma 6.}
\smallskip \noindent $\bullet \quad $ 
If $\,a \,=\, 0$, then $ \displaystyle \,\beta \,=\, \alpha \, = \, {{1}\over{2 \,
\tau}} \,  $ and we have simply~: 

\smallskip \noindent $  \displaystyle
{{1}\over{3}} \, (\alpha^3+ \beta^3) \,\,\tau \, \Delta t^2 \,  \,= \,\, 
{{2}\over{3}} \, \Bigl( {{1}\over{2\,\tau}} \Bigr)^{\!3} \,\tau \, \Delta t^2 \, \,= \,\,
{{1}\over{12\,\tau^2}} \, \Delta t^2 \, $

\smallskip \noindent $  \displaystyle
{{1}\over{2}} \, (\alpha^4 + \beta^4) \,\,\tau \, \Delta t^3 \,  \,= \,\, 
\Bigl( {{1}\over{2\,\tau}} \Bigr)^4 \,\tau \, \Delta t^3 \, \,= \,\,
{{1}\over{16\,\tau^3}} \, \Delta t^3 \,\, \leq \,\, {{1}\over{32\,\tau^2}} \, \Delta
t^2 \,\,$ due to the relation (2.43). Then the relation (2.48) is a direct of previous
estimation (2.44) established in Lemma 4. That completes  the proof of Lemma 6.    $
\hfill \square \kern0.1mm $

\bigskip  \bigskip \noindent  {\bf Lemma 7. $\quad$  }

 \noindent 
We suppose that function $\,\, t \longmapsto z(t) \,\,$ is defined at relation (2.34)
and that the numerical initial condition $\, d_{\Delta} \,$ satisfies hypothesis
(2.27). We denote by $\,\, h({\scriptstyle \bullet})\,\,$ the expression introduced in
(2.33).  Then we have 
\smallskip \noindent  (2.49) $\quad \displaystyle
\abs{z(0)} \,\,\, \leq \,\,\, {{1}\over{\eta}} \, \, {\rm min} \, \Bigl( C \,,\,
{{1}\over{k \, \tau}} - \eta \Bigr) \,$
\smallskip \noindent  (2.50) $\quad \displaystyle
\abs{h'\bigl(z(0)\bigr)} \,\,\, \leq \,\,\,  k \, \tau \, \Bigl ( {\rm min} \, \Bigl(
C \,,\, {{1}\over{k \, \tau}} - \eta \Bigr) \Bigr)^{\!2}   \,. $

\bigskip \noindent  {\bf Proof of Lemma 7.}
\smallskip \noindent $\bullet \quad $ 
We have, due to relation (2.27)~:   $ \displaystyle \qquad  -C  \, \,\,  \leq \,\,\, 
x^* \, - \, d\ib{\Delta} \, \,\,  \leq \,\,\, {{1} \over {k\,\tau}}  \,-\,  \eta \,.
\quad $ Then following (2.19) and (2.24) we deduce~: 

\smallskip \noindent $ \displaystyle
-C \,-\, {{1} \over {k\,\tau}}  \,\, \leq \,\, x^* \,-\, {{1} \over {k\,\tau}} \,-\,
 d\ib{\Delta}  \,\, = \,\, x_- \,-\,  d\ib{\Delta}\, \, \leq \,\,\, - \eta \,\qquad $ 
% \smallskip \noindent
and in consequence  

\smallskip \noindent  (2.51) $\quad \displaystyle
{{1}\over{ \abs{ x_- - d_{\Delta} }}} \,\,\, \leq \,\,\, {{1}\over{\eta}} \,. $

\smallskip \noindent
From relation (2.34) we have  $\quad \displaystyle z(0) \,=\,\, {{d_{\Delta} -
x^*}\over{x_- - d_{\Delta}}} \,\,$ and inequality (2.49) is a direct consequence of
(2.51) and of hypothesis (2.27). We derive now the expression (2.33) relatively to
variable $\,\, \xi \,\,$ and we have easily~: 

\smallskip \noindent $ \displaystyle
\abs{h'(z(0))} \,\,\, \leq \,\,\,  {{\abs{x_- \,- \, x^*}} \over{(1 \,+\, z(0))^{\!2}}}
\,\,\, = \,\,\, {{  \abs{x_- \,- \,d_{\Delta} }^2 }\over{ \abs{x_- \,- \, x^*} }}
\,\,\, \leq \,\,\, k \, \tau \, \Bigl ( {\rm min} \, \Bigl( C \,,\, {{1}\over{k \,
\tau}} - \eta \Bigr) \Bigr)^{\!2}     \,$

\smallskip \noindent
due to the above expression of $\, z(0) \,$ and  relation (2.51). The estimation
(2.50) is established and Lemma 7 is proven.  
$ \hfill \square \kern0.1mm $

\bigskip  \noindent  {\bf Proof of Theorem 1.}

\smallskip \noindent $\bullet \quad $ 
We cut the expression inside the absolute value of left hand side of (2.28) into two
parts. The first one is the error  for the continuous solution of the
ordinary differential equation when changing initial data and the second one is to the
error $\,\epsilon\ib{j}\, $ between the solution of the ordinary differential equation
and the discrete solution given  by the scheme for the same initial condition~: 

\smallskip \noindent $ \displaystyle
\abs{ x(j \Delta t\,;\, d) \,-\, x_j(\Delta t \,;\, d\ib{\Delta})} \,\, \,\leq \,\,\, 
\abs{ x(j \Delta t\,;\, d) \,-\,  x(j \Delta t\,;\, d\ib{\Delta})} \,+\,\, 
\abs{\epsilon\ib{j}} \,$

\smallskip \noindent 
and due to definition (2.30), 

\smallskip \noindent  (2.52) $\quad \displaystyle
\abs{ x(j \Delta t\,;\, d) \,-\, x_j(\Delta t \,;\, d\ib{\Delta})} \,\, \,\leq \,\,\, 
\abs{ y(j \Delta t\,;\, d) \,-\,  y(j \Delta t\,;\, d\ib{\Delta})} \,+\,\, 
\abs{\epsilon\ib{j}} \,. \, $

\bigskip \noindent  $\bullet \quad$
We first study the term $\,\abs{ y(j \Delta t\,;\, d) \,-\,  y(j \Delta t\,;\,
d\ib{\Delta})} \,$  in right hand side of (2.52).  We first minorate  the absolute
value  of the denominator of the right hand side of (2.31).  Under the hypothesis
(2.27), we have~:

\smallskip \noindent  $ \displaystyle
1\,+\,k\,\tau\, (d\ib{\Delta}\,-\,x^*)\, (1\,-\,e^{\displaystyle -{t/\tau}})\, \geq
\,\,\, 1 \quad  $ if $ \,\,  d\ib{\Delta}\,-\,x^* \,>\, 0 \,,$

\smallskip \noindent
and in the other case when $ \,\, d\ib{\Delta}-x^* \,\leq \, 0 \,\, $ then
$\,\displaystyle \, |d\ib{\Delta}\,-\,x^*| \,  \leq \, {{1} \over
{k\tau}}\,-\,\eta\,\,\, $  due to (2.27) and   we have in consequence 

\smallskip \noindent  $ \displaystyle
1\,+\,k \, \tau \, (d\ib{\Delta}\,-\,x^*) \, (1\,-\,e^{\displaystyle -{t/\tau}})
\,\,\,  \geq \, \,\, 1\,-\,(1 -k\, \tau \,  \eta)  \,=\,k\, \tau \, \eta \,\, > 0
\,\,\,   $ if $ \,\,  d\ib{\Delta}\,-\,x^* \, \leq \, 0 \,.$

\smallskip \noindent
Therefore the denominator of right hand side of (2.31) is always
strictly positive and is in all cases minorated by $\,k \, \tau \, \eta\, \, $~:

\smallskip \noindent  $ \displaystyle
1\,+\,k \, \tau \, (d\ib{\Delta}\,-\,x^*) \, (1\,-\,e^{\displaystyle -{t/\tau}})
\,\,\,  \geq \,\,k \, \tau \, \eta\, \qquad \forall \,  t \,>\, 0 \,.$

\smallskip \noindent 
In consequence of (2.33) and previous algebra,   $\,y(t\,;\,d)\, $ and
$\,y(t\,;\,d\ib{\Delta})\, $ are well defined for each $\,t\, $ such that  $\, 0 \, 
\leq \, t < \infty\, $ and we have~: 

\smallskip \noindent  $ \displaystyle
x(j \, \Delta t\,;\,d)\,-\,x(j \Delta t\,;\,d\ib{\Delta}) \,\,= \,\,
y(j \, \Delta t\,;\,d)\,-\,y(j \Delta t\,;\,d\ib{\Delta}) \,$
\smallskip \noindent  $ \displaystyle \qquad \quad 
= \,\,  {{(d \,-\, d\ib{\Delta}) \, e^{\displaystyle -{t/\tau}}} \over { \bigl( \, 
1\,+\,k \, \tau \, (d\,-\,x^*) \, (1\,-\,e^{\displaystyle -{t/\tau}}) \bigr) \,
\bigl(  1\,+\,k \, \tau \, (d\ib{\Delta}\,-\,x^*) \, (1\,-\,e^{\displaystyle
-{t/\tau}}) \bigr) }} \,,$

\smallskip \noindent
then~:

\smallskip \noindent  (2.53) $\quad \displaystyle
\abs{ x(j \Delta t,d)\,-\,x(j \Delta t,d\ib{\Delta}) } \quad \leq \,\,  \biggl( {{1}
\over {k\tau \eta}} \biggr)^{2} \, \abs {d \,-\, d\ib{\Delta}} \,$			

\smallskip \noindent 
and the estimation of the first term  introduced in (2.52) is then controlled .

\bigskip \noindent  $\bullet \quad$
We now study the error $\,\epsilon\ib{j}\, $. Due to lemmas 2 and 3 and in particular
relations  (2.36) and (2.42), we have

\setbox10=\hbox{$   \displaystyle    \, 1 \,-\, {\rm exp}  \Bigl( \theta_j  \, 
\varphi(\Delta t) \, \Bigr)  \,  $}
\smallskip \noindent  (2.54) $\quad \displaystyle
\abs{ \epsilon\ib{j} } \,\,\, \leq \,\,\,  \abs {  z(0)   } \,\, 
\abs { h' \bigl( z(0) \bigr)  } \,\, \,   e^{\displaystyle \,- \theta_j}  \,\,\,\,    
\mod{\box10}\,.$

\smallskip \noindent
{\bf (i)} $\quad$  If $\, a = 0 ,\,$ then  $\, \varphi(\Delta t) \,$ is negative due
to relation (2.48) and $\displaystyle \, \theta_j\,\,=\,\, {{j\, \Delta t}\over
{\tau}} \,\,\,  $ remains positive.  We deduce that in this case 

\setbox10=\hbox{$   \displaystyle    \, 1 \,-\, {\rm exp}  \Bigl( \theta_j  \, 
\varphi(\Delta t) \, \Bigr)  \,  $}
\smallskip \noindent $ \displaystyle
\mod{\box10} \,\,\, \leq \,\,\,  \,   1 \,-\, {\rm exp}  \Bigl( \theta_j  \, 
\varphi(\Delta t) \, \Bigr) \,\,\, \leq \,\,\,\,\, \abs{ \theta_j  \,  \varphi(\Delta
t) } \,\,\, = \,\,\, \theta_j \,\, \abs{  \varphi(\Delta t) } \, \, $
 
\smallskip \noindent
and in consequence of (2.54) and (2.48), 

\smallskip \noindent $ \displaystyle
\abs{ \epsilon\ib{j} } \,\,\, \leq \,\,\,   \abs {  z(0)   } \,\,  \abs { h' \bigl(
z(0) \bigr)  } \,\, \, \Bigl(  \,  e^{\displaystyle \,- \theta_j}  \,\,\,\,  \theta_j
\, \Bigr)  \, \, \Bigl( {1\over12} \,+\, {1\over32} \Bigr) \, {{\Delta
t^2}\over{\tau^2}} \,$
\smallskip \noindent $ \displaystyle \qquad 
\leq \,\,\, {{11}\over{96}} \, \,  \abs {  z(0)   } \,\,  \abs { h' \bigl( z(0)
\bigr)  } \,\,  \sup_{ \displaystyle \theta  \geq 0} \Bigl( \theta \, e^{\displaystyle
\,- \theta} \Bigr)  \,\, {{\Delta t^2}\over{\tau^2}} \,$

\smallskip \noindent $ \displaystyle \qquad 
\leq \,\,\, {{11}\over{96}} \,\, {{1}\over{\eta}} \,\,  {\rm min} \, \Bigl( C \,,\,
{{1}\over{k \, \tau}} - \eta \Bigr) \,\, k \, \tau \,\,  \Bigl( {\rm min} \, \Bigl( C
\,,\, {{1}\over{k \, \tau}} - \eta \Bigr) \Bigr)^{\!2} \,\, {{1}\over{e}} \,\,
 {{\Delta t^2}\over{\tau^2}} \,$

\smallskip \noindent $ \displaystyle \qquad 
\leq \,\,\, {{11}\over{96 \, e}} \, {{k}\over{\tau  \eta}} \,  \Bigl( {\rm min} \,
\Bigl( C \,,\, {{1}\over{k \, \tau}} - \eta \Bigr) \Bigr)^{\!3} \,\,  \Delta t^2 \,$

\smallskip \noindent 
and relation (2.29) is established with 

\smallskip \noindent  (2.55) $\quad \displaystyle
A \,\,= \,\, {\rm inf} \, \biggl( \,  \Bigl( {{1} \over {k\tau \eta}} \Bigr)^{\!2}
\,,\,  {{11}\over{96 \, e}} \, {{k}\over{\tau  \eta}} \,  \Bigl( {\rm min} \,
\Bigl( C \,,\, {{1}\over{k \, \tau}} - \eta \Bigr) \Bigr)^{\!3} \,\, \,   \biggl) \,$ 
\smallskip \noindent  (2.56) $\quad \displaystyle
B \,\,= \,\, {{\tau}\over{2}} \,. \,$

\smallskip \noindent
{\bf (ii)} $\quad$  If $\, a \neq 0 ,\,$ then  $\, \varphi(\Delta t) \,$ is positive
due to relation (2.47) if  $\, \Delta t \,$ satisfies condition (2.46).  We suppose
moreover that time step satisfies also the condition 

\smallskip \noindent  (2.57) $\quad \displaystyle
\Delta t \,\, \leq \,\, {{6}\over{19}} \, {{1}\over{\abs{a}}} \,$

\smallskip \noindent
and due to (2.47), we have 

\smallskip \noindent  (2.58) $\quad \displaystyle
\varphi(\Delta t) \,\leq \, \bigl( 1 + {{7}\over{12}} \bigr) \, \abs{a} \, \Delta t
\,\leq \, {1\over2} \,. \, $ 

\smallskip \noindent
In order to majorate the expression $ \, \, \,   e^{\displaystyle \,- \theta_j} 
\,\,     \mod{$   \displaystyle    \, 1 \,-\, {\rm exp}  \Bigl( \theta_j  \, 
\varphi(\Delta t) \, \Bigr)  \,  $}\,$ we distinguish between two cases. On one hand,
when $\,\, \theta_j \, \varphi(\Delta t) \, \leq \, 1 \,, $ we have by convexity of
the exponential function over the interval $\, [0 ,\, 1] \,$~:

\smallskip \noindent  $ \displaystyle
0 \,\, \leq \,\,  e^{\displaystyle \,\theta_j \, \varphi(\Delta t)} \,-\,1 \,\,
\leq \,\, (e-1) \,\theta_j \, \varphi(\Delta t) \qquad  $ and we deduce
 
\smallskip \noindent  $ \displaystyle
e^{\displaystyle \,- \theta_j}  \,\,     \mod{$   \displaystyle    \, 1 \,-\, {\rm
exp}  \Bigl( \theta_j  \,  \varphi(\Delta t)  \Bigr)     $} \,\, \leq  \,\, 
(e-1) \, \bigl[ \theta_j \,  e^{\displaystyle \,-\theta_j} \bigr] \, \varphi(\Delta t)
\,\, \leq  \,\,  {{19}\over{12}} \, {{e-1}\over{e}} \, \abs{a} \, \Delta t \,. \, $

\smallskip \noindent 
The previous inequality and estimation (2.54) show that under hypotheses (2.46) and
(2.57) concerning the time step, we have 

\smallskip \noindent $ \displaystyle
\abs{ \epsilon\ib{j} } \,\,\, \leq \,\,\,  {{19}\over{12}} \, {{e-1}\over{e}}
\,\,  {{1}\over{\eta}} \,  {\rm min} \, \Bigl( C \,,\, {{1}\over{k \, \tau}} - \eta
\Bigr) \,\, k \, \tau \,\,  \Bigl( {\rm min} \, \Bigl( C \,,\, {{1}\over{k \, \tau}} -
\eta \Bigr) \Bigr)^{\!2} \,\,\abs{a} \, \Delta t \,  $

\smallskip \noindent $ \displaystyle
\qquad  \leq \,\,\,  {{19}\over{12}} \, {{e-1}\over{e}} \,\,  {{\abs{a} \,  k \, \tau
}\over{\eta}} \,  \Bigl( {\rm min} \, \Bigl( C \,,\, {{1}\over{k \, \tau}} - \eta
\Bigr) \Bigr)^{\!3} \,\, \Delta t \qquad   $ due to (2.49) and (2.50) 

\smallskip \noindent $ \displaystyle
\qquad  \leq \,\,\,  {{19}\over{24}} \, {{e-1}\over{e}} \,\,  {{k}\over{\eta}} \, 
\Bigl( {\rm min} \, \Bigl( C \,,\, {{1}\over{k \, \tau}} - \eta
\Bigr) \Bigr)^{\!3} \,\, \Delta t \qquad   \qquad  \,\,\, $ due to (2.24)

\smallskip \noindent 
and 

\smallskip \noindent  (2.59) $\quad \displaystyle
\abs{ \epsilon\ib{j} } \,\,\, \leq \,\,\,   {{19}\over{12}} \, {{1}\over{e}} \,\, 
{{k}\over{\eta}} \, 
\Bigl( {\rm min} \, \Bigl( C \,,\, {{1}\over{k \, \tau}} - \eta
\Bigr) \Bigr)^{\!3} \,\, \Delta t \quad
{\rm when} \quad   \theta_j \, \varphi(\Delta t) \, \leq \, 1 \,. \,  $

\smallskip \noindent
On the other hand when  $\,\, \theta_j \, \varphi(\Delta t) \, \geq \, 1 \,, $ we
have~: 

\smallskip \noindent  $ \displaystyle
e^{\displaystyle \,- \theta_j}  \,\,     \mod{$   \displaystyle    \, 1 \,-\, {\rm
exp}  \Bigl( \theta_j  \,  \varphi(\Delta t)  \Bigr)     $} \,\,  = $ 
\smallskip \noindent  $ \displaystyle
\qquad \qquad  =  \,\, 
e^{\displaystyle \,- \theta_j / 2 \, } \,\, e^{\displaystyle \, \theta_j \,
\bigl( \varphi(\Delta t) - 1/2 \bigr) } \,\,\mod{$   \displaystyle    \, 1 \,-\, {\rm
exp}  \Bigl( - \theta_j  \,  \varphi(\Delta t)  \Bigr) $} \,$ 
\smallskip \noindent  $ \displaystyle
\qquad \qquad \leq \quad  e^{\displaystyle \,- \theta_j / 2 \, } \,\,\,\mod{$  
\displaystyle    \, 1 \,-\, {\rm exp}  \Bigl( - \theta_j  \,  \varphi(\Delta t) 
\Bigr) $} \qquad \,$ due to (2.58) 
\smallskip \noindent  $ \displaystyle
\qquad \qquad \leq \quad  e^{\displaystyle \,- \theta_j / 2 \, } \,\,\theta_j  \, 
\varphi(\Delta t) \qquad \,$ because $\,\, \theta_j \,\,$ and $\,\,  \varphi(\Delta t)
\,\,$ are both positive 
\smallskip \noindent  $ \displaystyle
\qquad \qquad \leq \quad 2 \,  \sup_{ \displaystyle \theta  \geq 0} \Bigl( \theta 
\, e^{\displaystyle \,- \theta} \Bigr)  \,\, \varphi(\Delta t) 
\qquad \leq \qquad {{19}\over{12}} \, {{2}\over{e}} \,\abs{a} \, \Delta t \qquad
\qquad \,$

\smallskip \noindent
thanks to relation (2.58). Following inequality (2.54), we obtain in this second case 

\smallskip \noindent $ \displaystyle
\abs{ \epsilon\ib{j} } \,\, \,  \leq \,\,\,  {{k\, \tau }\over{\eta}} \,  \Bigl( {\rm
min} \, \Bigl( C \,,\, {{1}\over{k \, \tau}} - \eta \Bigr) \Bigr)^{\!3} \,\, 
{{19}\over{6\, e}} \, \abs{a} \, \Delta t \,$ 

\smallskip \noindent  (2.60) $\quad \displaystyle
\abs{ \epsilon\ib{j} } \,\, \leq \,\,  {{19}\over{12\, e}} \,  {{k}\over{\eta}} \,  
\Bigl( {\rm min} \, \Bigl( C \,,\, {{1}\over{k \, \tau}} - \eta \Bigr) \Bigr)^{\!3} \,\, 
\Delta t \quad {\rm when} \quad   \theta_j \, \varphi(\Delta t) \, \geq \, 1 \, \,  $

\smallskip \noindent
and relation (2.29) is proved for this case with 

\smallskip \noindent  (2.61) $\quad \displaystyle
A \,\,= \,\,  {\rm inf} \, \biggl( \,  \Bigl( {{1 } \over {k\tau \eta}}
\Bigr)^{\!2} \,,\,  {{19}\over{12\, e}} \,  {{k}\over{\eta}} \,  
\Bigl( {\rm min} \, \Bigl( C \,,\, {{1}\over{k \, \tau}} - \eta \Bigr) \Bigr)^{\!3} 
\,\, \, \,   \bigg)
\,$ 

\smallskip \noindent  (2.62) $\quad \displaystyle
B \, \,= \,\,   {\rm inf} \, \biggl( \, {{\tau}\over{2}}\,,\, \abs{a} \,\tau^2 \,,\,
{{6}\over{19}} \, {{1}\over{\abs{a}}} \, \biggr) \,$

\smallskip \noindent 
which ends the proof of Theorem 1. $ \hfill \square \kern0.1mm $

%%%%%%%%%%%%%%%%%%%%%%%%%%%%%%%%%%%%%%%%%%%%%%%%%%%%%%%%%%%%%%%%%%%%%%%%%%%%%%%%%%%%%%%%
%%%%%%%%%%%%%%%%%%%%%%%%%%%%%%%%%%%%%%%%%%%%%%%%%%%%%%%%%%%%%%%%%%%%%%%%%%%%%%%%%%%%%%%%
 \bigskip \bigskip %% \vfill \eject  
\noindent  {\smcaps 3) $ \quad$ Matrix Riccati equation.}

\smallskip \noindent $\bullet \quad$
In order to define a numerical scheme to solve the Riccati differential equation
(1.22) with initial condition (1.19) we first introduce a real number 
$\,\mu$, which is chosen positive and such that the real matrix 
$\,\, \bigl[ \mu I \,-\, (A \,+\, A^{\rm \displaystyle t}) \bigr] \,\, $ is 
definite positive~:

\smallskip \noindent  (3.1) $\quad \displaystyle
\mu \, > \, 0\,, \qquad  {{1} \over {2}}(\mu \, x \,,\, x) \,-\, (A \, x \,,\, x) > 0
\, , \quad \forall \,  x \,  \neq \, 0 \,.$

\smallskip \noindent 
Then we define a definite positive matrix $\,M \, $ wich depends on strictly positive 
scalar  $\,\mu$ and matrix $\,A$ :

\smallskip \noindent  (3.2) $\quad \displaystyle
M \,\,=\,\, {{1} \over {2}} \mu \, I \,-\, A  \, .$

\smallskip \noindent 
The numerical scheme is then defined in analogy with relation (2.6). We have the
following decomposition :

\smallskip \noindent  (3.3) $\quad \displaystyle
A \, \,\,=\,\, \, A^+\, \,-\, A^-  \,\, ,  $

\smallskip \noindent 
with

\smallskip \noindent  (3.4) $\quad \displaystyle
A^+\, \,\,=\,\, \mu I  \, , \quad \, A^- \,\,=\,\, M, \quad  \mu > 0 ,\quad M \,\,$
 positive definite .
 
\smallskip \noindent
Taking as an explicit part the positive contribution $\,A^+\,$ of the decomposition
(3.3) of matrix $\,A$ and in the implicit part the negative contribution $\,A^-
\,=\, M$ of the  decomposition (3.3), we get the following {\bf harmonic scheme}  : 

\setbox20=\hbox{$\displaystyle  
{{1} \over {\Delta t}} (X_{j+1} \,-\, X_j) \,+\, {{1} \over {2}} (X_j K X_{j+1}
\,+\, X_{j+1} K X_j) \,+\,  $}
\setbox21=\hbox{$\displaystyle  \qquad  \qquad  \qquad  \qquad  \qquad \qquad  
\,+\, (M^{\rm \displaystyle t} X_{j+1} \,+\, X_{j+1} M) \,\,=\,\, \mu
X_j \,+\, Q  \,.$}
\setbox23= \vbox {\halign{#\cr \box20 \cr \box21 \cr}}
\setbox24= \hbox{ $\vcenter {\box23} $}
\setbox25=\hbox{\noindent (3.5) $\, \, \left\{  \box24 \right. $}
\smallskip \noindent $ \box25 $

\smallskip \noindent 
The numerical solution $\,X_{j+1}$  given by the scheme at time step $\,(j+1) \Delta
t \, $ is then defined as a solution of Lyapunov matrix equation with matrix $\,X$ as
unknown : 

\smallskip \noindent  (3.6) $\quad \displaystyle
S^{\rm \displaystyle t}_j \, X \,\,+\, \, X \, S_j \,\,=\,\, Y_j \, . $

\smallskip \noindent 
with :

\smallskip \noindent  (3.7) $\quad \displaystyle
S_j \, \,\,=\,\, \, {{1} \over {2}}  I \,+\, {{\Delta t} \over {2}} K X_j \,+\,
\Delta t \, M \,,$

\smallskip \noindent 
and :

\smallskip \noindent  (3.8) $\quad \displaystyle
Y_j \, \,\,=\,\, \, X_j \,+\, \mu \, \Delta t \, X_j \,+\, \Delta t \, Q  \,.$

\smallskip \noindent 
We note that $\,S_j \, $ is a (non necessarily symmetric) positive matrix and that
$\,Y_j \, $ is a symmetric definite positive matrix if it is the case for $\,X_j$. In
all cases, matrix $\, Y_j \, $ is a symmetric positive matrix.

\bigskip  \noindent  {\bf  Definition 2. \quad Symmetric matrices. }

 \noindent 
Let $n$ be an integer greater or equal to $1$. We define by $\,{\cal{S}}_n(\R)$, 
(respectively $\,{\cal{S}}^+_n(\R)$ , $\,{\cal{S}}^{+*}_n(\R)$ ) the linear
space (respectively the closed cone, the open cone) of symmetric-matrices (respectively
symmetric positive and symmetric definite positive matrices) ; we have 
\smallskip \noindent  (3.9) $\quad \,\, \displaystyle
(x\,,\,Sy) \,\,=\,\, (Sx,y) \,,  \quad  \forall \,  x , y \in \R^n \,\,, \quad  
\forall \,  S \in {\cal{S}}_n(\R)\,,$
\smallskip \noindent  (3.10) $\quad \displaystyle
(x\,,\,Sx) \,\,\geq \,\, 0 \,, \quad  \forall \,  x  \in \R^n \,,
\quad   \forall \,  S \in {\cal{S}}^{+}_n(\R)\,,$
\smallskip \noindent  (3.11) $\quad \displaystyle
(x\,,\,Sx) \,\,\,{\rm >}  \,\,\, 0 \,, \quad  \forall \,  x  \in \R^n \,, x \neq 0 \,, 
\quad   \forall \,  S \in {\cal{S}}^{+*}_n(\R)\,.$
\smallskip   \noindent
The following inclusions  $\,\, {\cal{S}}^{+*}_n(\R) \, \subset \, {\cal{S}}^+_n(\R) 
\, \subset \, {\cal{S}}_n(\R) \,\,\, $ are natural.

\bigskip 
% \vfill \eject 
\noindent  {\bf Proposition 3. \quad Property of the Lyapunov equation.}

\noindent $\bullet \quad$
Let $\,S$ be a matrix which is not necessary symmetric, such that the associated
quadratic form:  $\,\R^n \ni x \longmapsto (x,Sx) \in \R $, is strictly positive : 
\smallskip \noindent  (3.12) $\quad \displaystyle
S \,+\, S^{\rm \displaystyle t} \,\, \in {\cal{S}}^{+*}_n(\R) \,. $

\smallskip \noindent 
Then the application $\,\, \varphi\ib{S} \,\, $ defined by :
\smallskip \noindent  (3.13) $\quad \displaystyle
{\cal{S}}_n(\R) \, \ni \, X \, \longmapsto \, \varphi\ib{S}(X) \, \,\,=\,\, \,
S^{\rm \displaystyle t}\, X \,\,+\, \, X \, S \, \,\, \in \, {\cal{S}}_n(\R) \,,$
\smallskip \noindent 
is a one to one bijective application on the space $\,{\cal{S}}_n(\R)$ of real
 symmetric matrices of order $\,n$.

% \smallskip }\noindent {\vrule height \ht12 depth 0pt width 0,05cm} \quad \box12
% \vskip 2pt \setbox12=\vbox{\hsize=12,0cm  \noindent 
 
\smallskip \noindent $\bullet \quad$
Morever, if matrix $\,\varphi\ib{S}(X)\, $   is positive (respectively definite 
positive)  then the matrix $\,X\, $ is also  positive (respectively definite
positive)~: 
\smallskip \noindent  $ \displaystyle
\bigl(\varphi\ib{S}(X) \in {\cal{S}}^{+}_n(\R)  \Longrightarrow   X \in
{\cal{S}}^{+}_n(\R) \bigr) \, \,{\rm and}   \,\,  \bigl(\varphi\ib{S}(X) \in
{\cal{S}}^{+*}_n(\R)  \Longrightarrow  X \in {\cal{S}}^{+*}_n(\R) \bigr) \, . \,$

\bigskip \noindent  {\bf Proof of  proposition 3.}

\noindent $\bullet \quad $ 
We first observe that $\,\varphi\ib{S}\,$ is a linear map. In the case where $\,S\,$
 is a symmetric matrix, we can immediatly deduce from (3.12) that $\,S\,$ is definite
positive. Let $\,X\,$ be a matrix such that :

\smallskip \noindent  (3.14) $\quad \displaystyle
\varphi\ib{S}(X) \, \,\,=\,\, \, 0 \,.$

\smallskip \noindent 
We prove that $\,X \,=\, 0$, ie that $\,{\rm ker} \, \varphi\ib{S} \,\,=\,\,\{0\}\, $.
Let $\,x\, $ be an eigenvector of matrix $\,S\, $ associated to eigenvalue 
$\,\lambda\, $ :

\smallskip \noindent  (3.15) $\quad \displaystyle
S \, x \, \,\,=\,\, \, \lambda \, x \, , \, x \ne 0\,.$

\smallskip \noindent 
We deduce from relations (3.14) and (3.15) and the fact that matrix $\, S \,$ is
supposed to be symmetric the equality 

\smallskip \noindent  (3.16) $\quad \displaystyle
S ( X x ) \, \,\,=\,\, \, - \, \lambda \, ( X x ) \, .$

\smallskip \noindent 
According to the previous hypothesis the negative scalar $\, - \lambda \, $ cannot be 
an eigenvalue of $\,S\, $, then the vector $\, ( X \, x ) \, $ must be equal to zero.
This is true  for each eigenvalue of matrix $\,S$, wich prove the property in this
case, because  $\,\varphi\ib{S}\, $ is also an endomorphism from $\,{\cal{S}}_n(\R)\, $ to
$\,{\cal{S}}_n(\R)\, $.

\bigskip \noindent  $\bullet \quad\,$
In the case where $\,S\,$ is not symmetric we suppose first that $\,S\,$ is composed
into blocks of Jordan type. We first study the case where $\,S= \Lambda \, $ is
composed  of exactly one Jordan block associated to eigenvalue $\,\lambda\,$ :
 
\smallskip \noindent  (3.17) $\quad \displaystyle
\Lambda \, \,\,=\,\,  \pmatrix {\lambda & 1 & 0 & \cdots & 0  \cr
0 & \lambda & 1 &  \cdots & 0 \cr  \vdots & \ddots & \ddots & \ddots  & \vdots  \cr  
0 & \cdots &  0 & \lambda & 1\cr 0 & 0 &  \cdots & 0 & \lambda \cr }\, $

\smallskip \noindent 
and according to the hypothesis (3.12), we have : 

\smallskip \noindent  (3.18) $\quad \displaystyle
{\rm Re} \, \lambda \, > \, 0 \, .$

\smallskip \noindent 
Let $\,X_{i,j}\,$ be an element of the matrix $\,X\,$, we have from the following 
relation :

\smallskip \noindent  (3.19) $\quad \displaystyle
(\Lambda^{\rm \displaystyle t} X \,+\, X \Lambda)_{i,j} \,\,=\,\, 2 \lambda X_{i,j}
\,+\, X_{i-1,j} \, +\, X_{i,j-1} \,\qquad  {\rm if} \, i \, \geq \, 2 \,\,\, {\rm  and}
\, j \, \geq \, 2  \, , \,$

\smallskip \noindent
and :

\smallskip \noindent  (3.20) $\quad \displaystyle
(\Lambda^{\rm \displaystyle t} X \,+\, X \Lambda)_{1,1}   \,\,=\,\, 2 \, \lambda \,
X_{1,1} \, $
\smallskip \noindent  (3.21) $\quad \displaystyle
(\Lambda^{\rm \displaystyle t} X \,+\, X \Lambda)_{1,j}  \,\,=\,\,  2 \, \lambda \,
X_{1,j} \,\, +\,\,  X_{1,j-1} \quad \qquad j \ge 2  \,.$

\smallskip \noindent
Because $\, \lambda \ne 0, \,X_{1,1} \,=\, 0\,$ from (3.20), then $\,X_{1,j}
\,\,=\,\, 0\,$ if $\, j \ge 2\,$ from (3.21).  By induction, $\,X_{i,j} \,=\, 0\,$
taking into account (3.13).

\bigskip \noindent  $\bullet \quad$
When matrix $\, S \, $ is composed by a family of Jordan blocks of the previous type,
{ \it i.e.} $\, S \,= \, \Lambda \,=\, ({\rm diag} \, \Lambda_j) \,$ where $\, 1 \leq  j
\leq  p \, $ and $\,  \Lambda_j \,$ is a Jordan block of the type (3.17) and of order
$\, n_j  \, $ 

\smallskip \noindent  (3.22) $\quad \displaystyle
\Lambda \,\, = \,\,   \pmatrix {\Lambda_1 & 0 & 0 & \cdots & 0  \cr
0 & \Lambda_2 & 0 &  \cdots & 0 \cr  \vdots & \ddots & \ddots & \ddots  & \vdots  \cr  
0 & \cdots &  0 & \Lambda_{p-1} & 0 \cr 0 & 0 &  \cdots & 0 & \Lambda_p \cr }\,, \,  $

\smallskip \noindent
we decompose the matrix $\, X \,$ into blocks $\, X_{i,j} \,$ of order $\,
n_i \times n_j \,$ : 

\smallskip \noindent  (3.23) $\quad \displaystyle
X \,\, = \,\,   \pmatrix { X_{1,1} &  X_{1,2} &  X_{1,3} & \cdots &  X_{1,p}  \cr
 X_{2,1} & X_{2,2} &  X_{2,3} &\cdots & X_{2,p} \cr  \vdots & \ddots & \ddots & \ddots 
& \vdots  \cr   X_{p-1,1}  & \cdots &   \cdots &  X_{p-1,p-1} &  X_{p-1,p} \cr 
X_{p,1} & X_{p,2} &  \cdots & X_{p,p-1} & X_{p,p} & \cr }\,. \,  $

\smallskip \noindent 
Then the block number $\, (i,j) \,$ of the expression 
$Ê\, S^{\rm \displaystyle t}\, X \,+ \, X \, S \, \,$ is equal to $Ê\,
\Lambda_{i}^{\rm \displaystyle t}\, X_{i,j} \,\,+\, \, X_{i,j} \, \Lambda_{j} \, \,$
and for $\, i\ne j \,$ we  have to prove that the nondiagonal matrix $\, 
X_{i,j} \,$  is identically null.

\bigskip \noindent  $\bullet \quad$
We establish  that if $\, \Lambda \,$ is a Jordan block  of order $\,n \,$ of the
type (3.17) satisfying inequality (3.18), if    $\, M \,$ is a second   Jordan block 
of order $\, m \,$ of the previous type, 

\smallskip \noindent  (3.24) $\quad \displaystyle
M \, \,\,=\,\,  \pmatrix {\mu & 1 & 0 & \cdots & 0  \cr
0 & \mu & 1 &  \cdots & 0 \cr  \vdots & \ddots & \ddots & \ddots  & \vdots  \cr  
0 & \cdots &  0 & \mu & 1\cr 0 & 0 &  \cdots & 0 & \mu \cr }\, $

\smallskip \noindent  
 such that an inequality analogous to (3.18) is valid : 

\smallskip \noindent  (3.25) $\quad \displaystyle
{\rm Re} \, \mu \, > \, 0 \, .$

\smallskip \noindent  
and if $\,X \,$ is a real matrix of order $\, n \times m \,$ chosen such that 

\smallskip \noindent  (3.26) $\quad \displaystyle
\Lambda^{\rm \displaystyle t}\, X \,\,+\, \, X \, M \, \, = \,\, 0 \,,\, $

\smallskip \noindent
then matrix $\, X \,$ is identically equal to zero. As in relations (3.19)-(3.21), we
have~: 

\smallskip \noindent  (3.27) $\quad \displaystyle
(\Lambda^{\rm \displaystyle t} X \,+\, X M)_{1,1}  \,\,=\,\,  \lambda \,  X_{1,1}
\,+\, \mu \,  X_{1,1}  \,$
\smallskip \noindent  (3.28) $\quad \displaystyle
(\Lambda^{\rm \displaystyle t} X \,+\, X M)_{1,j} \,\,=\,\, \lambda \,  X_{1,j}
\,+\, \mu \,  X_{1,j}  \,+ \, X_{1,j-1} \,\,  \qquad j \, \geq \, 2 \, $
\smallskip \noindent  (3.29) $\quad \displaystyle
(\Lambda^{\rm \displaystyle t} X \,+\, X M)_{i,1} \,\,=\,\, \lambda \,  X_{i,1}
\,+\, \mu \,  X_{i,1}  \,+ \, X_{i-1,1} \,\,  \qquad  \,\, i \, \geq \, 2 \, $
\smallskip \noindent  (3.30) $\quad \displaystyle
(\Lambda^{\rm \displaystyle t} X \,+\, X M)_{i,j} \,\,=\,\, \lambda \,  X_{i,j}
\,+\, \mu \,  X_{i,j}  \,+ \, X_{i-1,j} + \, X_{i,j-1} \,\,  \qquad i \,,\,\, 
j \, \geq \, 2 \,. \,  $

\smallskip \noindent
Then due to (3.18), (3.25) and (3.27) we have $\, X_{1,1} = 0 .\,$ By induction on
$\,j \,$ and according to relation (3.28) we have $\, X_{1,j} = 0 .\,$ By induction on
$\,i \,$ and due to relation (3.29) we have analogously $\, X_{i,1} = 0 .\,$ Finally
relation (3.30) prove by induction that $\, X_{i,j} = 0 $ when $ \, i \,$ and $\, j
\, \geq 2 ,  \, $ and matrix $\, X \,$ is identically null.

\bigskip \noindent  $\bullet \quad$
When matrix $\, S \,= \, \Lambda \,= \,  ({\rm diag} \, \Lambda_j) \,$ is composed by
a family of Jordan blocks of the previous type, the nondiagonal blocks $\, X_{i,j} \,$
of decomposition (3.23) are null due to the previous point. Moreover, the diagonal
matrices $\, X_{i,i} \,$ are also identically null due to the property established at 
relations (3.17) to (3.21). Then the proposition is established when matrix $\,S \,
= \, \Lambda  \,$ is composed by Jordan blocks as in  relation (3.22).

\bigskip \noindent  $\bullet \quad$
In the general case where real matrix $\,S\,$ satisfy relation (3.14) there exists
always a complex matrix $\,Q\,$ such that :

\smallskip \noindent  (3.31) $\quad \displaystyle
S \, \,\,=\,\, \, Q^{-1} \, \Lambda \, Q \, , \, $

\smallskip \noindent 
where matrix $\,\Lambda\,$ has a bloc Jordan form given {\it e.g.} by the right hand side
of the  relation (3.22). We have also  the following elementary calculus~:

\setbox20=\hbox{$ \displaystyle \varphi\ib{S}(X) $}
\setbox21=\hbox{$ \displaystyle  S^{\rm \displaystyle t} \, X \,+\, X \, S$}
\setbox22=\hbox{$ \displaystyle  Q^{\rm \displaystyle t} \, \Lambda^{\rm \displaystyle
t} \,  Q^{\rm \displaystyle -t} \,   X \,+\, X\,  Q^{-1}\,  \Lambda \, Q $}
\setbox23=\hbox{$\displaystyle  Q^{\rm \displaystyle t} \, \Lambda^{\rm \displaystyle
t} \, ( Q^{\rm \displaystyle -t} X  Q^{-1} ) \, Q \,+\, Q^{\rm \displaystyle t} \, (
Q^{\rm \displaystyle -t} X Q^{-1})\, \Lambda \,  Q $}
\setbox24=\hbox{$\displaystyle  Q^{\rm \displaystyle t} \, 
( \Lambda^{\rm \displaystyle t} \,  Y \,+\, Y \Lambda  ) \,  Q \,$}
\setbox40= \vbox {\halign{#&#&# \cr $\box20$  & $\,\, = \,\, $ &  $\box21$ \cr
& $\,\, = \,\, $ & $\box22$ \cr  & $\,\, = \,\, $ & $\box23$ \cr
& $\,\, = \,\, $ & $\box24$ \cr}}
\smallskip \noindent $ \box40 $

\smallskip \noindent 
with $\,Y \,\,=\,\, Q^{-t} X Q^{-1}\,$. Because the matrix $\,Q\,$ is inversible and
$\,\varphi\ib{S}(X)\,$ is equal to zero, the matrix 
$\, \smash { (  \Lambda^{\rm \displaystyle t} } \,  Y \,+\, Y \Lambda ) \,$ 
is equal to zero. Matrix $\,\Lambda\,$ is a diagonal bloc Jordan form, and from the
previous points we deduce that  $\,Y\,$ is equal to zero, then $\,X \,=\, 0\,$ and
the proof of first assertion of proposition 3 is established in  the general case.   

\bigskip \noindent  $\bullet \quad$
We suppose now that the matrix $\,\varphi\ib{S}(X)\,$ is symmetric and  positive,
that is~:  $\,\varphi\ib{S}(X) \in {\cal{S}}^+_n(\R) . \, $ Let  $\,x\,$ be an
eigenvector of matrix $\,X\,$ associated to the real eigenvalue $\,\lambda\,$~:

\smallskip \noindent  (3.32) $\quad \displaystyle
X \, x \, \,\,=\,\, \, \lambda \, x \, , \quad  x \,  \neq \, 0 \, .$

\smallskip \noindent 
From the definition of application $\,\varphi\ib{S}\,$, we have the following relation :

\smallskip \noindent  (3.33) $\quad \displaystyle
(\, x \, , \varphi\ib{S}(X)\, x \, ) \, \,\,=\,\, \, 2 \, \lambda \, (x , S x \,)\,$

\smallskip \noindent 
and from hypothesis on $\,\varphi\ib{S}(X)\,$ the left hand side of relation (3.33) is
  positive. The expression $\,(x,Sx)\,$ is also strictly positive since 
 relation (3.8) holds. We  deduce that the number $\,\lambda\,$ is  positive  because
$\,\, X \in {\cal{S}}_n(\R) \,\,$ has an orthogonal decomposition in eigenvector
spaces.

\bigskip \noindent  $\bullet \quad$
If matrix $\,\varphi\ib{S}(X)\,$ is symmetric and  positive definite 
$\,(\varphi\ib{S}(X) \in {\cal{S}}^{+*}_n(\R)) , \, $ we introduce eigenvector $\, x
\,$ of matrix $\, X\, $ as previously (relation (3.32)). Then relation (3.33) remains
true and  the left hand side of this relation is strictly  positive. Then the
eigenvalue  $\,\lambda\,$ of matrix $\, X \,$ is   strictly positive and
proposition 3 is proven.  $ \hfill \square \kern0.1mm $

\bigskip \bigskip \noindent $\bullet \quad $ 
The numerical scheme has been writen as an  equation of unknown $\,X \,=\, X_{j+1}\,$
that takes the form :

\smallskip \noindent  (3.34) $\quad \displaystyle
\varphi\ib{S_j}(X) \, \,\,=\,\, \, Y_j \,.$

\smallskip \noindent 
with $\,\varphi\ib{S_j}\,$ given by a relation of the type (3.13)  with the help of
matrix $\,S_j\,$ defined in (3.7) and $\, Y_j\,$ in (3.8).  Then we have the following
proposition.

\bigskip  \bigskip 
\noindent  {\bf Proposition 4. }

\noindent  {\bf Homographic scheme computes  a definite positive matrix.}

 \noindent 
\smallskip \noindent $\bullet \quad $ 
The matrix $\,X_{j}\,$ defined by  numerical scheme (3.5)  with the initial
condition 
\smallskip \noindent  (3.35) $\quad \displaystyle
X_0 \,\,=\,\, 0\,$
\smallskip \noindent
is  positive for each time step $\,\Delta t > 0  \,$ : 
\smallskip \noindent  (3.36) $\quad \displaystyle
 X_j \in {\cal{S}}^+_n(\R), \,\,\qquad \quad  \forall \,  j \ge 0 . \,$

\smallskip \noindent $\bullet \quad $ 
If there exists some integer $ \,  m \, $ such that $\,  X_m \,$ belongs to the open
cone  $\,  {\cal{S}}^{+*}_n(\R)\,$, then  matrix $\,  X_{m+j} \,$ belongs to the open
cone  $\,  {\cal{S}}^{+*}_n(\R)\,$ for each   $\,j $~: 

\smallskip \noindent  (3.37) $\quad \displaystyle
\Bigl( \, \exists \, m \in \N , \,\, X_m \in {\cal{S}}^{+*}_n(\R)\, \Bigr) \,\,
\Longrightarrow \Bigl( \,  \forall \,  j \ge 0 , \,\,   X_{m+j} \in 
{\cal{S}}^{+*}_n(\R)\, \, \Bigr) \,\,  . \,$

\bigskip \noindent  {\bf Proof of  proposition 4.}
\smallskip \noindent $\bullet \quad $ 
First we have $\,Y_0 \,=\,  \Delta t \, Q\,$ and $\,S_0 \,=\,  {{1} \over {2}} I
\,+\, \Delta t \, M ,\,$ then $\,X_1\,$ is a symmetric  positive matrix $\, (X \in
{\cal{S}}^{+}_n(\R) ) \,$  according to proposition 3 since matrix $\,S_0\,$ is
symmetric  positive and $\,M\,$ has been chosen such that 

\smallskip \noindent  (3.38) $\quad \displaystyle
M \,+Ê\, M^{\displaystyle \rm t} \,\,$ is  positive definite.

\smallskip \noindent
The end of the first point follows by induction. 

\smallskip \noindent  $\bullet \quad $ 
If real symmetric  positive definite matrix  $\, X_{j+1} \,$ is given, relation (3.8)
clearly indicates that matrix $\, Y_{j}\,$ is symmetric positive definite and matrix
$\, S_{j} \,$ introduced at relation (3.7) has a symmetric part which is positive
definite if we verify that the following matrix 

\smallskip \noindent  (3.39) $\quad \displaystyle
K \, X_j \,+\, X_j  \, K \,$ 

\smallskip \noindent 
is positive. But this property is a consequence of the following. On one
hand, matrix $\, K \,$ is positive definite and map $\, X \longmapsto K \, X \,+ \, X
\, K \,$ transforms the closed cone of positive  symmetric matrices onto
himself. On the other hand, matrix $\, X_j \,$ is symmetric  positive definite then
$\,\,  K \, X_j \,+ \, X_j \, K  \, \in {\cal{S}}^{+*}_n(\R) \,\,$ due to proposition
3.  Proposition 4 is established.   $ \hfill \square \kern0.1mm $

\bigskip \noindent  $\bullet \quad$
We have defined a numerical scheme for solving in a approximate way the Riccati
equation (1.20) with the help of relation  (3.5)  and the initial condition 

\smallskip \noindent  (3.40) $\quad \displaystyle
X_0 \,\, = \,\, 0 \,$  

\smallskip \noindent
naturally associated with initial condition (1.19). Recall that time step $\, \Delta t
\,$ is  {\bf not}  limited by any stability condition : matrix $\, X_j \,$ is allways
symmetric  positive and positive definite if matrix $\, Q \,$ is symmetric positive
definite.  Moreover the equation (3.5) that allows the
computation of $\, X_{j+1}\,$ from data is a {\bf linear}  equation whose   unknown
is a   symmetric matrix. But $\, X_{j+1}\,$ is a nonlinear (homographic !) function of
previous iteration matrix $\, X_j .\,$ 

\bigskip \noindent  $\bullet \quad$
We now study  the convergence of the iterate matrix $\, X_j\,$ as long as discrete time
$\, j \, \Delta t \,$ tends to infinity. We know from Proposition 1 that the solution
of the differential Riccati equation tends to the solution of  stationary Riccati
equation~:  

\smallskip \noindent  (3.41) $\quad \displaystyle
X_{\infty}K X_{\infty} \,-\, (A^{\rm \displaystyle t} X_{\infty} \,+\, X_{\infty} A ) 
\,-\, Q \,\,=\,\, 0  \,.$

\smallskip \noindent 
We first study the monotonicity of our numerical scheme. Recall first that if $\, A
\,$ and $\, B \,$ are two real symmetric matrices, the condition

\smallskip \noindent  (3.42) $\quad \displaystyle
A \,\, \leq \,\, B \,$

\smallskip \noindent 
and  respectively the condition 

\smallskip \noindent  (3.43) $\quad \displaystyle
A \,\, < \,\, B \,$

\smallskip \noindent 
means that matrix $\, B-A \,$ is positive $\, \bigl( B-A \, \in {\cal{S}}^{+}_n(\R)
\bigr) \, $ 

\smallskip \noindent  (3.44) $\quad \displaystyle
(x \,,\, (B-A) \, x \bigr) \,\, \geq \,\, 0 \,\,, \quad \forall \, x \in \R^n \,, $

\smallskip \noindent 
and respectively that  matrix $\, B-A \,$ is positive definite  $\, \,
\bigl( B-A \, \in {\cal{S}}^{+*}_n(\R) \bigr)  \,:  $

\smallskip \noindent  (3.45) $\quad \displaystyle
(x \,,\, (B-A) \, x \bigr) \,\, > \,\, 0 \,\,, \quad \forall \, x \in \R^n \,. \,$

\bigskip  \bigskip \noindent  {\bf Proposition 5. Monotonicity.}

\noindent $\bullet \quad $ 
Under the two conditions :

\smallskip \noindent   $ (3.46)  \qquad \quad \displaystyle
Q \,\,$ is a definite positive symmetric matrix 
\smallskip \noindent   $ (3.47)  \qquad \quad \displaystyle
{{1} \over {2} }\, \bigl( K X_{\infty} + X_{\infty} K \bigr)  \,
< \, \bigl( \mu \,+\, {{1} \over {\Delta t}}  \bigr) \, I \,, \, $

\smallskip \noindent 
the scheme (3.5) is monotone and we have  more precisely : 
\smallskip \noindent   $  (3.48)  \qquad \quad \displaystyle
\Bigl( \, 0 \, \,  \leq \, \, X_j \, \,  \leq \, \, X_{\infty} \,\Bigr) \,
\Longrightarrow  \, \Bigl( \, 0 \, \,   \leq \, \, X_j \, \,  \leq \, \, X_{j+1} \, \, 
\leq \, \, X_{\infty} \Bigr) \, .$

\bigskip \noindent  {\bf Proof of  proposition 5.}

\smallskip \noindent $\bullet \quad $ 
We know from Lewis [Le86]  that for symmetric  definite positive matrix  $\, K \,$ 
and symmetric   positive matrix  $\, Q , \,$  the algebraic Riccati equation
(3.41)  has a unique symmetric positive solution $\, X_{\infty} .\,$ Moreover, matrix
$\,  X_{\infty} \,$ is positive definite if  matrix $\, Q \,$  is  positive definite.

\bigskip \noindent $\bullet \quad $ 
We first establish that matrix $\, \Theta \, \equiv \, X_{\infty} -  X_{j+1} \,$ is
positive if the left hand side of implication (3.48) is satisfied. We substract
relation (3.41) from numerical scheme (3.5), observe that 

\smallskip \noindent  $ \displaystyle
X_j \, K \, X_{j+1} \, - \, X_{\infty}\, K \, X_{\infty} \,\, = \,\, 
X_j \, K \,(X_{j+1}  - X_{\infty}) \,+\, (X_{j}  - X_{\infty}) \, K \, X_{\infty} \,$
\smallskip \noindent  $ \displaystyle \qquad \qquad \qquad \qquad \qquad \quad
\,\,\, = \,\, (X_{j}  - X_{\infty}) \, K \, X_{\infty} \,- \,  X_j \, K \, \Theta
\,$ 

\smallskip \noindent 
and that 

\smallskip \noindent  $ \displaystyle
X_{j+1} \, K \, X_{j} \, - \, X_{\infty}\, K \, X_{\infty} \,\, = \,\, 
(X_{j+1}  - X_{\infty}) \, K \,  X_{j}  \,+\, X_{\infty} \, K \,( X_{j} - 
X_{\infty}) \,$ 
\smallskip \noindent  $ \displaystyle \qquad \qquad \qquad \qquad \qquad \quad
\,\,\, = \,\, X_{\infty} \, K \, ( X_{j} \, - \, X_{\infty}) \,-\, \Theta \, K \, X_j
\,$ 

\smallskip \noindent 
then we obtain the following equation for $\, \Theta \,$ : 

\setbox20=\hbox{$\displaystyle  
\varphi\ib{\Sigma_1} (\Theta) \quad \equiv \quad 
{{1}\over{\Delta t}} \, \Theta \,+\, {1\over2} \, \bigl( X_j \, K \, \Theta + \Theta
\, K \, X_j \bigr) \,+\, (M^{\displaystyle \rm t} \, \Theta + \Theta \, M) \,\,= \,\, 
$}
\setbox21=\hbox{$\displaystyle  \qquad \qquad \,\,\, 
= \,\,{{1}\over{\Delta t}} \,\bigl(X_{\infty} - X_j \bigr) \,+ \, \mu \,
\bigl(X_{\infty} - X_j \bigr) \, + \, {1\over2} \, \bigl[ ( X_{j}  - X_{\infty})
\,K \, X_{\infty}  $ }
\setbox22=\hbox{$\displaystyle  \qquad \qquad \qquad \qquad 
+ \, X_{\infty} \, K \,  ( X_{j}  - X_{\infty}) \bigr] \quad \equiv \quad 
\varphi\ib{\Sigma_2} (X_{\infty} - X_j) \,  $ }
\setbox23= \vbox {\halign{#\cr \box20 \cr \box21 \cr \box22 \cr }}
\setbox24= \hbox{ $\vcenter {\box23} $}
\setbox44=\hbox{\noindent (3.49) $\, \, \left\{  \box24 \right. $}
\smallskip \noindent $ \box44 $

\smallskip \noindent 
with

\smallskip \noindent  $ \displaystyle
\Sigma_1 \,\, = \,\, {{1}\over{2 \, \Delta t}} \,I \,+\,  {1\over2} \,  K \, X_j 
 \, + \, M \,$
\smallskip \noindent  $ \displaystyle
\Sigma_2 \,\, = \,\, {{1}\over{2}} \bigl( \mu + {{1}\over{\Delta t}} \bigr) \,-\, 
{{1}\over{2}} \, K \, X_{\infty} \, . \,$

\smallskip \noindent 
The matrix $\, \Sigma_1 \,$ admits a positive definite symmetric part $\, \Sigma_1
\,+\, \Sigma_1^{\displaystyle \rm t} \,$  because it is the case for matrix $\, I .\,$
Moreover, since matrix $\, X_{\infty} - X_j \,$ is symmetric positive by hypothesis,
it is sufficient to establish that the symmetric matrix  $\, \Sigma_2
\,+\, \Sigma_2^{\displaystyle \rm t} \,$  is positive definite and to
apply the proposition 3. This last property is exactly expressed by hypothesis (3.47)
and the first point is proven. 

\bigskip \noindent $\bullet \quad $ 
We consider now the matrix $\, S_\infty \,$ defined by : 

\smallskip \noindent   $   (3.50) \qquad \quad \displaystyle
S_\infty \, \,\,=\,\, \, {{1} \over {2}} (X_{\infty} \, K \,+\, K \,  X_{\infty} )
\,-\, (A^{\rm \displaystyle t} \,+\, A) \,.$

\smallskip \noindent 
The matrix $\,S_\infty\,$ is symmetric and we have the following calculus : 

\smallskip \noindent  $ \displaystyle
S_\infty \,\,= \,\, {1\over2} \Bigl[ (X_{\infty} \, K  \,  X_{\infty} ) \,
X_{\infty}^{-1} \,+\, X_{\infty}^{-1} \,  (X_{\infty} \, K  \,  X_{\infty} ) \Bigr] $

\smallskip     $ \displaystyle \hfill 
\,\,- \,\,  \Bigl[  X_{\infty}^{-1} \,(X_{\infty} \, A) \,+\, (A^{\rm \displaystyle t}
\,X_{\infty}) \,  X_{\infty}^{-1} \Bigr]  \,$

\smallskip \noindent   $   (3.51) \qquad \quad \displaystyle
S_\infty \,\,= \,\,\varphi\ib{\Sigma_3} \bigl( X_{\infty}^{-1} \bigr) \,$

\smallskip \noindent 
with   $ \displaystyle \qquad  \Sigma_3 \,\,= \,\,  {1\over2} \, X_{\infty} \, K  \, 
X_{\infty} \,-\,A^{\rm \displaystyle t} \, X_{\infty} \,. \quad $ Then due to
relation ( 3.41)  we have   $  \displaystyle \quad \Sigma_3 \, + \, \Sigma_3^{\rm
\displaystyle t} \,\,= \,\, Q  \quad > \,\, 0  \quad  $ due to hypothesis (3.46). Then
the propostion 3 joined with (3.51) and the fact that matrix $\, X_{\infty}^{-1} \,$
is symmetric positive definite establish that matrix $\, S_\infty \,$ is symmetric
definite positive.

\bigskip \noindent $\bullet \quad $ 
We  establish that under the same hypothesis (3.48), the matrix $\, Z \, \equiv \,
X_{j+1} -  X_{j} \,$ is positive. We start from the numerical scheme (3.5) 
and  replace matrix $\, Q \,$ by its value obtained from relation (3.41). It comes~:  

\smallskip \noindent   $  \displaystyle
\varphi\ib{\Sigma_1} (Z) \quad \equiv \quad 
{{1}\over{\Delta t}} \, Z \,+\, {1\over2} \, (X_{j} \, K \, Z \,+\, Z \, K \, X_{j} )
\, + \, (M^{\rm \displaystyle t} \, Z \,+\, Z \,M ) \,\,= \,$
\smallskip \noindent   $  \displaystyle \qquad \qquad \,\, 
= \,\, \mu \, X_j \,+ \, Q \,-\, X_{j} \, K \, X_{j} \, - \,   (M^{\rm \displaystyle t}
\, X_{j} \,+\, X_{j} \,M )  \,$
\smallskip \noindent   $  \displaystyle \qquad \qquad \,\, 
= \,\,  A^{\rm \displaystyle t} \, X_{j} \,+\, X_{j} \,A \, + \, Q \, -  X_{j} \, K \,
X_{j} \,  \,$
\smallskip \noindent   $  \displaystyle \qquad \qquad \,\, 
= \,\,  X_{\infty} \, K  \,  X_{\infty} \,-Ê\,  X_{j} \, K \, X_{j} \,- \bigl[
A^{\rm \displaystyle t} \,(X_{\infty} \,-Ê\,  X_{j} ) \,+\, (X_{\infty} \,-Ê\,  X_{j}
) \, A \bigr]   \,$
\smallskip \noindent   $  \displaystyle \qquad \qquad \,\, 
= \,\, {1\over2} \, \Bigl[ X_{\infty} \, K  \,(X_{\infty} \,-Ê\,  X_{j} ) \,+\,
(X_{\infty} \,-Ê\,  X_{j} ) \, K \, X_{\infty} \,+Ê\, X_{j} \, K  \,(X_{\infty}
\,-Ê\,  X_{j} ) \,+\,  \,$
\smallskip \noindent   $  \displaystyle \qquad  \qquad \qquad \,\, 
\,+\, (X_{\infty} \,-Ê\,  X_{j} ) \, K \, X_{j}   \Bigr] \,- \, \bigl[ A^{\rm
\displaystyle t} \, (X_{\infty} \,-Ê\,  X_{j} ) \,+\, (X_{\infty} \,-Ê\,  X_{j} ) \,A
\bigr] \,$ 
\smallskip \noindent   $  \displaystyle \qquad \qquad \,\, 
= \,\,   \varphi\ib{\Sigma_4 } \bigl(X_{\infty} \,-Ê\,  X_{j} \bigr) \quad $ with
$\displaystyle \quad \Sigma_4 \,\, =  \,\,  {{1}\over{2}} \, (K \, X_{\infty} \,+
\, K \, X_j) \,- \, A  \,. \quad $   

\smallskip \noindent  
The matrix   $  \displaystyle  \quad  \Sigma_4 \,+\, \Sigma_4^{\displaystyle \rm t}
\,=\,  {{1}\over{2}} \, (X_j  K \,+ \, K  X_j) \,+ \, S_\infty  \quad  $  is
positive definite due to the previous point  ;  in consequence matrix  $\,
\varphi\ib{\Sigma_4} \bigl(X_{\infty} \,-Ê\,  X_{j} \bigr) \,$ is symmetric positive.
The end of the proof is a consequence of propositon 3 and of the fact that  the
matrix  $ \quad  \displaystyle \Sigma_1 \,\, = \,\, {{1}\over{2}} \, I \,+\, 
{{1}\over{2}} \,K \, X_j + \, M \quad $  has clearly a symmetric part  $\displaystyle
\quad \, \Sigma_1 \,+\, \Sigma_1^{\displaystyle \rm t} \quad $  which is positive
definite.   $ \hfill \square \kern0.1mm $

\bigskip 
\bigskip \noindent  {\bf Proposition  6. $\quad$  Convergence when discrete time 
tends to infinity.  }

 \noindent 
We suppose that the data $\,K,\, A \,,\, Q \, $ of Riccati equation (1.20) and
parameters $\, \mu \,$ and $\, \Delta t \,$ of harmonic scheme (3.5)  satisfy the 
conditions (3.1), (3.46) and (3.47). Let  $\,  X_{j} \,$ and $\,  X_{\infty}
\,$ be the solution of scheme (3.5) at discrete time $\, j\, \Delta t \,$ and 
the symmetric definite positive matrix solution of the so-called algebraic
Riccati equation (3.41). Then $\,  X_{j} \,$ tends to $\,  X_{\infty} \,$ when
$\,j \,$ tends to infinity : 
\smallskip \noindent   $   (3.52) \qquad \quad \displaystyle
 X_{j} \,\, \longrightarrow \,\,  X_{\infty} \,. \,$

\bigskip \noindent  {\bf Proof of  proposition 6.}

\noindent $\bullet \quad $ 
Let $\, {\cal E}_k \,$ be the Grassmannian manifold composed by all the   linear
subspaces of space $\, \R^n \,$ with dimension exactly equal to $\,k \,$ : 

\smallskip \noindent   $   (3.53) \qquad \quad \displaystyle
{\cal E}_k \,\, = \,\, \bigl\{ \, W, \, W \,$ subspace of $ \,\, \R^n, \, {\rm dim}
\, W \,=\, k \, \bigr\} \,. \,$

\smallskip \noindent 
Then we have the classical characterization of the $\, k^{\rm o}$ eigenvalue of
symmetric matrix $\, A \,$ with the so-called inf-sup condition (see {\it e.g.} Lascaux and
Th\'eodor  [LT86])  : 

\smallskip \noindent   $   (3.54) \qquad \quad \displaystyle
\mu_k \,\, =\,\, \inf_{\displaystyle W \!\! \in  \!{\cal E}_k} \,\,\, 
\sup_{\displaystyle v  \!\! \in  \! W} \,\, {{(Av \,,\, v)}\over{(v\,,\, v)}} \, . \,$

\smallskip \noindent
Let $\, \lambda^k_{j} \,$ be the $\, k^{\rm o}$ eigenvalue of matrix $\, X_j ,\,$ 
$\, \lambda^k_{\infty} \,$  the $\, k^{\rm o}$ eigenvalue $\,\, \mu_k \,\,(\mu_1
\,\leq\, \mu_2 \,\leq\, \cdots \, \,\leq\, \mu_n ) \,\, $ of matrix $\, X_{\infty}
\,$  and $\,W\,$ a fixed subspace of $\, \R^n \,$ of dimension exactly equal
to $ \, k .\,$  We deduce from inequality  (3.48)  of  Proposition 5 :  

\smallskip \noindent   $   \displaystyle
{{(X_j \, v \,,\, v)}\over{(v\,,\, v)}} \,\,  \leq  \,\,   {{(X_{j+1} \, v \,,\,
v)}\over {(v\,,\, v)}}  \,\,   \leq  \,\,    {{(X_{\infty } \, v \,,\, v)}\over
{(v\,,\, v)}}  \,,\qquad \forall \, v \neq 0 \,,\,\, v  \in W . \qquad $ 

\smallskip \noindent 
Then when subspace $\,W\,$ is arbitrarily given in set $\, {\cal E}_k \,$ we have : 

\smallskip \noindent   $   \displaystyle
\sup_{\displaystyle v  \!\! \in  \! W} \, \biggl[ {{(X_j \, v \,,\, v)}\over{(v\,,\,
v)}} \, \biggr]   \,\,   \leq  \,\,   \sup_{\displaystyle v  \!\! \in  \! W} \,
\biggl[  {{(X_{j+1} \, v \,,\, v)}\over {(v\,,\, v)}}  \, \biggr]  
\,\,   \leq  \,\,   \sup_{\displaystyle v  \!\! \in  \! W} \,
\biggl[  {{(X_{\infty} \, v \,,\, v)}\over {(v\,,\, v)}}  \, \biggr]  \,,\,\, 
\forall \, W  \! \in \!  {\cal E}_k \, \qquad $

\smallskip \noindent
and taking the infimum bound of previous line as subspace $\,W \,$ belongs to
Grassmannian manifold $\, {\cal E}_k \,$ we deduce, thanks to (3.54) 

\smallskip \noindent   $   (3.55) \qquad \quad \displaystyle
\lambda^k_{j} \,\, \leq \,\, \lambda^k_{j+1} \,\, \leq \,\, \lambda^k_{\infty}  
\,\, . \,$

\smallskip \noindent 
By monotonicity, eigenvalue $\, \lambda^k_{j} \,$ is converging towards some
scaler $\, \mu^k \,$ as $\, j \,$ tends to infinity : 

\smallskip \noindent   $   (3.56) \qquad \quad \displaystyle
\lambda^k_{j} \,\, \longrightarrow \,\,   \mu^k  \,\quad {\rm as} \,\,\, j \,\,
\longrightarrow \,\, \infty \,,  \qquad 1 \leq k \leq n  \,. \,$

\smallskip \noindent $\bullet \quad $ 
Consider now the unitary eigenvector $\, v^k_{j} \,$ of matrix $\, X_j \,$ associated
with eigenvalue $\, \lambda^k_{j} \,$ : 

\smallskip \noindent   $   (3.57) \qquad \quad \displaystyle
X_j \,  v^k_{j} \,\,= \,\, \lambda^k_{j} \, v^k_{j} \,\,, \quad \parallel v^k_{j}
\parallel \,=\, 1 \,\,, \quad 1 \leq k \leq n \,,\, \,\, j \geq 0 \,. \,$

\smallskip \noindent
Because $\, X_j \,$ is a symmetric matrix, the family of eigenvectors $\, \bigl( 
v^k_{j} \bigr)_{1 \leq k \leq n} \,$ is orthonormal and defines an orthogonal
operator  $\, \rho_j \,$ of space $\, \R^n \,$ defined as acting on the canonical
basis $\, \bigl(  e_{j} \bigr)_{1 \leq k \leq n} \,$ of space $\, \R^n \,$  by the
conditions 

\smallskip \noindent   $   (3.58) \qquad \quad \displaystyle
\rho_j \,{\scriptstyle \bullet}\,  e_k \,\,= \,\, v^k_{j} \,\,,\quad 1 \leq k \leq n
\,,\, \,\, j \geq 0 \,.\,$

\smallskip \noindent  
Rotation $\, \rho_j \,$ belongs to the compact group $\, O(n) \,$ of orthogonal
linear transformations  of space $\, \R^n \,$ . Then after extraction of a
convergent  subsequence $\, \rho'_j \,$ of the initial sequence $\,  \bigl( \rho_j
\bigr)_{j \geq 0}  \,$ we know that there exists an orthogonal mapping  $\,
\rho_{\infty} \in  O(n) \,$ such that 

\smallskip \noindent   $   (3.59) \qquad \quad \displaystyle
 \rho'_j  \,\, \longrightarrow \,\,   \rho_{\infty} \qquad   {\rm as } \,\, j \,\,
\longrightarrow \,\, +\infty \,.$ 

\smallskip \noindent 
We introduce the family $\, \bigl(   w^k_{\infty} \bigr)_{1 \leq k \leq n} \,$ of
vectors in $\, \R^n \,$ by the conditions 

\smallskip \noindent   $   (3.60) \qquad \quad \displaystyle
w^k_{\infty} \,\,= \,\, \rho_{\infty}  \,{\scriptstyle \bullet}\,  e_k \,\,,\quad 1
\leq k \leq n \,\,. \, $

\smallskip \noindent 
It constitutes an orthogonal basis of space $\, \R^n \,$ and for each integer $\,k
\,(1 \leq k \leq n),  $  the extracted sequence of vectors $ \, v'^{k}_{j} \,$ is
converging towards vector $\, w^k_{\infty} \,$ : 

\smallskip \noindent   $   (3.61) \qquad \quad \displaystyle
 v'^{k}_{j}  \,\, \longrightarrow \,\,   w^k_{\infty} \qquad ,\,\,  1 \leq k \leq n \,,\,
\,\, j \longrightarrow \,\, +\infty \,.$ 

 \smallskip   \noindent  $\bullet \quad$
We introduce the symmetric positive definite  operator $\, Y_{\infty} \,$  by the
conditions 

\smallskip \noindent   $   (3.62) \qquad \quad \displaystyle
Y_{\infty} \,{\scriptstyle \bullet}\,  w^k_{\infty} \,\,= \,\, \mu^k \,  w^k_{\infty}
\,\qquad  ,\,\, 1 \leq k \leq n \,.\,  $ 

\smallskip \noindent
We study now the convergence of the subsequence of matrices $\, X'_j \,$ towards
$\, Y_{\infty} .\, $ We first remark that the sequence of matrices $\, \bigl(  X_j
\bigr)_{j \geq 0}  \,$ is bounded in space  $\,{\cal{S}}_n(\R) \, $~: 

\smallskip \noindent   $   (3.63) \qquad \quad \displaystyle
\parallel  X_j \parallel \,\,\,  \leq \,\, \lambda^n_{\infty} \, \quad  ,\,\,\forall \, j
\in \N \, $ 

\smallskip \noindent 
and we have also the following set of identities : 

\smallskip \noindent   $  \displaystyle
\bigl( X_j \,- \, Y_{\infty} \bigr) \,   w^k_{\infty} \,\, = \,\, X_j \, (
w^k_{\infty} \,-\,  v^{k}_{j} ) \,\, - \,\, \lambda^k_{\infty} \,  w^k_{\infty}
\,\,+\,\, \lambda^k_{j} \,  v^{k}_{j} \,$ 
\smallskip \noindent   $  \displaystyle
\qquad \qquad \qquad \quad \,\, \, = \,\, X_j \, ( w^k_{\infty} \,-\,  v^{k}_{j} ) \,\,
- \,\, \lambda^k_{\infty} \, ( w^k_{\infty} \,-\,  v^{k}_{j}) \,\,+\,\, (\lambda^k_{j}
\,-\, \lambda^k_{\infty} ) \,  v^{k}_{j} \,. \, $

\smallskip \noindent 
So we deduce, due also to  (3.63) and (3.57) : 

\smallskip \noindent   $   (3.64) \,\,  \displaystyle
\parallel \bigl( X_j \,- \, Y_{\infty} \bigr) \,   w^k_{\infty} \parallel \, \leq
\, \lambda^n_{\infty} \, \parallel  w^k_{\infty} \,-\,  v^{k}_{j} \parallel \, +
\,  \lambda^k_{\infty}  \, \parallel w^k_{\infty} \,-\,  v^{k}_{j} \parallel \,+
\, \mid \lambda^k_{j} \,-\, \lambda^k_{\infty} \mid \, $

\smallskip \noindent 
and according of the convergence  (3.61) of subsequences $\,  v^{k}_{j} \,$ 
and  (3.56) of sequences $\, \lambda^k_{j} \,$ as index $\,j \,$ tends to infinity,
we get from estimation (3.64) : 

\smallskip \noindent   $   (3.65) \qquad \quad \displaystyle
 X'_j \,- \, Y_{\infty}  \,\, \longrightarrow \,\,0 \,\, , \qquad j  \longrightarrow
\,+\infty \,. \, $

\smallskip \noindent  $\bullet \quad$
Due to the definition (3.5) of the numerical scheme, sub-sequence $\,  X'_j  \,$
necessarily converges to the unique positive definite matrix $\,  X_{\infty} \, $ of
the algebraic Riccati equation and in consequence we have necessarily 

\smallskip \noindent   $   (3.66) \qquad \quad \displaystyle
Y_{\infty} \,\,= \,\, X_{\infty} \,. \,$ 

\smallskip \noindent  
We deduce that for any arbitrary subsequence of the family $\, \bigl(  X_j \bigr)_{j
\geq 0}  \,$ there exists an extracted sub-subsequence converging towards $\,
X_{\infty} \,$ as $\,j \,$ tends to infinity. Then the entire family   $\, \bigl(  X_j 
\bigr)_{j \geq 0}  \,$ converges towards  $\, X_{\infty}\,$ as $\,j \,$ tends to
infinity and the property is established.  $ \hfill \square \kern0.1mm $

%%%%%%%%%%%%%%%%%%%%%%%%%%%%%%%%%%%%%%%%%%%%%%%%%%%%%%%%%%%%%%%%%%%%%%%%%%%%%%%%%%%%%%
%%%%%%%%%%%%%%%%%%%%%%%%%%%%%%%%%%%%%%%%%%%%%%%%%%%%%%%%%%%%%%%%%%%%%%%%%%%%%%%%%%%%%%
% \vfill \eject 
\bigskip \bigskip \noindent  {\smcaps 4) $ \quad$ Numerical experiments.}

\smallskip \noindent  {\bf 4.1) $ \quad$ Square root function.}
\smallskip \noindent $\bullet \quad$
The first example studied is the resolution of the equation :

\smallskip \noindent  (4.1) $  \qquad \quad  \displaystyle
{{{\rm d}X} \over {{\rm d}t}} \,+\, X^2 \,-\, Q \,\,=\,\, 0 , \quad X(0)
\,\,=\,\, 0 \,$

\smallskip \noindent
with $ \, n = 2 ,\,\, A = 0 ,\,\,  K = I \,\,$ and  matrix $\, Q \,$ equal to

\smallskip \noindent  (4.2) $  \qquad \quad \displaystyle
Q \,\,=\,\, {1\over2} \,\, \pmatrix{1 & -1 \cr 1 & 1\cr } \, \, \pmatrix{1
& 0 \cr 0 & 100  \cr } \, \,\, \pmatrix{1 & 1 \cr -1 & 1\cr }  \,.  $

\smallskip \noindent $\bullet \quad$
We have tested our numerical scheme for fixed value  $\,
\Delta t = 1/100 \,$ and   different values of parameter $\, \mu \, : \,$
$\,  \mu = 0.1 , \, 10^{-6} ,\,  10^{+6} . \,$ For small values of parameter $ \, \mu
,\,$ the behaviour of the scheme does not change between $\,   \mu = 0.1 \,$ and
$\,   \mu = 10^{-6} . \,$ Figures 1 to 4 show the evolution with time of the
eigenvalues of matrix $\, X_{j} \,$ and the convergence is achieved to the square root
of matrix $\, Q .\,$ For large value of parameter $\, \mu  \, \, ( \mu =  10^{+6} ),
\,$ we loose completely consistency of the scheme (see figures 5 and 6).

%%%%%%%%%%%%%%%%%%%%%%%%%%%%%%%%%%%%%%%%%%%%%%%%%%%%%%%%%%%%%%%%%%%%%%%%%%%%%%%%%%%%%%%
 \bigskip \noindent 
  { \epsfysize=4,0cm    \epsfbox   {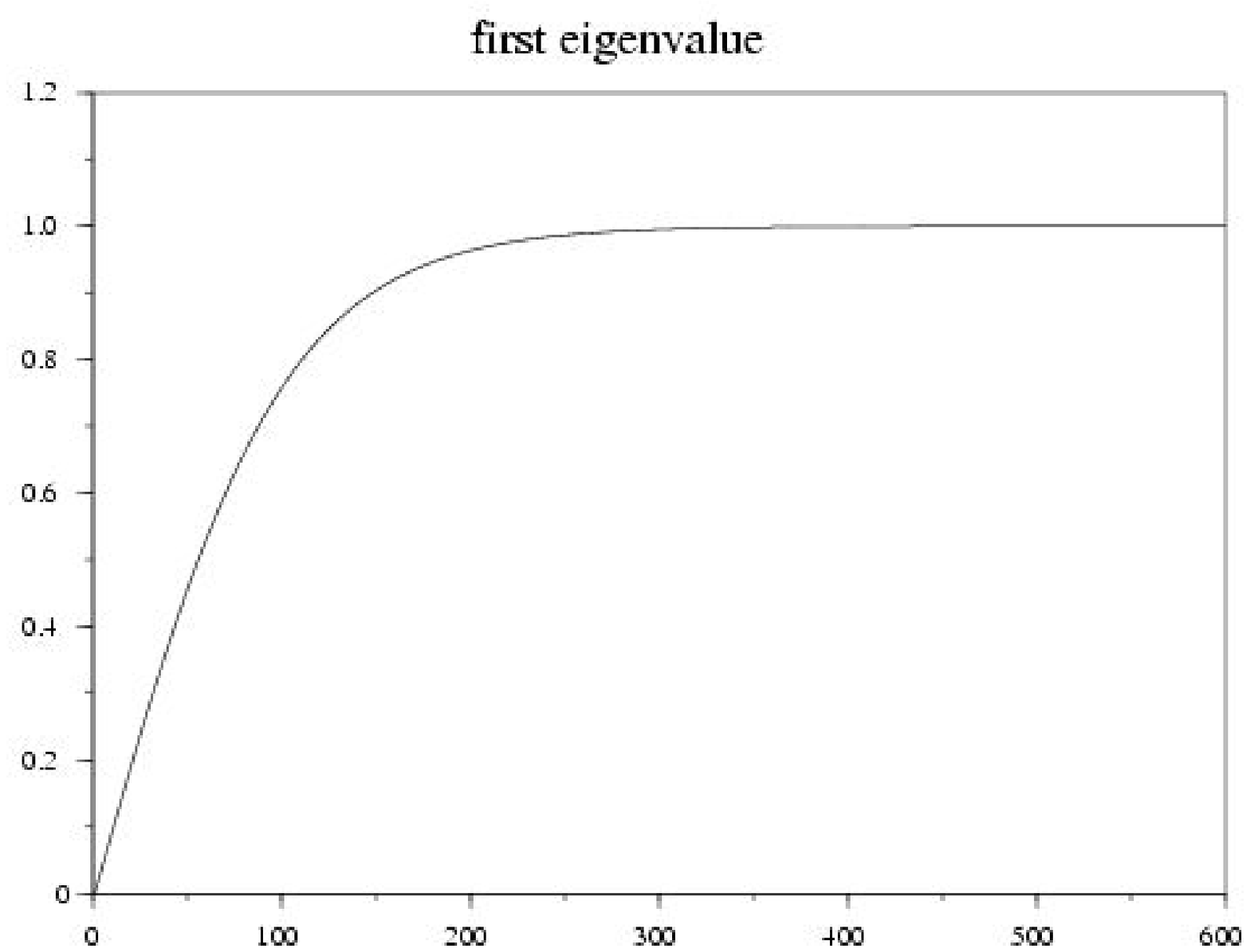 } 
  { \epsfysize=4,0cm     \epsfbox  {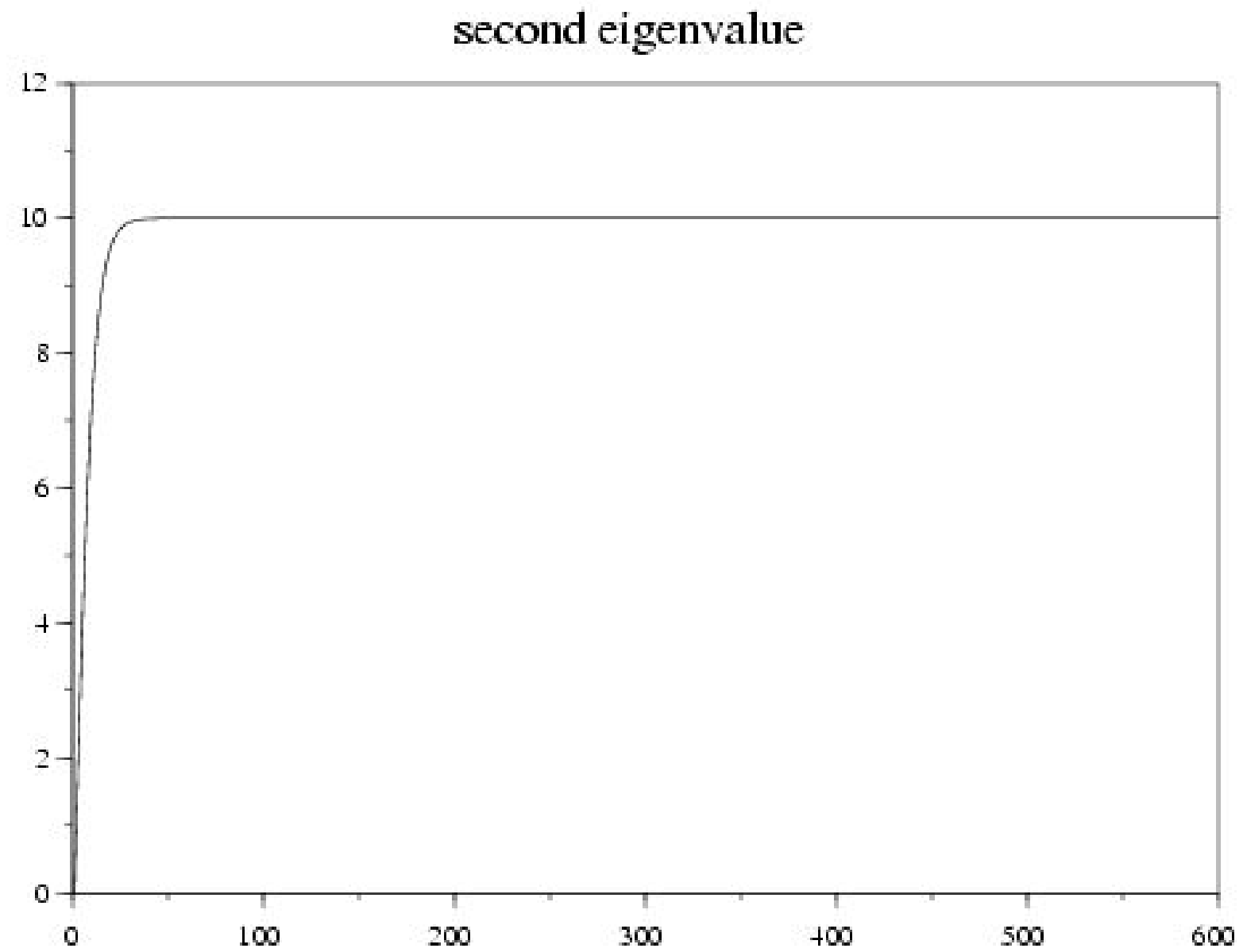 } } 
 \smallskip  

 \centerline  { { \bf Figures 1 and 2.} Square root function test. }
 \centerline  { Two first eigenvalues of numerical solution ($\mu \,=\,0.1$). }
 \smallskip  \smallskip
%%%%%%%%%%%%%%%%%%%%%%%%%%%%%%%%%%%%%%%%%%%%%%%%%%%%%%%%%%%%%%%%%%%%%%%%%%%%%%%%%%%%%%%

%%%%%%%%%%%%%%%%%%%%%%%%%%%%%%%%%%%%%%%%%%%%%%%%%%%%%%%%%%%%%%%%%%%%%%%%%%%%%%%%%%%%%%%
 \bigskip \noindent 
  { \epsfysize=4,0cm    \epsfbox   {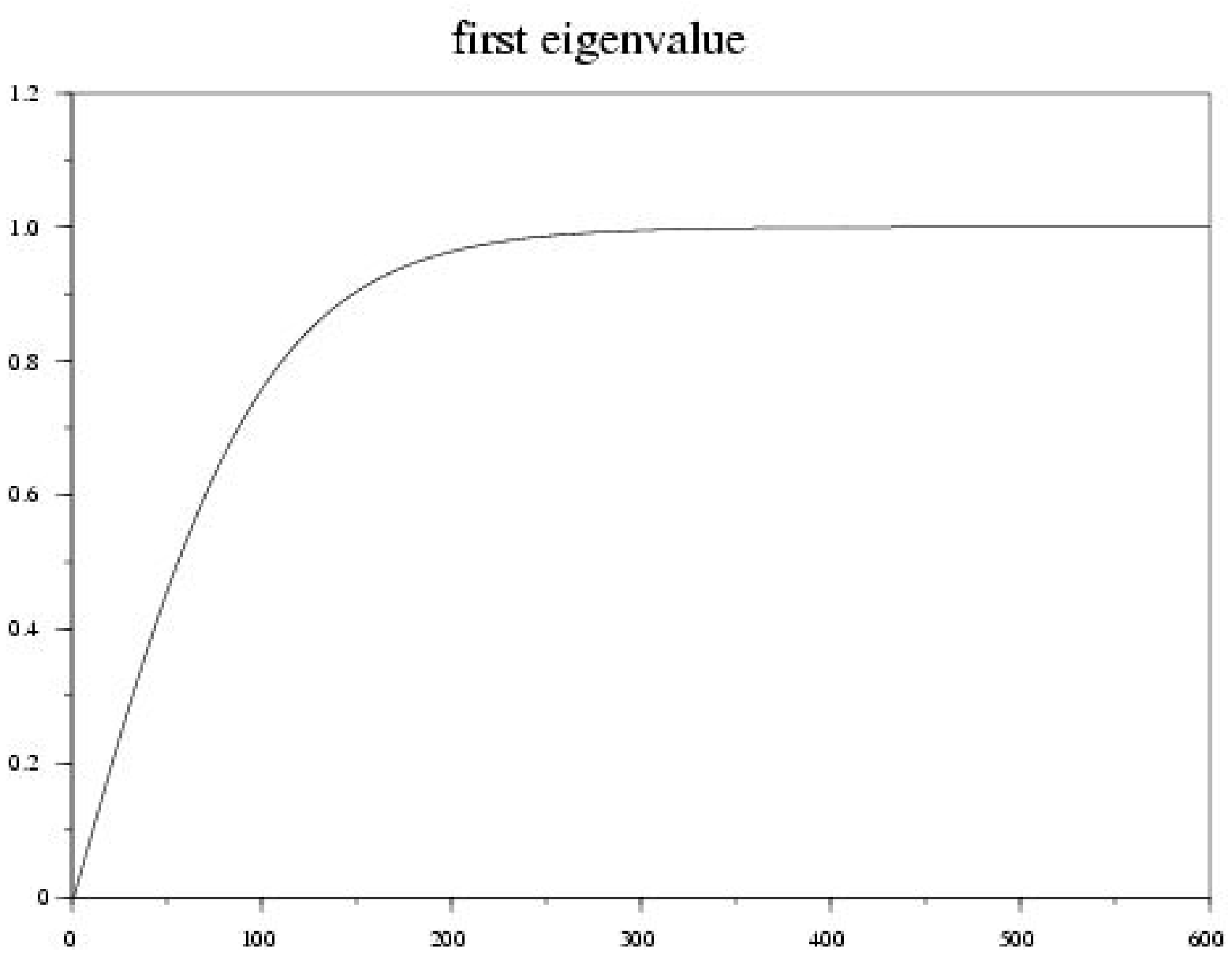 } 
  { \epsfysize=4,0cm     \epsfbox  {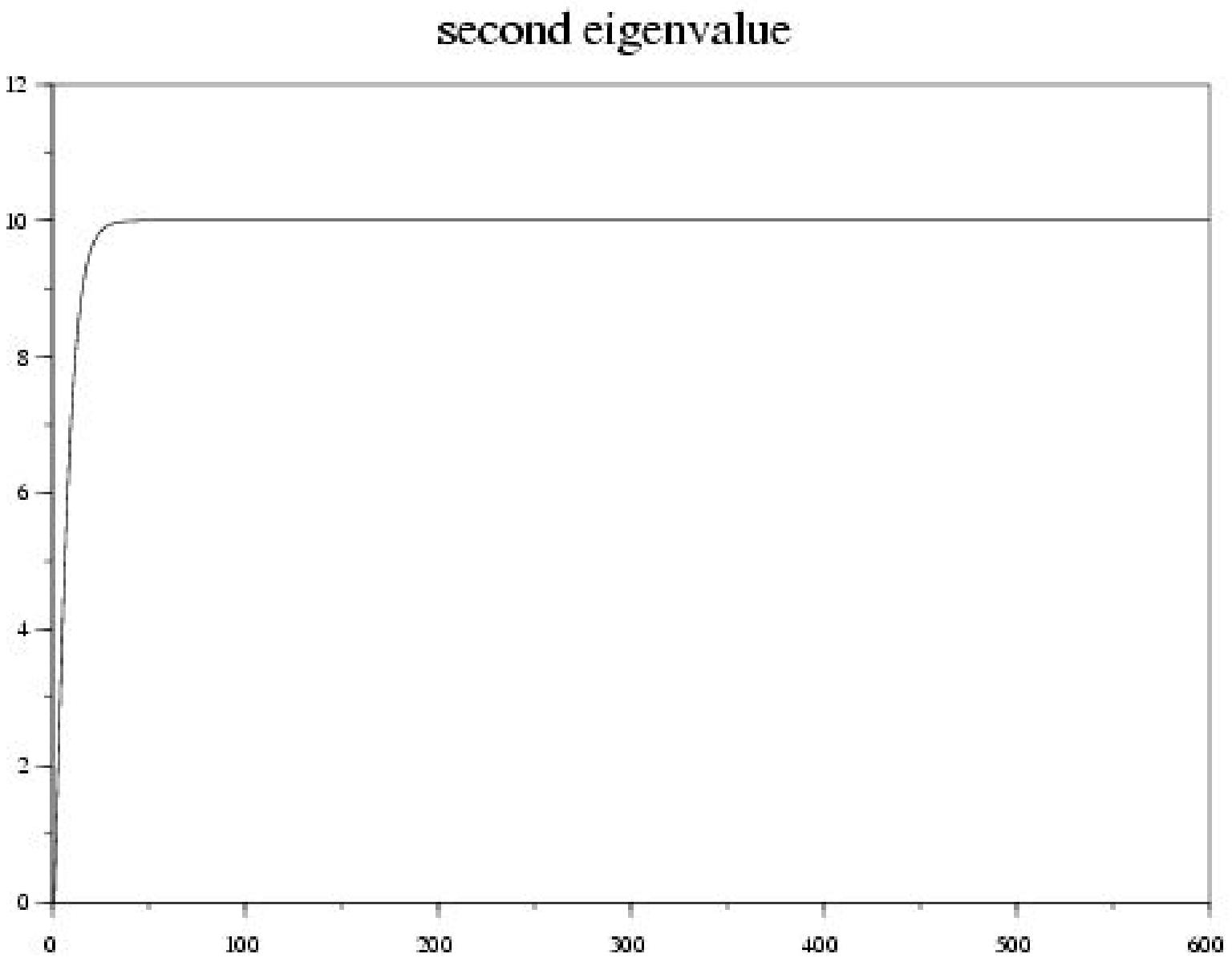 } } 
 \smallskip  

 \centerline  { { \bf Figures 3 and 4.} Square root function test. }
 \centerline  { Two first eigenvalues of numerical solution ($\mu \,=\,10^{-6}$). }
 \smallskip  \smallskip
%%%%%%%%%%%%%%%%%%%%%%%%%%%%%%%%%%%%%%%%%%%%%%%%%%%%%%%%%%%%%%%%%%%%%%%%%%%%%%%%%%%%%%%

%%%%%%%%%%%%%%%%%%%%%%%%%%%%%%%%%%%%%%%%%%%%%%%%%%%%%%%%%%%%%%%%%%%%%%%%%%%%%%%%%%%%%%%
 \bigskip \noindent 
  { \epsfysize=4,0cm    \epsfbox   {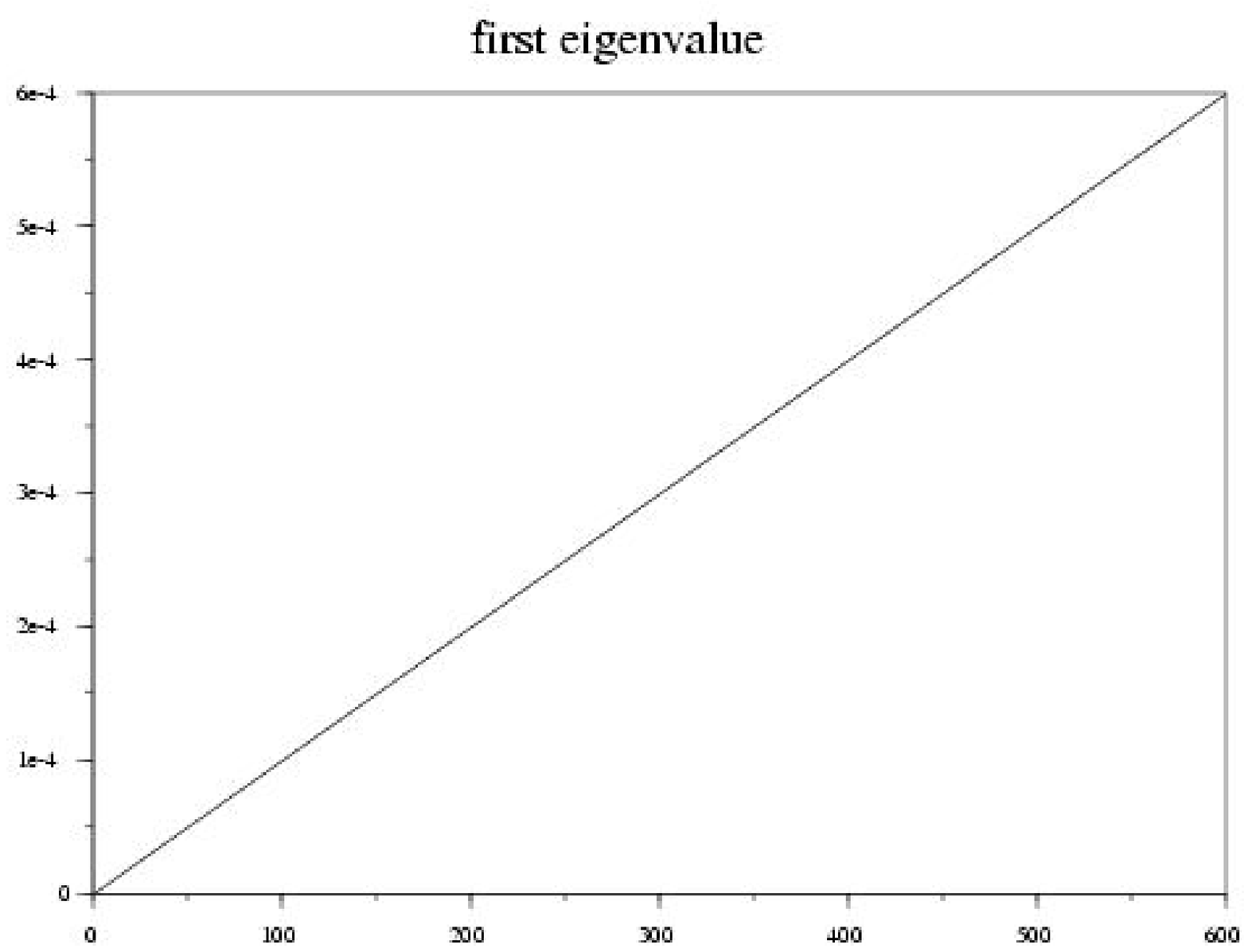 } 
  { \epsfysize=4,0cm     \epsfbox  {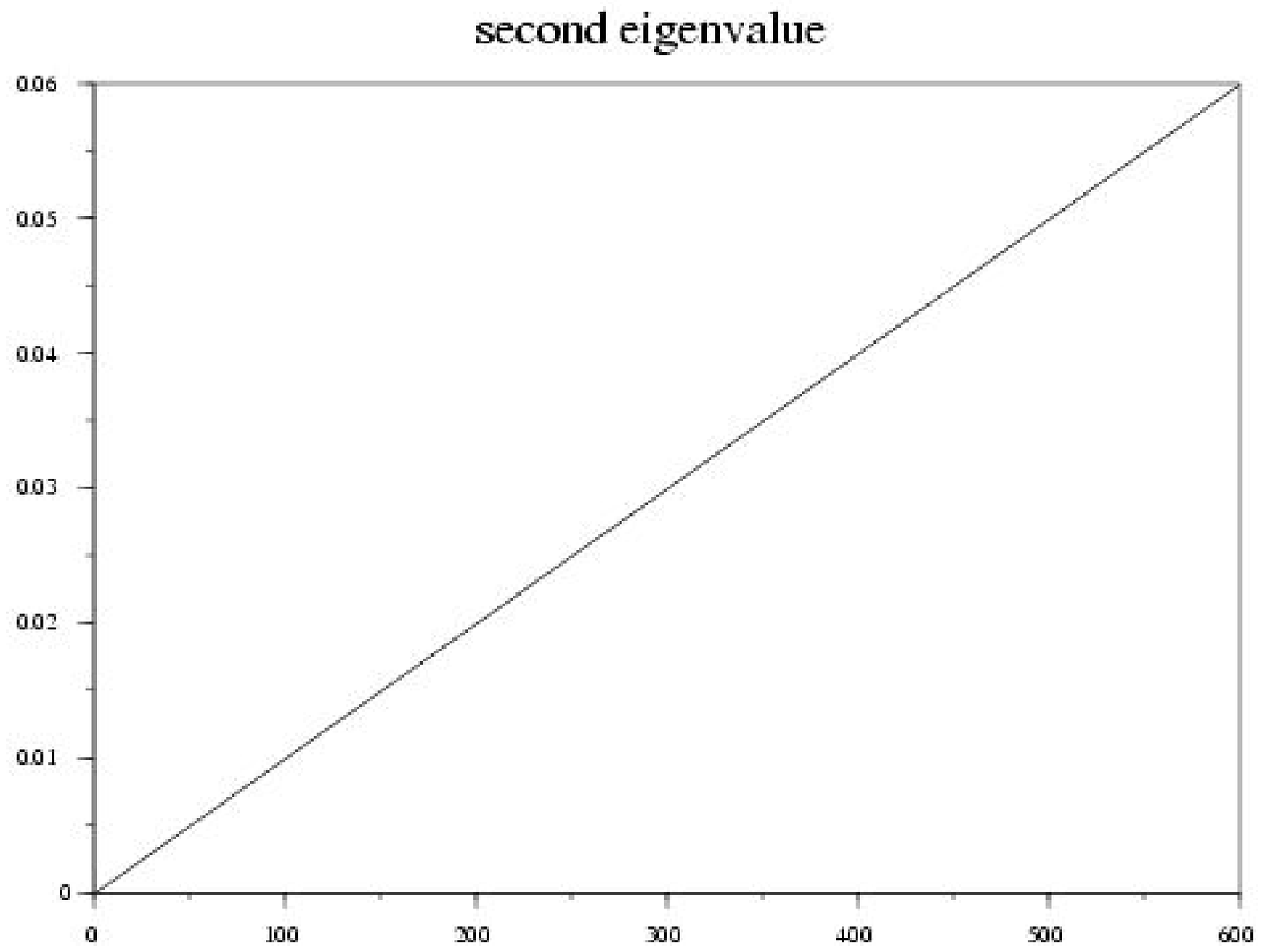 } } 

 \centerline  { { \bf Figures 5 and 6.} Square root function test. }
 \centerline  { Two first eigenvalues of numerical solution ($\mu \,=\,10^{+6}$). }
%  \smallskip  \smallskip
%%%%%%%%%%%%%%%%%%%%%%%%%%%%%%%%%%%%%%%%%%%%%%%%%%%%%%%%%%%%%%%%%%%%%%%%%%%%%%%%%%%%%%%

\bigskip \noindent  {\bf 4.2) $ \quad$ Harmonic oscillator.}

\smallskip  \noindent $\bullet \quad$
The second exemple is the classical harmonic oscillator. Dynamical system
$\, y(t) \,$ is governed by the second order differential equation with  command $\,
v(t) \,$ :

\smallskip \noindent  (4.3) $  \qquad \quad  \displaystyle
{{{\rm d}^2 y(t)} \over {{\rm d}t^2}} \,+ \, 2 \, \delta \,
{{{\rm d} y(t)} \over {{\rm d}t}} \, + \, \omega^2 \, y(t) \,\, = \,\, b \,
v(t) \,. $

\smallskip \noindent
This equation is written as a first order system of differential equations :

\smallskip \noindent  (4.4) $  \qquad \quad  \displaystyle
Y \,= \, \pmatrix {y(t) \cr  {\displaystyle{{\rm d} y(t)} \over
{\displaystyle {\rm d}t}} \cr } \,, \quad  { {{ \rm d} Y} \over { {\rm d}t}} \,\,= \,\,
\pmatrix {0 & 1 \cr - \omega^2 & -2 \, \delta  } \,Y(t) \,\, + \,\, 
\pmatrix{0 \cr b \, v(t)\,} \,.  $

\smallskip \noindent
The parameters $\, R, \, Q ,\, \omega ,\, \delta \, $ and $\, b \, $  of the ordinary
differential equation (4.4) and the cost function (1.6) are given by : 

\smallskip \noindent  (4.5) $  \qquad \quad  \displaystyle
  R \,=\, \pmatrix{\alpha & 0 \cr 0 & \alpha } \,\,, \,\, 
Q \,=\, \pmatrix{ {1\over2} & 0 \cr 0 &{1\over2} }  \,\,, \,\,  
\omega \,=\, \sqrt{250} \,\,, \,\, \delta \,=\, 0  \,\,, \,\, b \,=\, 1   \,. \, $

%%%%%%%%%%%%%%%%%%%%%%%%%%%%%%%%%%%%%%%%%%%%%%%%%%%%%%%%%%%%%%%%%%%%%%%%%%%%%%%%%%%%%%%
\bigskip 
\noindent 
  { \epsfysize=4,0cm    \epsfbox   {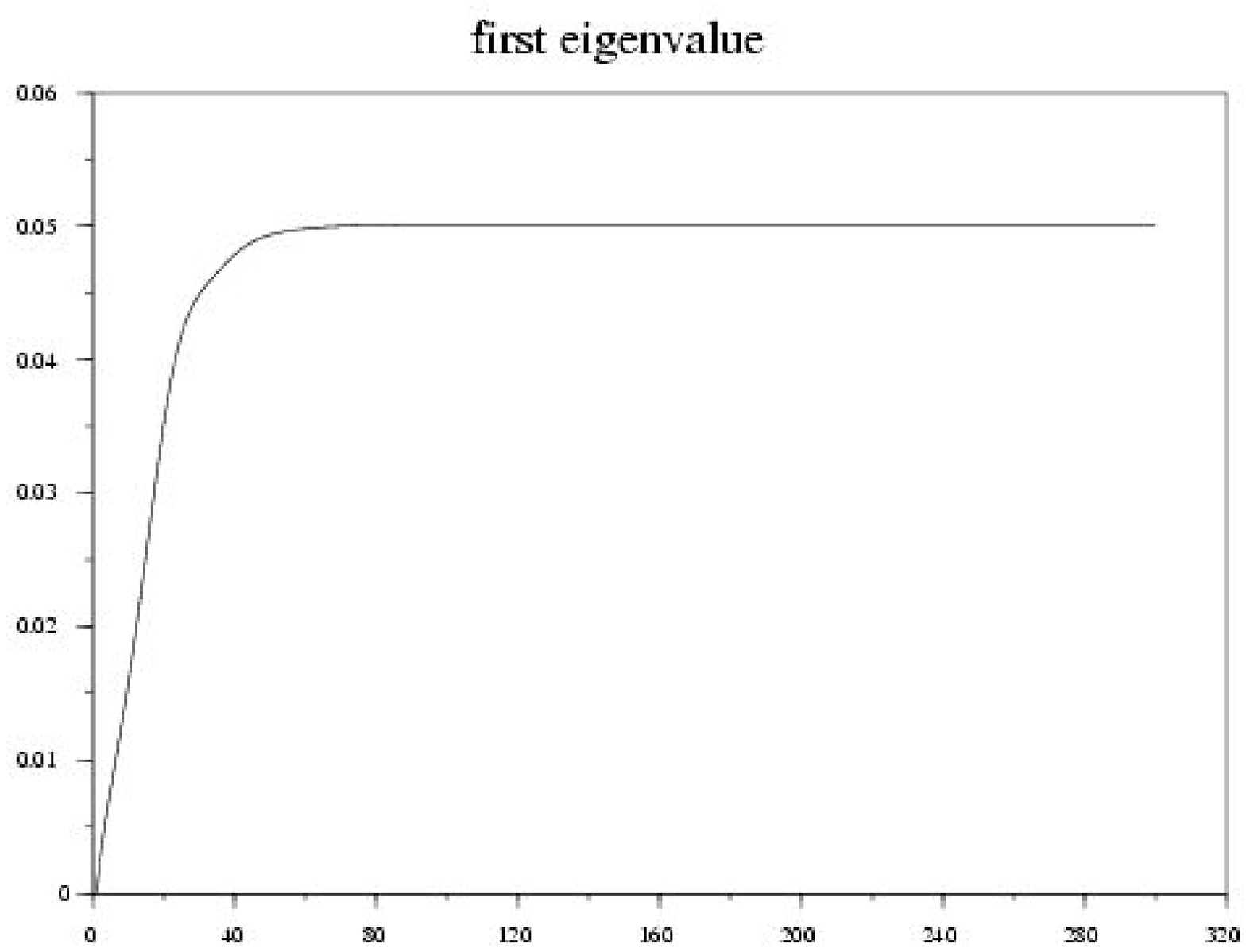 } 
  { \epsfysize=4,0cm     \epsfbox  {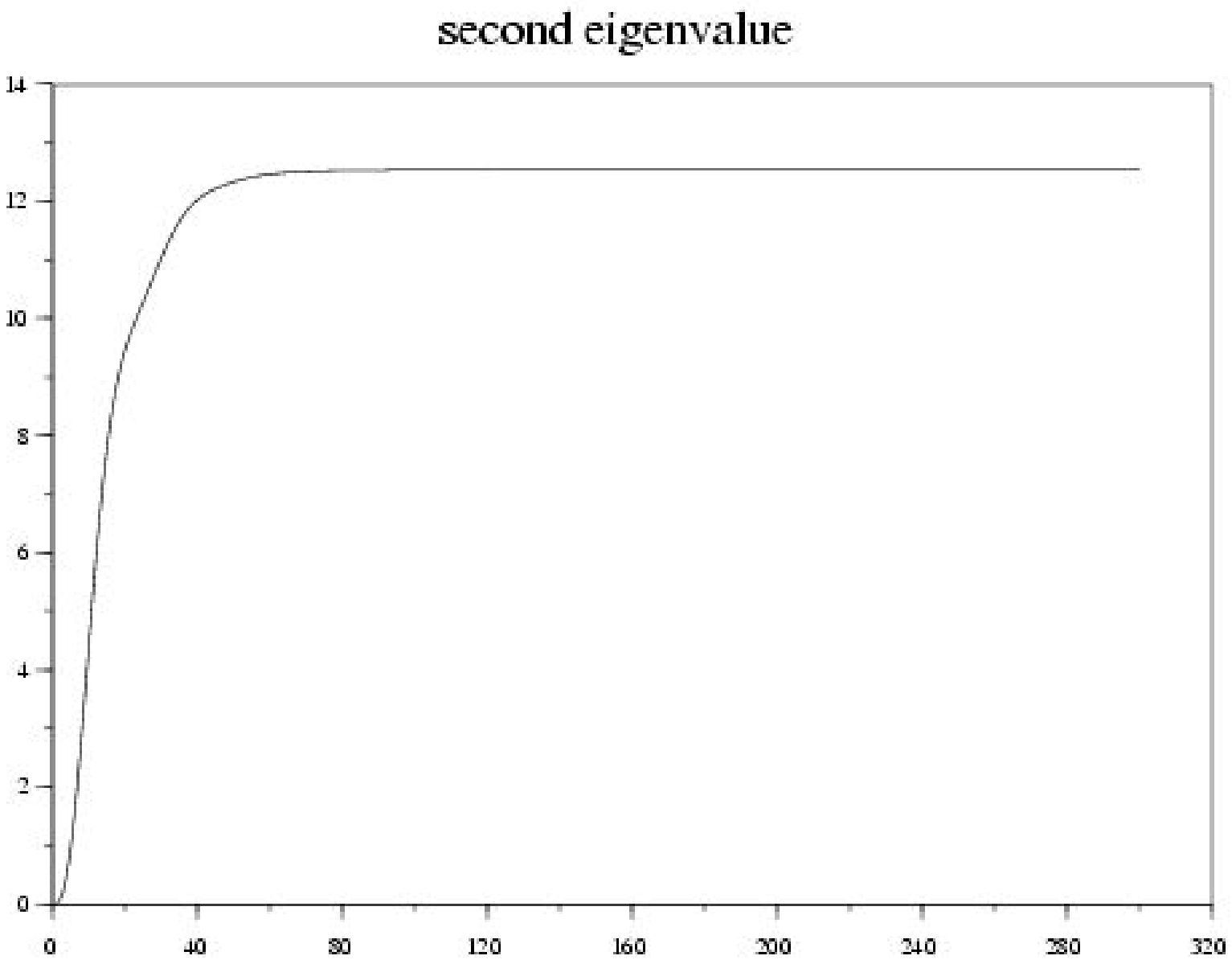 } } 

 \centerline  { { \bf Figures 7 and 8.}  Harmonic oscillator. }
 \centerline  { Two first eigenvalues of numerical solution (
$\mu \,=\, 0.1 , \,\, \alpha \,=\,  0.01 , \,\,  \Delta t \,=\,  0.01 $). }
  \smallskip  \smallskip
%%%%%%%%%%%%%%%%%%%%%%%%%%%%%%%%%%%%%%%%%%%%%%%%%%%%%%%%%%%%%%%%%%%%%%%%%%%%%%%%%%%%%%%

\smallskip \noindent $\bullet \quad$
In this case, we have tested the stability of the scheme for fixed value of
parameter $\, \mu \,( \mu = 0.1) \,$ and  different values of time step $\, \Delta t
.\,$  We have chosen three sets of parameters : $ \, \alpha = \Delta t = 1/100 $
(reference experiment, figures 7 and 8),  $ \, \alpha = 10^{-6} \,, \,  \Delta t =
1/100 $ (very small value for $\alpha$ , figures 9 and 10)  and  $ \, \alpha = 1/100
\,, \,  \Delta t = 100 $ (too large value for time step,  figures 11 and 12). Note
that for the last set of parameters, classical explicit schemes fail to give any
answer.  As in previous test case, we have represented the two eigenvalues of discrete
matrix solution $\, X_j \,$  as time is increasing. On reference experiment (figures~7
and 8), we  have convergence of the solution to the solution of algebraic Riccati
equation. If control parameter $\, \alpha \,$ is chosen too small, the
first eigenvalue of Riccati matrix oscillates during the first time steps but reach
finally the correct values of limit matrix, the solution of  algebraic Riccati
equation. If time step is too large, we still have stability but we loose also
monotonicity. Nevertheless, convergence is achieved as in previous case. 

%%%%%%%%%%%%%%%%%%%%%%%%%%%%%%%%%%%%%%%%%%%%%%%%%%%%%%%%%%%%%%%%%%%%%%%%%%%%%%%%%%%%%%%
  \bigskip \noindent
% \smallskip \noindent 
  { \epsfysize=4,0cm    \epsfbox   {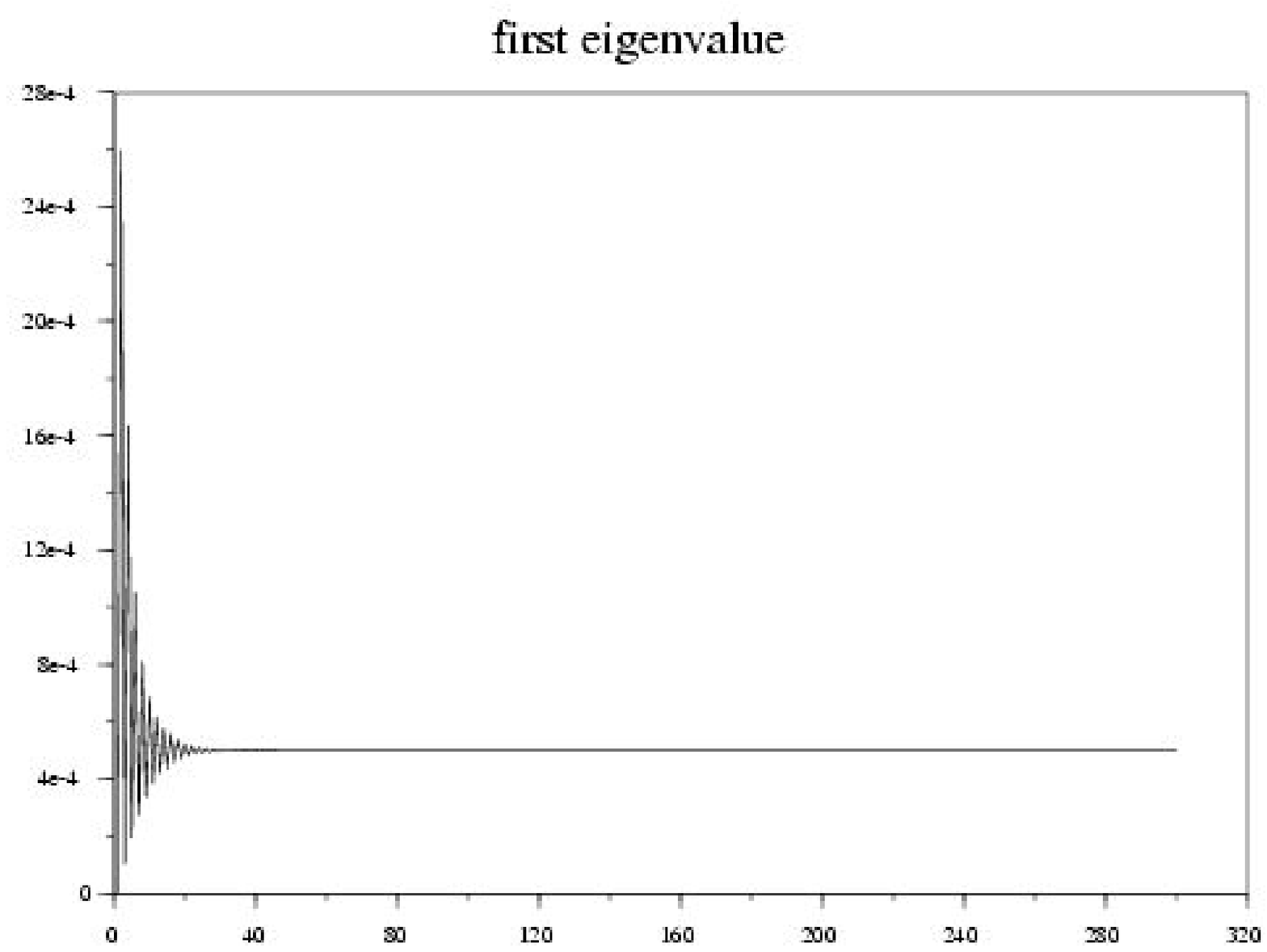 } 
  { \epsfysize=4,0cm     \epsfbox  {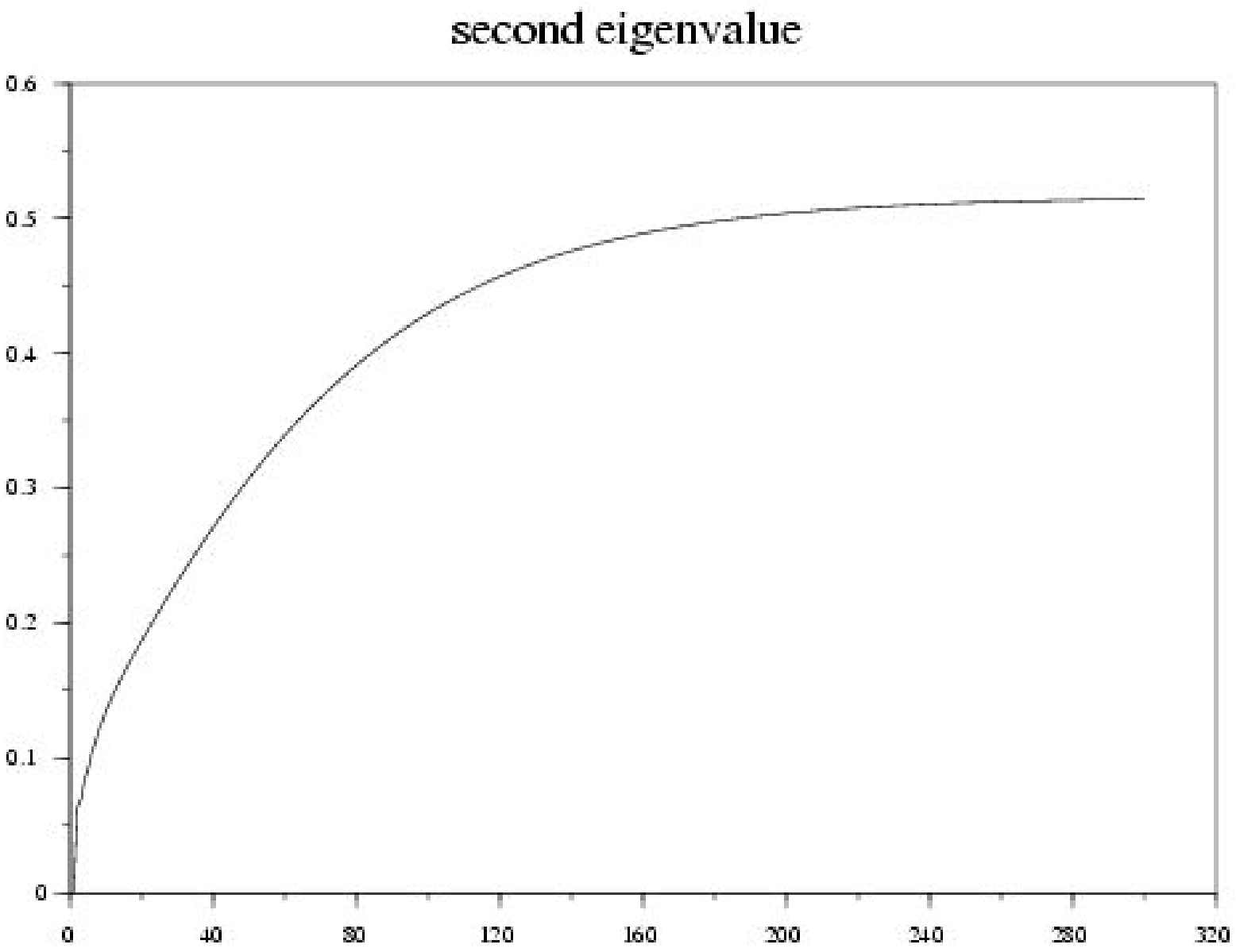 } } 

 \centerline  { { \bf Figures 9 and 10.}  Harmonic oscillator. }
 \centerline  { Two first eigenvalues of numerical solution (
$\mu \,=\, 0.1 , \,\, \alpha \,=\, 10^{-6} , \,\,  \Delta t \,=\,  0.01 $). }
%  \smallskip  \smallskip
%%%%%%%%%%%%%%%%%%%%%%%%%%%%%%%%%%%%%%%%%%%%%%%%%%%%%%%%%%%%%%%%%%%%%%%%%%%%%%%%%%%%%%%

%%%%%%%%%%%%%%%%%%%%%%%%%%%%%%%%%%%%%%%%%%%%%%%%%%%%%%%%%%%%%%%%%%%%%%%%%%%%%%%%%%%%%%%
  \bigskip \noindent
% \smallskip \noindent  
  { \epsfysize=4,0cm    \epsfbox   {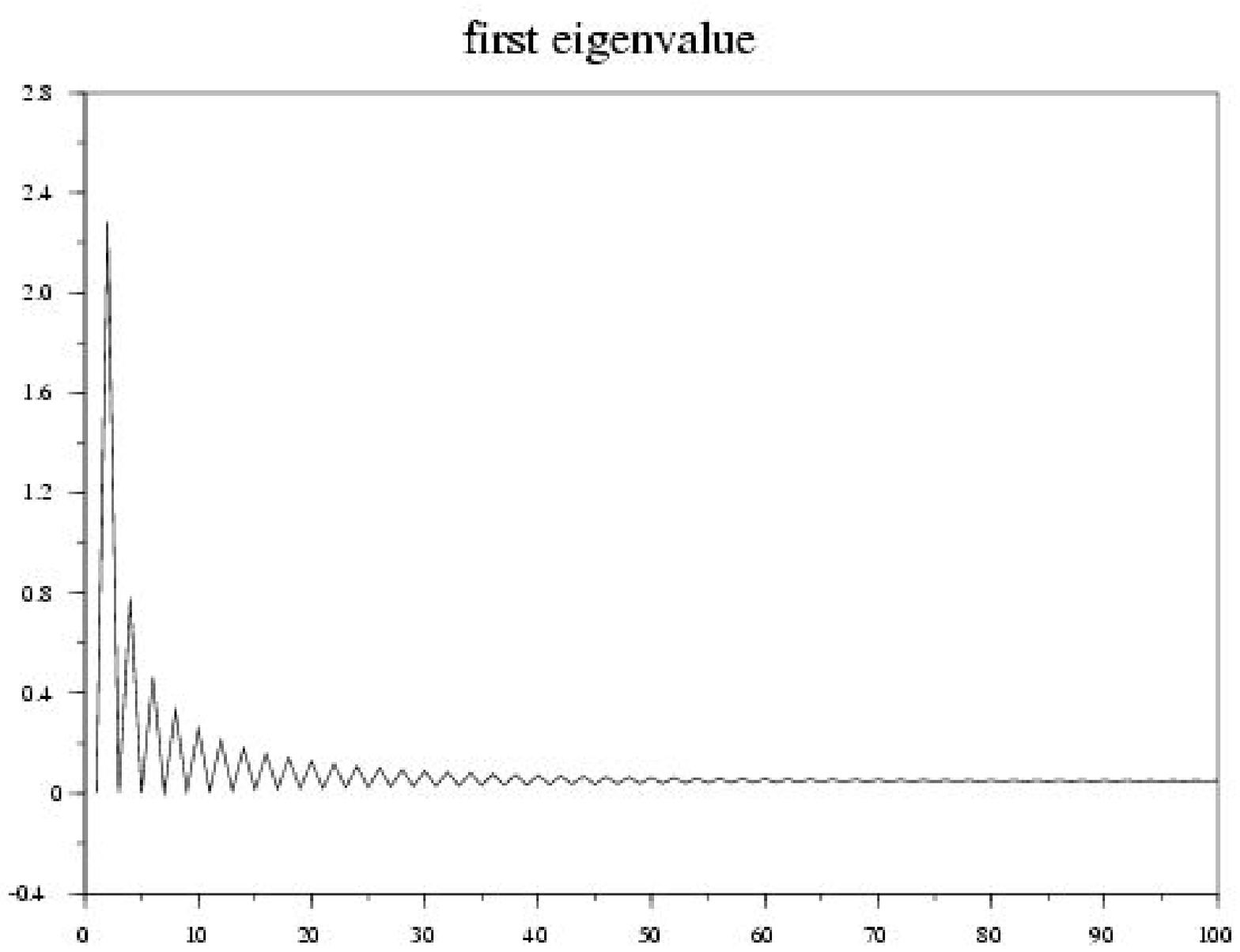 } 
  { \epsfysize=4,0cm     \epsfbox  {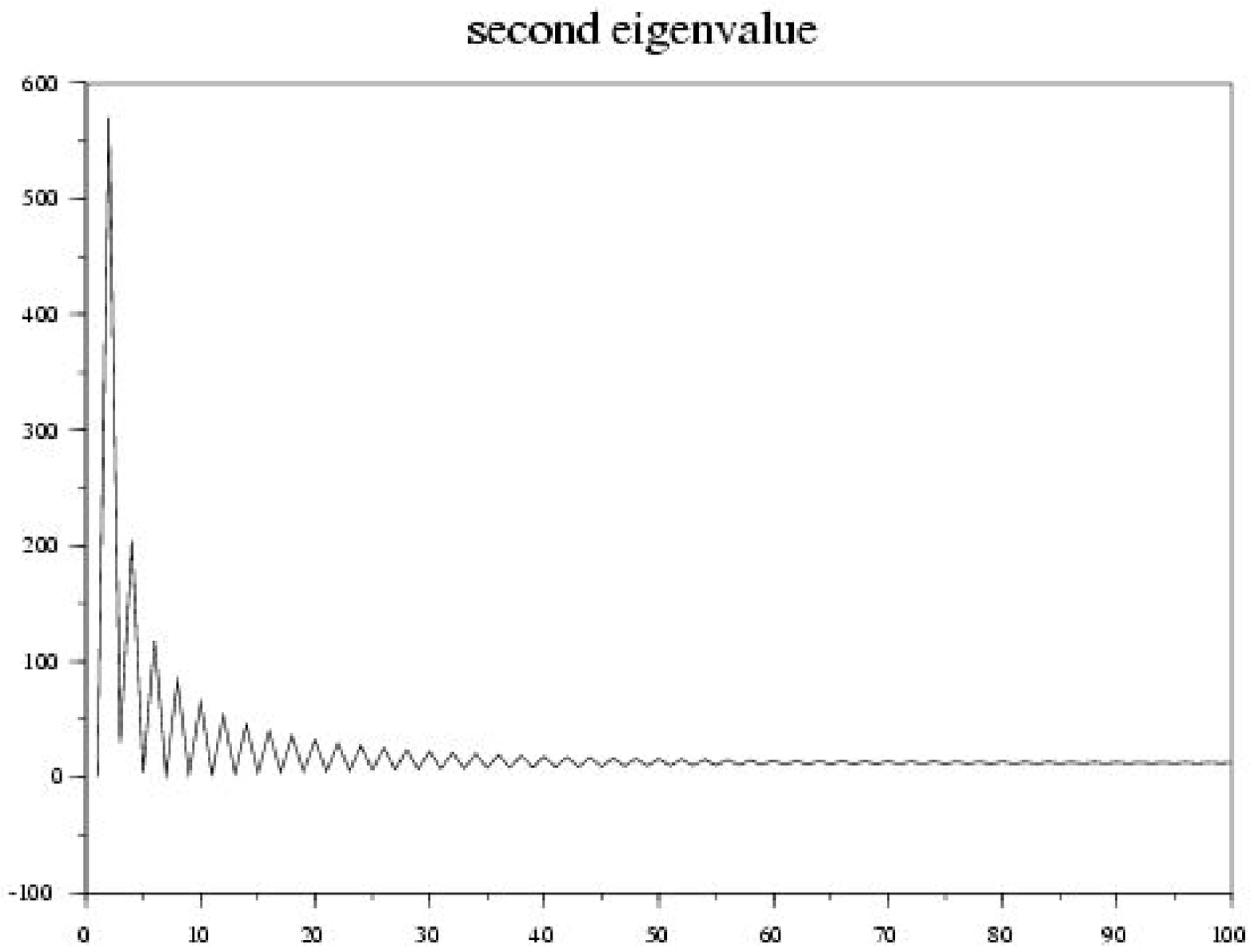 } } 

 \centerline  { { \bf Figures 11 and 12.}  Harmonic oscillator. }
 \centerline  { Two first eigenvalues of numerical solution (
$\mu \,=\, 0.1 , \,\, \alpha \,=\,  0.01 , \,\,  \Delta t \,=\,  100 $). }
%  \smallskip  \smallskip
%%%%%%%%%%%%%%%%%%%%%%%%%%%%%%%%%%%%%%%%%%%%%%%%%%%%%%%%%%%%%%%%%%%%%%%%%%%%%%%%%%%%%%%

 \bigskip  \bigskip   \noindent  {\bf 4.3) $ \quad$ String of high speed vehicles.}
\smallskip \noindent $\bullet \quad$
This example has been considered by  Athans, Levine and  Levis [ALL67] in modelling
position and velocity control for a string of high speed vehicles. Let $\,N \,$ be
some integer and 

\smallskip \noindent  (4.6) $  \qquad \quad  \displaystyle
n \,\,= \,\, 2 \, N \,- \, 1 \,$ 

\smallskip \noindent 
be the order of the given matrices $\, A_N ,\, K_N \,$ and $\, Q_N . \,$ The matrices
$\, A_N \,$ $\,  K_N \,$ and  $\, Q_N \,$  admit the following structure : 

\smallskip \noindent  (4.7) $  \qquad \quad \displaystyle
 A_N \,\, = \,\,   \pmatrix {a_{11}  & a_{12}  & 0 & 0 & \cdots & 0  \cr
0 & a_{22}  & a_{23} &  0 &   \cdots & 0 \cr  
\vdots & \ddots & \ddots & \ddots  & \ddots  & \vdots  \cr  
0 & \cdots &  0 & a\ib{N-2\, N-2} & a\ib{N-2\, N-1} & 0 \cr
0 & \cdots &    & 0 & a\ib{N-1\, N-1} & -1 \cr
0 & \cdots &  \cdots  & 0 & 0 &-1  \cr }\,, \,   $

\smallskip \noindent  (4.8) $  \qquad \quad \displaystyle
K_N \,\, = \,\, {\rm diag} \, \, \bigl( \, 1 \,,\, 0 \,,\, 1 \,,\, 0 \,,\, \cdots 
\,,\,  1 \,,\, 0 \,,\, 1 \,,\, 0  \,,\, \cdots  1 \,,\, 0 \,,\, 1 \, \bigr) \,$

\smallskip \noindent  (4.9) $  \qquad \quad \displaystyle
Q_N \,\, = \,\, {\rm diag} \, \, \bigl( \, 10 \,,\, 0 \,,\, 10 \,,\, 0 \,,\, \cdots 
\,,\,  10 \,,\, 0 \,,\, 10 \,,\, 0  \,,\, \cdots  10 \,,\, 0 \,,\, 10 \, \bigr) \,.\, 
$

\smallskip \noindent 
The unknown positive definite matrix $\, X_N \,$ satisfies the algebraic Riccati
equation 

\smallskip \noindent  (4.10) $  \qquad \quad \displaystyle
X_N \, K_N \, X_N \,\,- \,\, \bigl( A_N^{\displaystyle \rm t} \, X_N \,+\, X_N \,
A_N \, \bigr) \,\,- \,\, Q_N \,\,\, = \,\,\, 0 \,. \,$

\smallskip \noindent
The solution of this equation is detailed on 'figure' 13 with 10  significative
decimals. The first six decimals are absolutly identical to the ones published by
Laub [La79].

%%%%%%%%%%%%%%%%%%%%%%%%%%%%%%%%%%%%%%%%%%%%%%%%%%%%%%%%%%%%%%%%%%%%%%%%%%%%%%%%%%%%%%%
\bigskip \bigskip 
\bigskip \noindent   columns $1$ to $3$ :

\setbox30= \vbox {\halign{#&#&#&#  \cr
$ +1.3630206938E+00 \,\, $ & $ \,\, +2.6172154724E+00 \,\, $ & $ \,\,  -7.0542734123E-01 $ \cr
$ +2.6172154724E+00 \,\, $ & $ \,\, +7.5925521955E+00 \,\, $ & $ \,\, -1.6803557707E+00 $ \cr
$ -7.0542734123E-01 \,\, $ & $ \,\, -1.6803557707E+00 \,\, $ & $ \,\, +1.7747816032E+00 $ \cr
$ +9.3685970173E-01 \,\, $ & $ \,\, +1.4752196951E+00 \,\, $ & $ \,\, +2.1577096313E+00 $ \cr
$ -2.9366643189E-01 \,\, $ & $ \,\, -4.5950584109E-01 \,\, $ & $ \,\, -6.0913599887E-01 $ \cr
$ +4.7735386064E-01 \,\, $ & $ \,\, +6.6514730720E-01 \,\, $ & $ \,\, +6.7071749345E-01 $ \cr
$ -1.9737508953E-01 \,\, $ & $ \,\, -2.6614220828E-01 \,\, $ & $ \,\, -2.6284317351E-01 $ \cr
$ +2.1121165236E-01 \,\, $ & $ \,\, +2.8065373289E-01 \,\, $ & $ \,\, +2.6614220828E-01 $ \cr
$ -1.6655183115E-01 \,\, $ & $ \,\, -2.1121165236E-01 \,\, $ & $ \,\, -1.9737508953E-01 $ \cr }}
\setbox31= \hbox{ $\vcenter {\box30}$ }
\setbox44=\hbox{\noindent  $\displaystyle  \quad    \box31  $}
\smallskip \noindent $ \box44 $

\smallskip \noindent  columns $4$ to $6$ :

\setbox30= \vbox {\halign{#&#&#&#  \cr
$  +9.3685970173E-01   \,\, $ & $ \,\,  -2.9366643189E-01    \,\, $ & $ \,\, +4.7735386064E-01   $ \cr
$  +1.4752196951E+00   \,\, $ & $ \,\,  -4.5950584109E-01    \,\, $ & $ \,\, +6.6514730720E-01  $ \cr
$ +2.1577096313E+00    \,\, $ & $ \,\,  -6.0913599887E-01    \,\, $ & $ \,\, +6.7071749345E-01  $ \cr
$  +8.2576995027E+00   \,\, $ & $ \,\,  -1.9464979789E+00    \,\, $ & $ \,\, +1.7558734280E+00  $ \cr
$  -1.9464979789E+00   \,\, $ & $ \,\,  +1.8056048615E+00    \,\, $ & $ \,\, +1.9464979789E+00  $ \cr
$  +1.7558734280E+00   \,\, $ & $ \,\,  +1.9464979789E+00    \,\, $ & $ \,\, +8.2576995027E+00  $ \cr
$  -6.7071749345E-01  \,\, $ & $ \,\,   -6.0913599887E-01    \,\, $ & $ \,\, -2.1577096313E+00  $ \cr
$  +6.6514730720E-01  \,\, $ & $ \,\,   +4.5950584109E-01    \,\, $ & $ \,\, +1.4752196951E+00   $ \cr
$ -4.7735386064E-01   \,\, $ & $ \,\,   -2.9366643189E-01    \,\, $ & $ \,\, -9.3685970173E-01   $ \cr}}
\setbox31= \hbox{ $\vcenter {\box30}$ }
\setbox44=\hbox{\noindent  $\displaystyle  \quad    \box31  $}
\smallskip \noindent $ \box44 $

 \vfill \eject  % \smallskip 
\noindent  columns $7$ to $9$ :

\setbox30= \vbox {\halign{#&#&#&#  \cr
$  -1.9737508953E-01 \,\, $ & $ \,\, +2.1121165236E-01  \,\, $ & $ \,\, -1.6655183115E-01   $ \cr  
$  -2.6614220828E-01 \,\, $ & $ \,\, +2.8065373289E-01  \,\, $ & $ \,\, -2.1121165236E-01   $ \cr
$  -2.6284317351E-01 \,\, $ & $ \,\, +2.6614220828E-01  \,\, $ & $ \,\, -1.9737508953E-01   $ \cr
$  -6.7071749345E-01 \,\, $ & $ \,\, +6.6514730720E-01  \,\, $ & $ \,\, -4.7735386064E-01   $ \cr
$  -6.0913599887E-01 \,\, $ & $ \,\, +4.5950584109E-01  \,\, $ & $ \,\, -2.9366643189E-01   $ \cr
$  -2.1577096313E+00 \,\, $ & $ \,\, +1.4752196951E+00  \,\, $ & $ \,\, -9.3685970173E-01   $ \cr
$  +1.7747816032E+00 \,\, $ & $ \,\, +1.6803557707E+00  \,\, $ & $ \,\, -7.0542734123E-01   $ \cr
$  +1.6803557707E+00 \,\, $ & $ \,\, +7.5925521955E+00  \,\, $ & $ \,\, -2.6172154724E+00   $ \cr
$  -7.0542734123E-01 \,\, $ & $ \,\, -2.6172154724E+00  \,\, $ & $ \,\, +1.3630206938E+00 $ \cr}}
\setbox31= \hbox{ $\vcenter {\box30}$ }
\setbox44=\hbox{\noindent  $\displaystyle  \quad    \box31  $}
\smallskip \noindent $ \box44 $

 \bigskip \noindent 

 \centerline  { { \bf Figures 13.}  String of high speed vehicles. Matrix X is $9$ by
 $9$.  }

 \centerline  {  Numerical solution of stationary Riccati equation (4.10). }

 \centerline  {   The   parameters $\mu=0.1$ and $\Delta t=0.1$ have been used in
 homographic scheme.  }
 \smallskip  \smallskip

\bigskip  \bigskip \noindent  {\bf 4.4) $ \quad$ Control of the wave equation.}
\smallskip \noindent $\bullet \quad$
The fourth example is the control of the wave equation in one space dimension 

\smallskip \noindent  (4.11) $  \qquad \quad \displaystyle
{{\partial^2 y}\over{\partial t^2}} \,\,-\,\, c^2 \, {{\partial^2 y}\over{\partial
x^2}} \,\,=\,\, \sum_{i=1}^{m} \, \gamma_i (x) \, u_i (t) \,\,,\quad x \in ] \,0
\,,\, L [ \, $ 

\smallskip \noindent
with homogeneous Dirichlet boundary conditions

\smallskip \noindent  (4.12) $  \qquad \quad \displaystyle
y(t,\,0) \,\, = \,\, y(t,\,L) \,\, = \,\, 0 \,. \, $ 

\smallskip \noindent
For solving problem (4.11)-(4.12), we use a spectral decomposition on the eigenmodes
$\, \Phi_j(x)\,$  that are solution of the stationay problem : 

\smallskip \noindent  (4.13) $  \qquad \quad \displaystyle
-\, {{\partial^2 \, \Phi_j(x)}\over{\partial x^2}} \,\,=\,\, \lambda_j\, \,
\Phi_j(x) \,\,, \qquad x \in ] \,0 \,,\, L [ \, $ 
\smallskip \noindent  (4.14) $  \qquad \quad \displaystyle
\Phi_j(0) \,\,= \,\, \Phi_j(L) \,\,= \,\, 0 \,$

\smallskip \noindent 
and are classically explicited by 

\smallskip \noindent  (4.15) $  \qquad \quad \displaystyle
\Phi_j(x) \,\,= \,\, \sqrt{2} \,\,\,  {\rm sin} \, \Bigl( {{j \pi
x}\over{L}} \Bigr) \,\,,\quad x \in ] \,0 \,,\, L [ \, $  
\smallskip \noindent  (4.16) $  \qquad \quad \displaystyle
\lambda_j \,\,= \,\, {{j^2 \, \pi^2}\over{L^2}} \,\,,\quad j \,=\, 1,\, 2 ,\, 3 \cdots
\,$ 

\smallskip \noindent $\bullet \quad$
In practice we restrict to  an approximation with the $\, N \,$ first modes :  
 
\smallskip \noindent  (4.17) $  \qquad \quad \displaystyle
y(t,\,x) \,\,= \,\, \sum_{j=1}^{N} \, y_j(t)  \, \Phi_j(x) \,\,, \quad t \, \geq \, 0
\,\,, \quad x \in ] \,0 \,,\, L [ \,\,, \, $

\smallskip \noindent 
and with such a spectral approximation, problem (4.11)-(4.12) is projected with $\,
L^2 \,$ scalar product $\, < {\scriptstyle \bullet} \,,\, {\scriptstyle \bullet} >
\,$

\smallskip \noindent  (4.18) $  \qquad \quad \displaystyle
< \, u \,,\, v \, > \,\,= \,\, {{1}\over{L}} \int_{0}^{L} \!\! u(t) \, v(t) \, 
{\rm d}t \,$ 

\smallskip \noindent 
and the discrete formulation stands as 

\smallskip \noindent  (4.19) $  \qquad \quad \displaystyle
{{{\rm d}^2 y_j} \over{{\rm d}t^2}} \,\,+\,\, c^2 \, \lambda_j^2 \,\,  y_j(t) 
\,\,\,= \,\,\, \sum_{i=1}^{m} \,  u_i (t) \, < \gamma_i ( {\scriptstyle \bullet})
\,,\, \Phi_j ( {\scriptstyle \bullet}) > \,\,,\quad 1 \leq j \leq N \,. \, $

\smallskip \noindent 
We reduce this discrete differential system to a first order one by setting 

\smallskip \noindent  (4.20) $  \qquad \quad \displaystyle   
Y(t) \,\,= \,\, \Bigl( \, y_1 \,,\, y_2 \,,\, \cdots \,,\, y_N \,,\, 
{{{\rm d} y_1} \over{{\rm d}t}} \,,\,{{{\rm d} y_2} \over{{\rm d}t}} \,,\, \cdots 
\,,\, {{{\rm d} y_N} \over{{\rm d}t}} \, \Bigr)^{\displaystyle \! \rm t} \quad  \in
\R^{2N} \,$
\smallskip \noindent  (4.21) $  \qquad \quad \displaystyle   
U(t) \,\,= \,\, \Bigl( \, u_1 (t) \,,\, u_2 (t) \,,\, \cdots \,,\,  u_m (t)
\,\Bigr)^{\displaystyle \!  \rm t} \quad    \in \R^{n} \,$
\smallskip \noindent  (4.22) $  \qquad \quad \displaystyle   
\Lambda \,\,= \,\,   c^2 \,  {\rm diag}  \, \bigl( \, \lambda_1 \,,\, \cdots \,,\,
\lambda_N \, \bigr) \,$
\smallskip \noindent  (4.23) $  \qquad \quad \displaystyle   
C_{j\, i} \,\,= \,\, < \gamma_i ( {\scriptstyle \bullet}) \,,\, \Phi_j ( {\scriptstyle
\bullet}) > \,\,, \qquad 1 \leq i \leq n \,\,, \quad 1 \leq j \leq N \,$

\smallskip \noindent
and system  (4.19) can be re-written as 

\smallskip \noindent  (4.24) $  \qquad \quad \displaystyle 
{{{\rm d}}\over{{\rm d}t}} \, Y(t) \,\,\,= \,\,\,\, \pmatrix {0 & {\rm I} \cr -\Lambda
& 0 \cr } \, Y(t) \,\,+\,\, \pmatrix {0  \cr C \cr } \, U(t) \,. \,$ 

\smallskip \noindent  
The matrices $\,R\,$ and $\, Q \,$ associated to the definition of the cost function
$\, J({\scriptstyle \bullet}) \,$ (relation (1.6) are of order $ \,n \,$ and $\, 2N \,$
 respectively. We have chosen the follo\-wing simple form parameterized by $\,
\alpha = 1 \,$ and $\, \beta = 10 \,$ in this particular test case : 

\smallskip \noindent  (4.25) $  \qquad \quad \displaystyle   
R \,\,= \,\,   {\rm diag}  \, \bigl( \, \alpha \,,\, \cdots \,,\, \alpha \, \bigr) \,$
\smallskip \noindent  (4.26) $  \qquad \quad \displaystyle   
Q \,\,= \,\,   {\rm diag}  \, \bigl( \, \beta \,,\, \cdots \,,\, \beta \, \bigr)
\,. \,$

\smallskip \noindent $\bullet \quad$
We have tested the scheme for matrices of order $\, n = 2N=10 \,$ and for time step $\,
\Delta t \,=\, 0.01 \,$ and a parameter $\, \mu \,=\, 0.001 .\,$ We represent on
Figures 14 to 23 the ten different eigenvalues of Riccati matrix $\, X(t)\,$
and for these parameters the scheme is stable as for as the explicit two stages 
Runge-Kutta scheme is unstable. 
                       
%%%%%%%%%%%%%%%%%%%%%%%%%%%%%%%%%%%%%%%%%%%%%%%%%%%%%%%%%%%%%%%%%%%%%%%%%%%%%%%%%%%%%%% 
%  \bigskip 
\smallskip \noindent 
  { \epsfysize=4,0cm    \epsfbox   {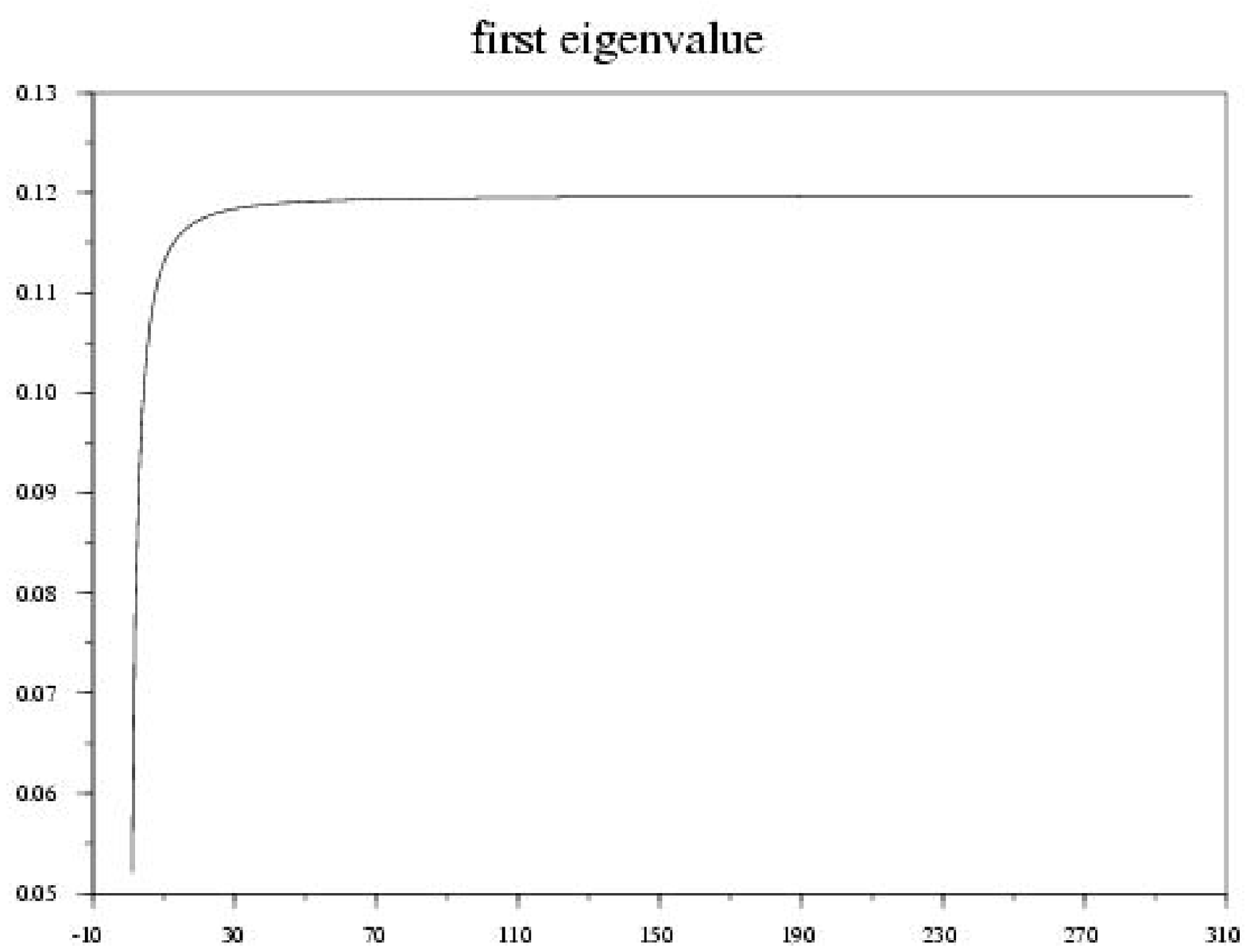 } 
  { \epsfysize=4,0cm     \epsfbox  {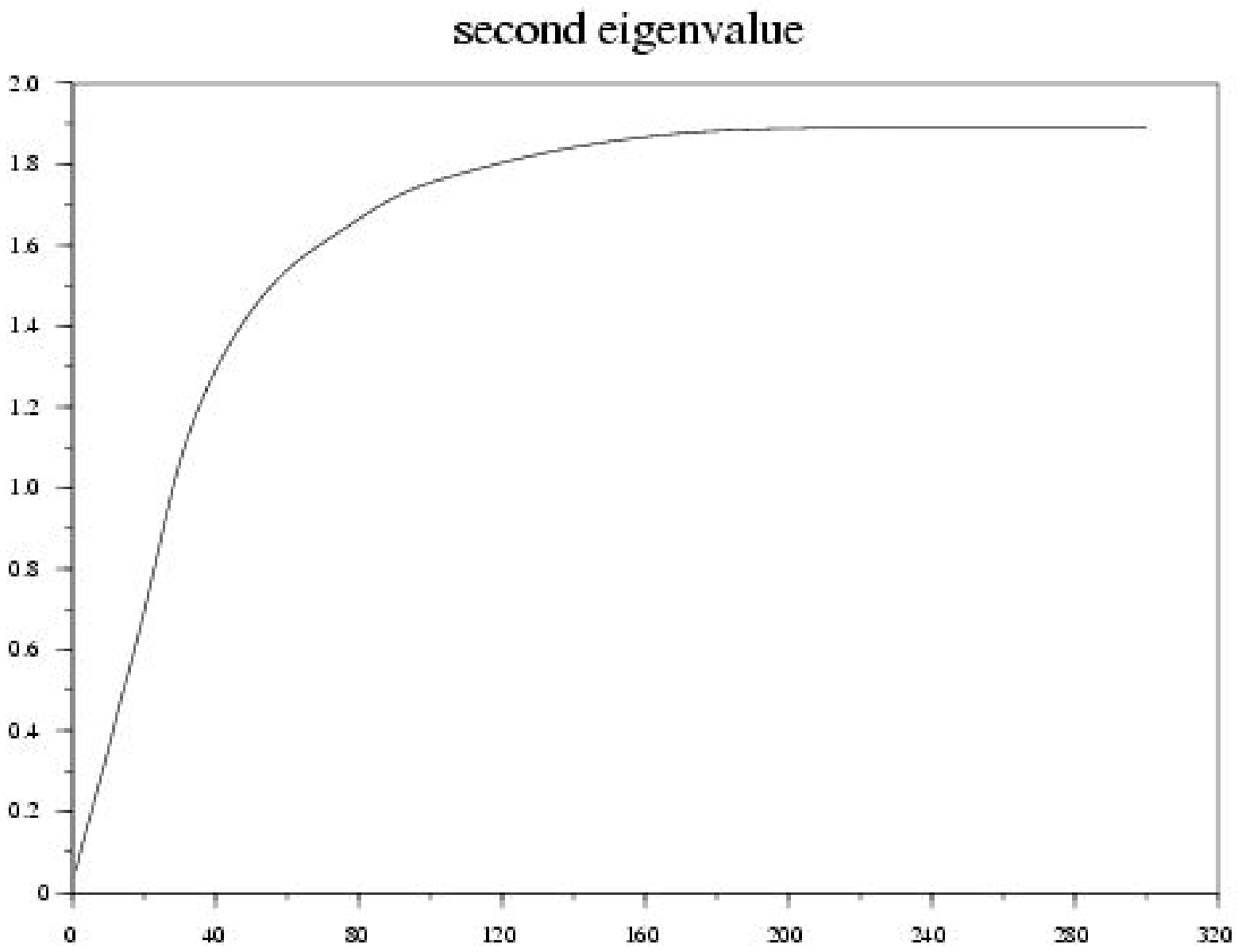 } } 

 \centerline  { { \bf Figures 14 and 15.}  Wave control test.  }
 \centerline  { Two first eigenvalues of numerical solution.}
%%%%%%%%%%%%%%%%%%%%%%%%%%%%%%%%%%%%%%%%%%%%%%%%%%%%%%%%%%%%%%%%%%%%%%%%%%%%%%%%%%%%%%%

%%%%%%%%%%%%%%%%%%%%%%%%%%%%%%%%%%%%%%%%%%%%%%%%%%%%%%%%%%%%%%%%%%%%%%%%%%%%%%%%%%%%%%%
%  \bigskip 
\smallskip \noindent  
  { \epsfysize=4,0cm    \epsfbox   {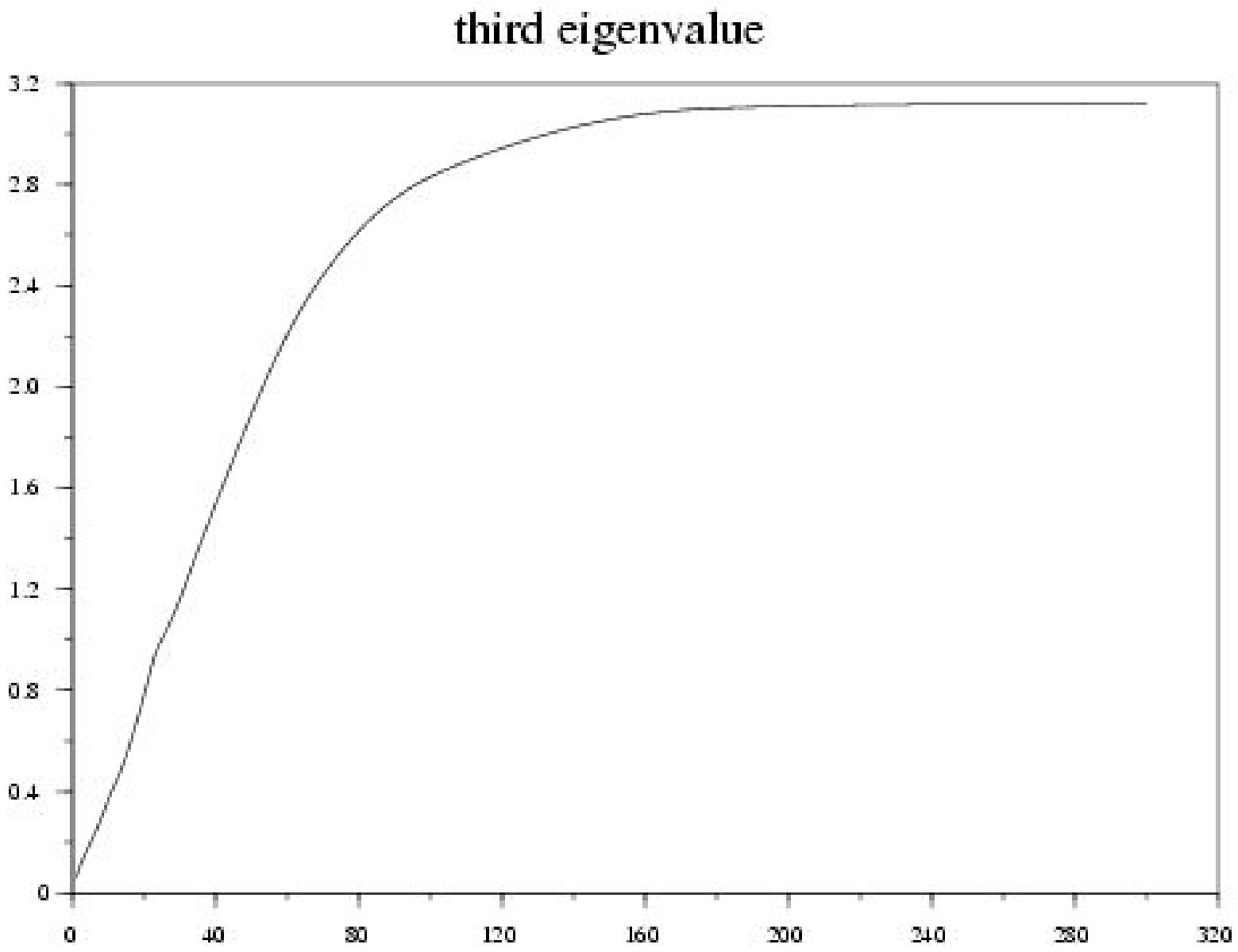 } 
  { \epsfysize=4,0cm     \epsfbox  {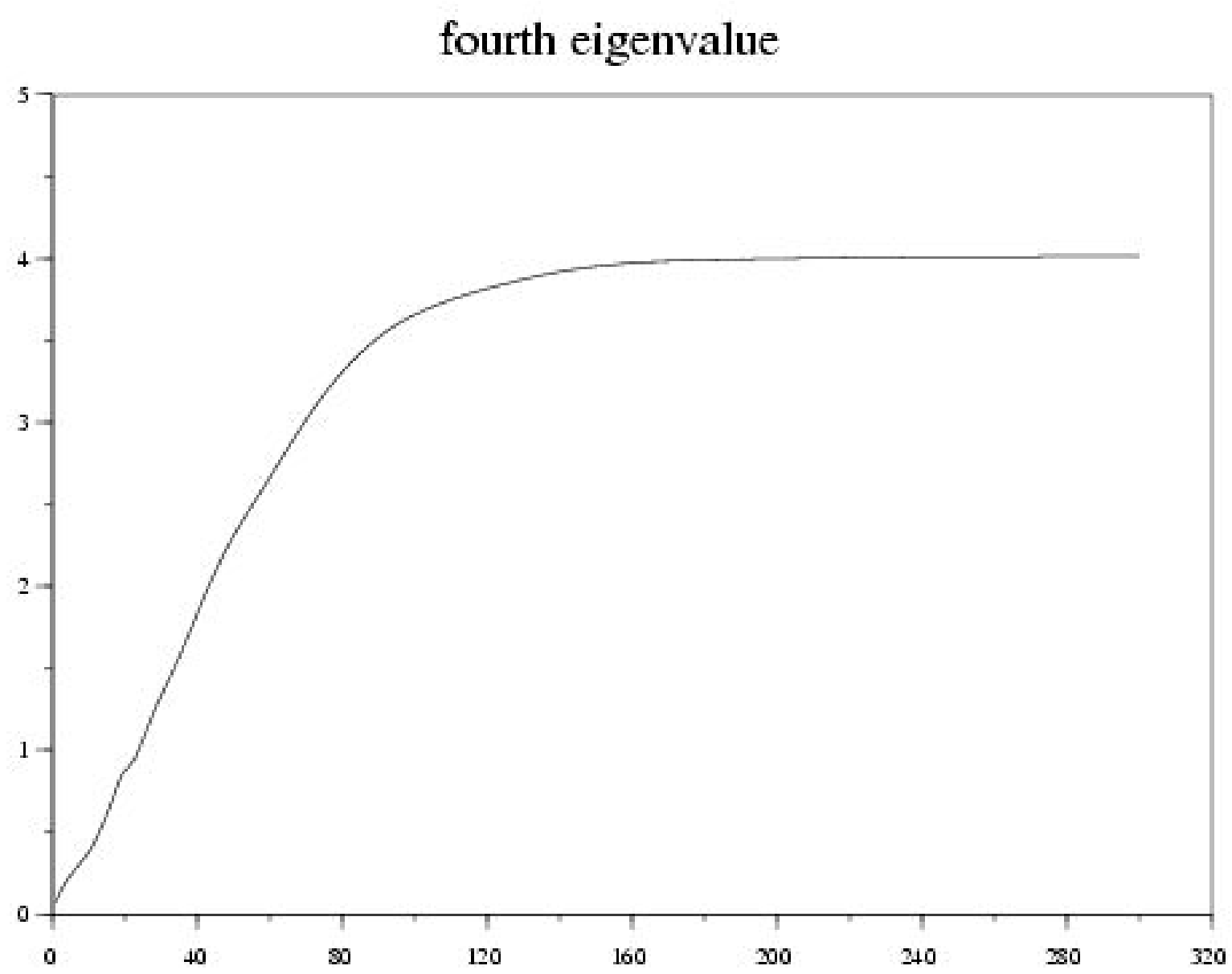 } } 

 \centerline  { { \bf Figures 16 and 17.}  Wave control test. }
 \centerline  { Third and fourth eigenvalues of numerical solution. }
%%%%%%%%%%%%%%%%%%%%%%%%%%%%%%%%%%%%%%%%%%%%%%%%%%%%%%%%%%%%%%%%%%%%%%%%%%%%%%%%%%%%%%%

%%%%%%%%%%%%%%%%%%%%%%%%%%%%%%%%%%%%%%%%%%%%%%%%%%%%%%%%%%%%%%%%%%%%%%%%%%%%%%%%%%%%%%%
%  \bigskip 
\smallskip \noindent  
  { \epsfysize=4,0cm    \epsfbox   {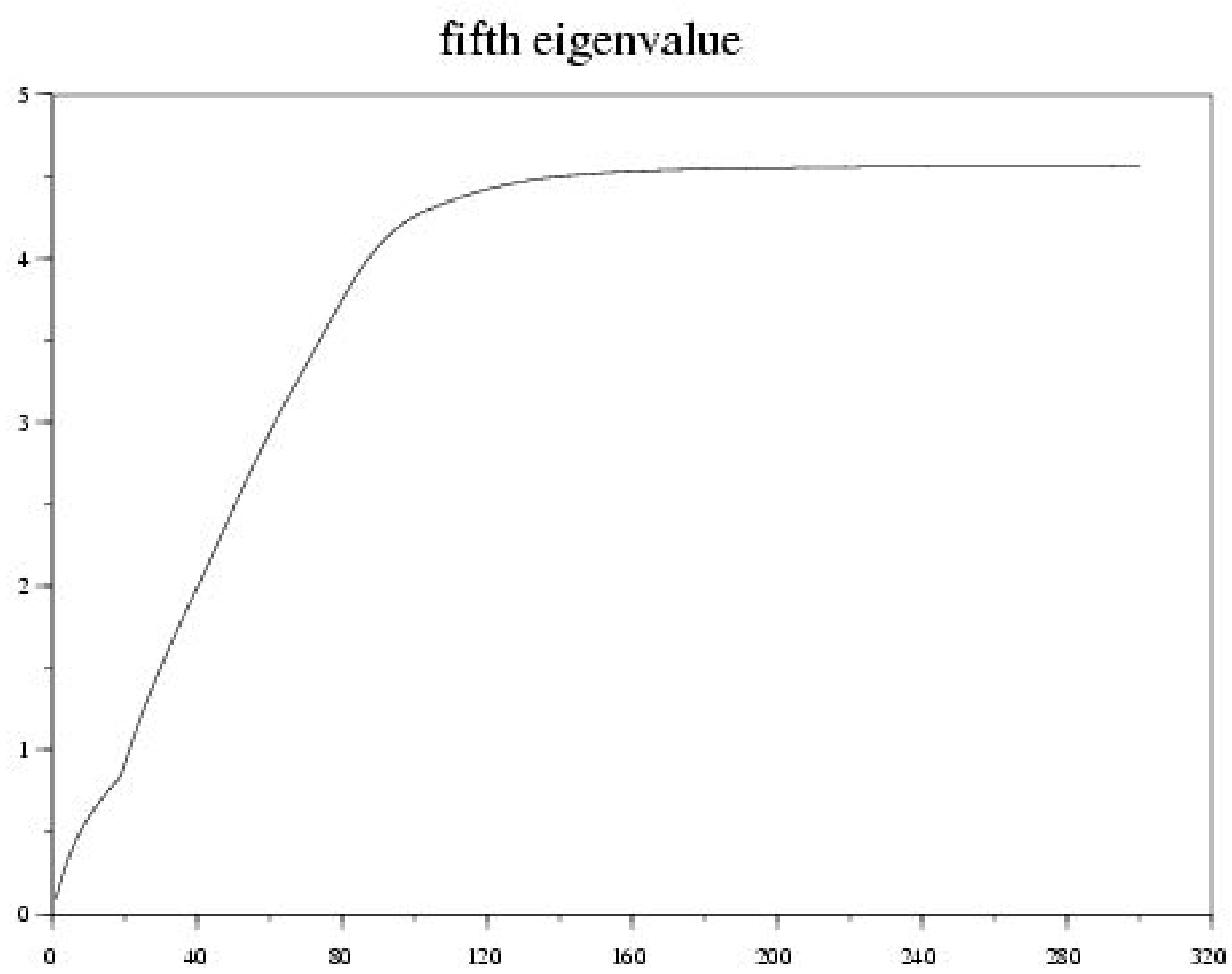 } 
  { \epsfysize=4,0cm     \epsfbox  {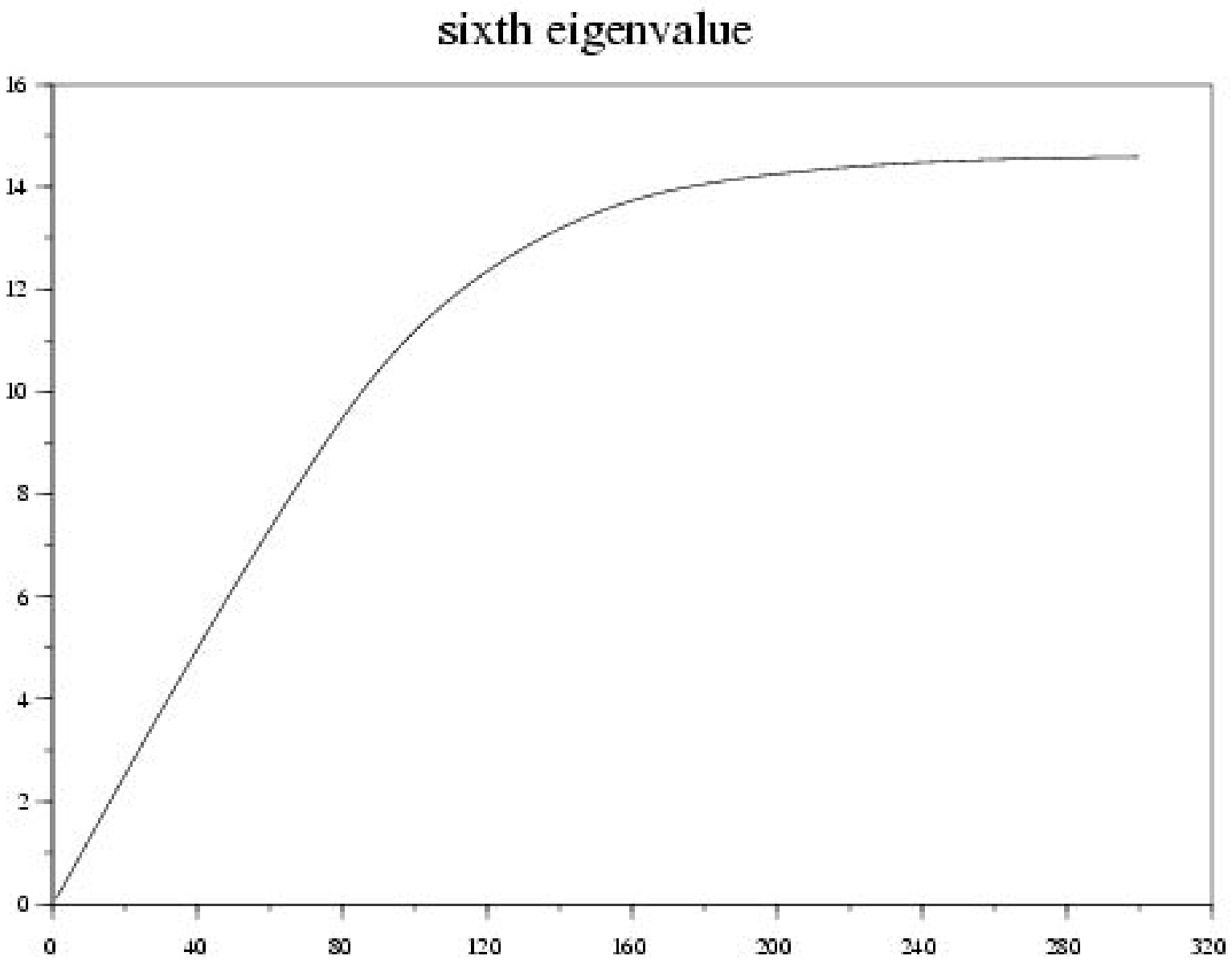 } } 

 \centerline  { { \bf Figures 18 and 19.}  Wave control test.  
Fifth and sixth eigenvalues. }}
%%%%%%%%%%%%%%%%%%%%%%%%%%%%%%%%%%%%%%%%%%%%%%%%%%%%%%%%%%%%%%%%%%%%%%%%%%%%%%%%%%%%%%%

%%%%%%%%%%%%%%%%%%%%%%%%%%%%%%%%%%%%%%%%%%%%%%%%%%%%%%%%%%%%%%%%%%%%%%%%%%%%%%%%%%%%%%%
 \bigskip \noindent 
  { \epsfysize=4,0cm    \epsfbox   {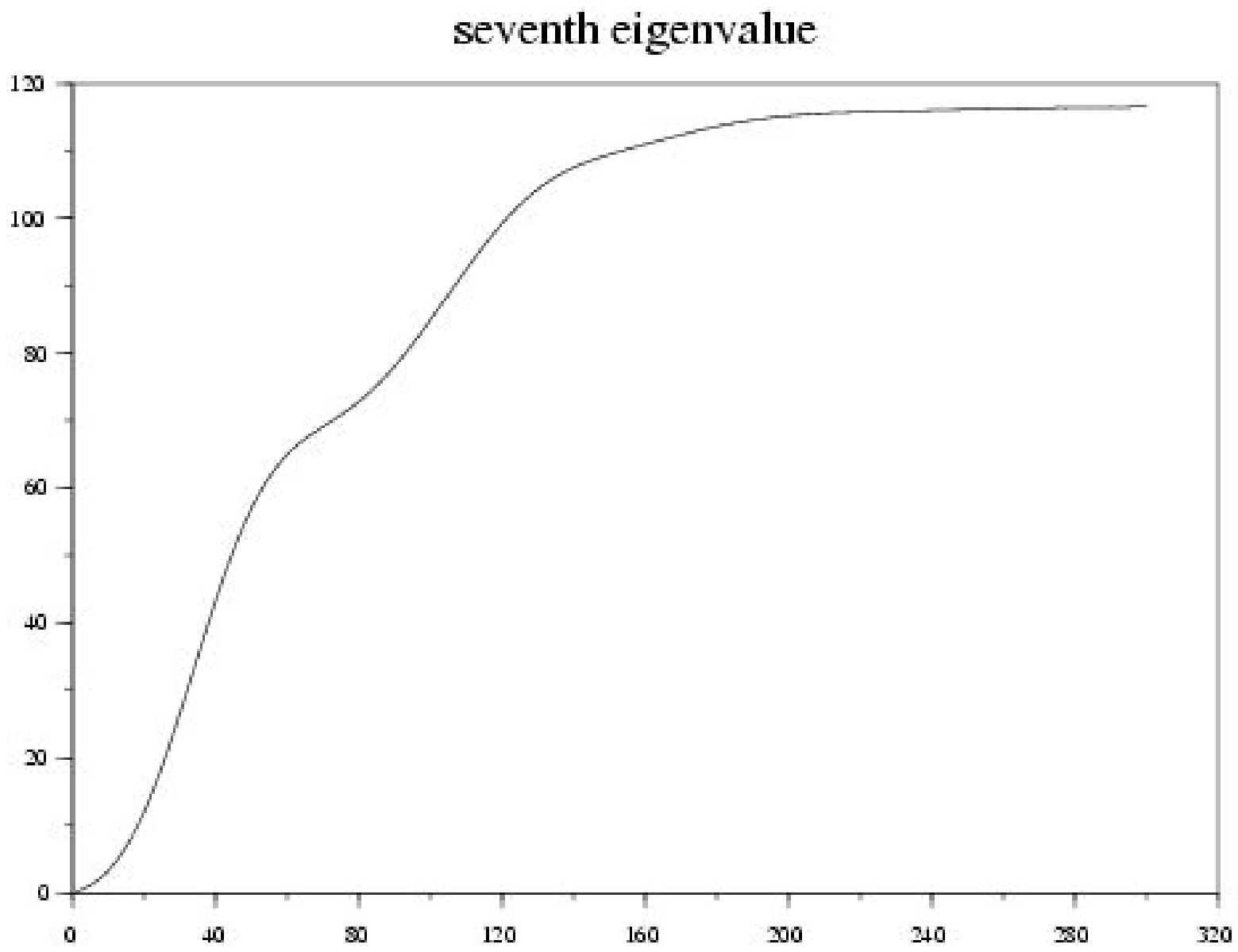 } 
  { \epsfysize=4,0cm     \epsfbox  {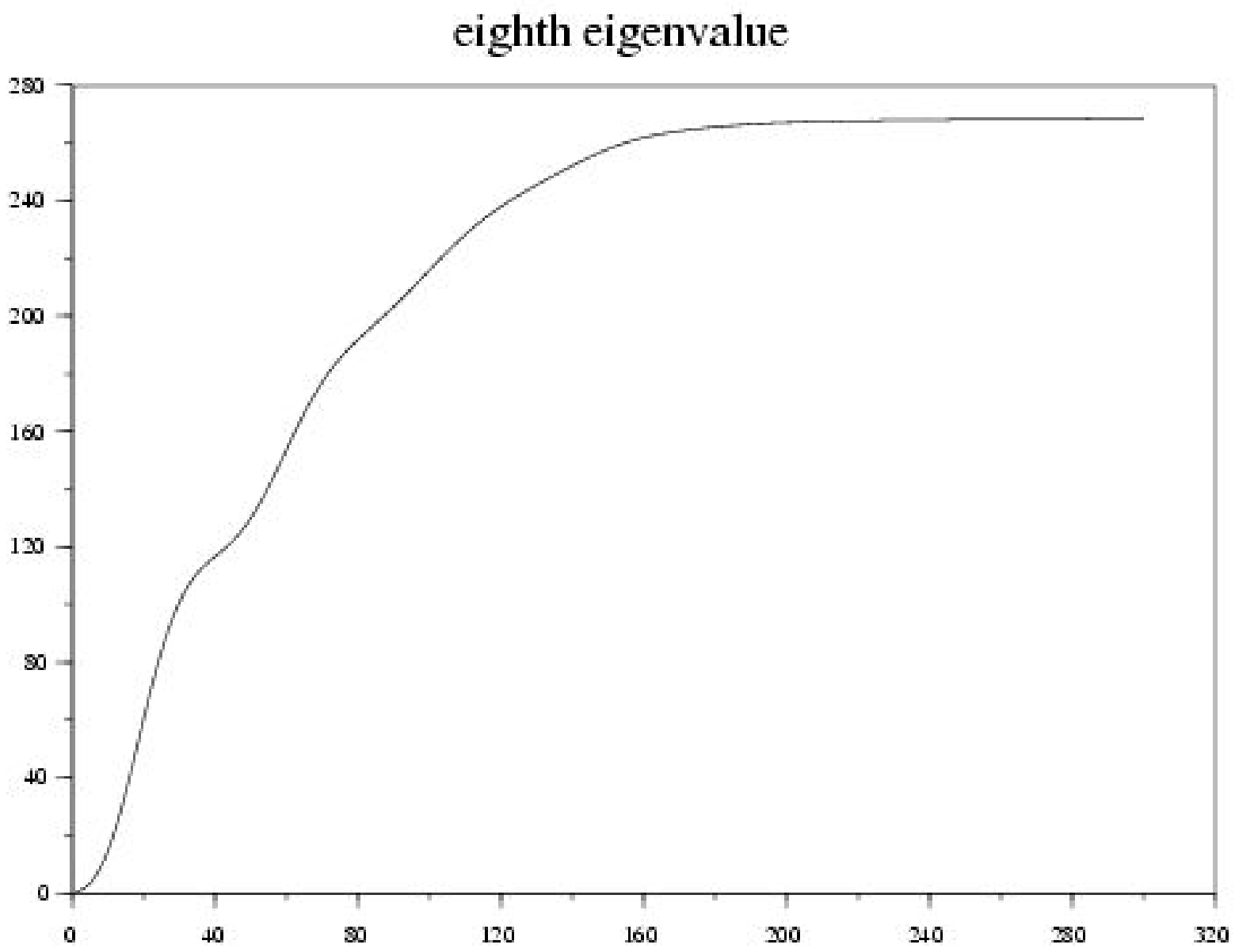 } } 

 \centerline  { { \bf Figures 20 and 21.}  Wave control test.  
Seventh and eighth eigenvalues. }}
%%%%%%%%%%%%%%%%%%%%%%%%%%%%%%%%%%%%%%%%%%%%%%%%%%%%%%%%%%%%%%%%%%%%%%%%%%%%%%%%%%%%%%%

%%%%%%%%%%%%%%%%%%%%%%%%%%%%%%%%%%%%%%%%%%%%%%%%%%%%%%%%%%%%%%%%%%%%%%%%%%%%%%%%%%%%%%%
 \bigskip \noindent 
  { \epsfysize=4,0cm    \epsfbox   {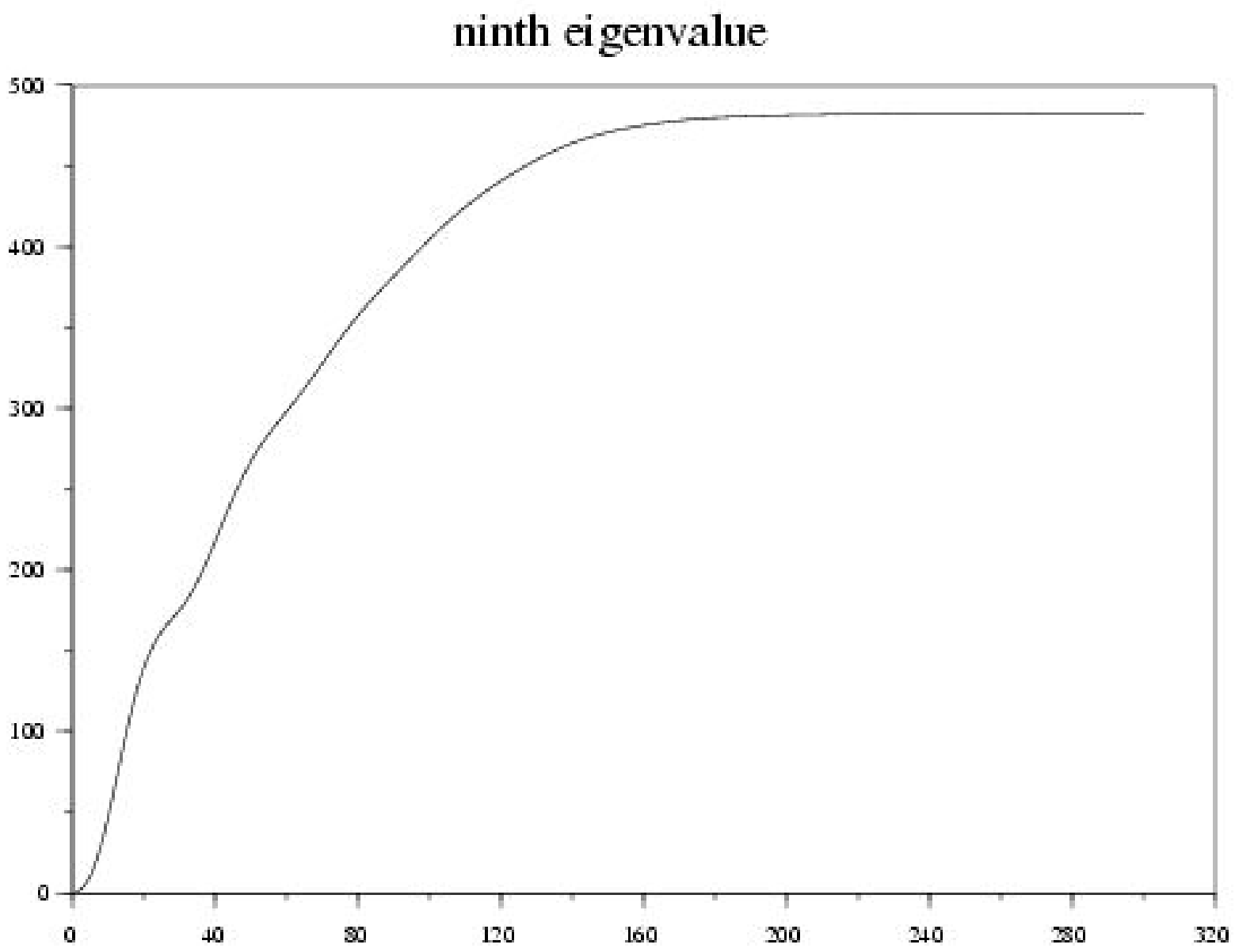 } 
  { \epsfysize=4,0cm     \epsfbox  {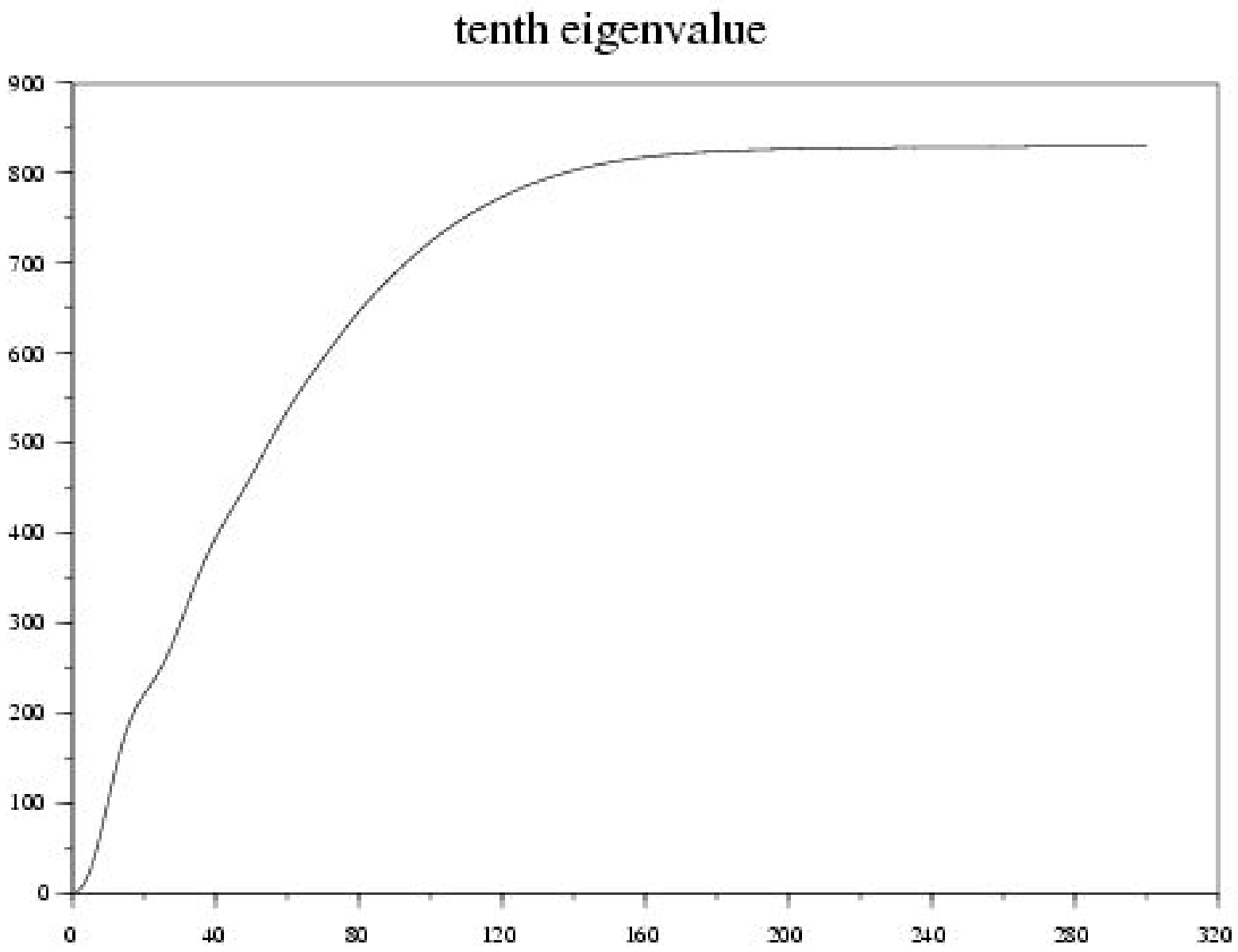 } } 

 \centerline  { { \bf Figures 22 and 23.}  Wave control test.  
Ninth and tenth eigenvalues. }}
%%%%%%%%%%%%%%%%%%%%%%%%%%%%%%%%%%%%%%%%%%%%%%%%%%%%%%%%%%%%%%%%%%%%%%%%%%%%%%%%%%%%%%%
   
\bigskip \bigskip   \noindent  {\smcaps  5)  $ \quad$  Conclusion.}         
\smallskip \noindent $\bullet \quad$
We have proposed an harmonic   scheme for the resolution of the matrix Riccati
equation. The scheme is implicit, unconditionnaly stable, needs to use one scalar 
parameter and to solve   a linear system of equations for each time step. This scheme
is convergent in the scalar case. In the matrix case, harmonic scheme  has good
monotonicity properties and discrete solution tends to the positive solution of
algebraic  Ricatti equation as discrete time tends to infinity. We have computed first 
 test cases of matrix square root, harmonic ocsillator, string of vehicles and 
discretized wave equation  where classical explicit schemes fail to give a
definite positive discrete solution. Our first numerical experiments show stability
and robustness when various parameters have large variations. We plan to develop this
work in two directions~: first prove the convergence of the harmonic scheme in the
case of Ricatti matrix equation and second construct a multistep version in order to
achieve second order accuracy.

\bigskip \bigskip   \noindent  {\smcaps  6)  $ \quad$    Acknowledgments.}  
\smallskip \noindent $\bullet \quad$
The authors thank Marius Tucsnak for helpfull comments on a preliminary draft 
 of this report.

\bigskip \bigskip  
\noindent  {\smcaps  7)  $ \quad$  References.} 

\smallskip \hangindent=9mm \hangafter=1 \noindent  
 [AF66]      M. Athans, P.L. Faulb. {\it Optimal Control.
An Introduction to the theory and Its Applications}, Mc Graw-Hill, New York, 1966.

\smallskip \hangindent=9mm \hangafter=1 \noindent  
 [ALL67]        M. Athans, W.S. Levine and A. Levis. A
system for the optimal and suboptimal position and velocity control for a string of
high speed of vehicles, in {\it Proc. 5th Int. Analogue Compytation Meeting}, Lausanne,
Switzerland, september 1967.

\smallskip \hangindent=9mm \hangafter=1 \noindent  
  [Ba91]        R. Baraille. D\'eveloppement de
sch\'{e}mas num\'{e}riques adapt\'{e}s \`a l'hydro\-dyna\-mique, {\it Th\`{e}se de
l'Universit\'{e} Bordeaux 1},  d\'{e}cembre 1991.

\smallskip \hangindent=9mm \hangafter=1 \noindent  
  [Ca79]       D. Cariolle. Mod\`{e}le
unidimentionnel de
chimie de l'ozone, Internal note, {\it Etablissement d'Etudes et de Recherches
M\'{e}t\'{e}orologiques}, Paris 1979.

\smallskip \hangindent=9mm \hangafter=1 \noindent  
 [DE96]        L. Dieci, T. Eirola. Preserving 
monotonicity in  the numerical solution of Riccati differential equations, 
{\it Numer. Math.,} vol$.\,$ 74, p. 35-47, 1996.

\smallskip \hangindent=9mm \hangafter=1 \noindent  
 [Du93]       F. Dubois. Un sch\'{e}ma implicite 
non lin\'{e}airement inconditionellement stable pour l'\'{e}quation de Riccati,
{\it unpublished manuscript}, april 1993.

\smallskip \hangindent=9mm \hangafter=1 \noindent  
 [DS95]      F. Dubois, A. Sa\"{\i}di.  Un sch\'{e}ma 
implicite non lin\'{e}airement inconditionellement stable pour l'\'{e}quation
de Riccati, {\it Congr\`{e}s d'Analyse Num\'{e}\-rique}, Super Besse, France, may 1995.

\smallskip \hangindent=9mm \hangafter=1 \noindent  
 [DS2k]       F. Dubois, A. Sa\"{\i}di.  Unconditionally
stable scheme for Riccati equation,   {\it in}   Control of systems governed by partial
differential equations, F.~Conrad and  M.~Tucsnak Editors,
  {\it  ESAIM: Proceedings}, vol. 8,  p. 39-52, 2000.

\smallskip \hangindent=9mm \hangafter=1 \noindent  
 [FR84]      P. Faurre, M. Robin.  {\it El\'ements
d'automatique}, Dunod, Paris, 1984.

\smallskip \hangindent=9mm \hangafter=1 \noindent  
 [KS72]       H. Kawakernaak, R. Sivan. {\it Linear optimal
control systems}, Wiley, 1972.

\smallskip \hangindent=9mm \hangafter=1 \noindent  
 [La74]        J.D. Lambert. {\it Computational methods in
ordinary differential equations}, J.~Wiley \& Sons, 1974.

\smallskip \hangindent=9mm \hangafter=1 \noindent  
  [La79]     A.J. Laub. A Schur Method for Solving
Algebraic Riccati Equations, {\it IEEE Trans. Aut. Control},  
vol$.\,$AC-24, p.~913-921, 1979.

\smallskip \hangindent=9mm \hangafter=1 \noindent  
 [Le86]     F.L. Lewis. {\it Optimal Control},
J. Wiley-Interscience, New York, 1986.

\smallskip \hangindent=9mm \hangafter=1 \noindent
  [Li68]     J.L. Lions. {\it Contr\^{o}le
optimal des syst\`emes gouvern\'es par des \'equations aux d\'eriv\'ees
partielles}, Dunod, Paris, 1968.

\smallskip \hangindent=9mm \hangafter=1 \noindent
 [LT86]      P. Lascaux, R. Th\'eodor.   {\it
Analyse num\'erique matricielle appliqu\'ee \`a l'art de l'ing\'enieur}, Masson,
Paris, 1986.

\smallskip \hangindent=9mm \hangafter=1 \noindent
  [Mi84]    J.C. Miellou. Existence globale
pour une classe de syst\`{e}mes paraboliques semi-lin\'{e}aires mod\'{e}lisant le
probl\`{e}me de la stratosph\`{e}re : la m\'{e}thode de la  fonction  agr\'{e}g\'{e}e,
{\it C.R.  Acad. Sci., Paris, Serie I}, t. 299, p.723-726, 1984.

\smallskip \hangindent=9mm \hangafter=1 \noindent
  [Sa97]     A. Sa\"{\i}di. Analyse math\'{e}matique et 
num\'{e}rique de mod\`{e}les de structures intelligentes
et de leur contr\^{o}le, {\it Th\`{e}se de l'Universit\'{e} Pierre et Marie Curie},
France, april 1997.

\bye